# Boolean and ortho fuzzy subset logics

Daniel J. Greenhoe

*Abstract*: Constructing a fuzzy subset logic **L** with Boolean properties is notoriously difficult because under a handful of "reasonable" conditions, we have the following three debilitating constraints: (1) Bellman and Giertz in 1973 showed that if **L** is *distributive*, then it must be *idempotent*. (2) Dubois and Padre in 1980 showed that if **L** has the *excluded middle* or the *non-contradiction* property or both, then it must be *non-idempotent*. (3) Bellman and Giertz also demonstrated in 1973 that even if **L** is *idempotent*, then the *only* choice available for the $(\wedge, \vee)$ logic operator pair is the $(\min, \max)$ operator pair. Thus it would seem impossible to construct a non-trivial fuzzy subset logic with *Boolean* properties. However, this paper examines these three results in detail, and shows that "hidden" in the hypotheses of the three is the assumption that the operator pair $(\wedge, \vee)$ is *pointwise evaluated*. It is further demonstrated that removing this constraint yields the following results: (A) It is indeed possible to construct *fuzzy subset logic*s that have all the *Boolean* properties, including that of *idempotency*, *non-contradiction*, *excluded middle*, and *distributivity*. (B) Even if *idempotency* holds, $(\min, \max)$ is *not* the only choice for $(\wedge, \vee)$.

## Contents

# Introduction

**The problem.**    This paper addresses a well known conflict in *fuzzy subset logic* theory. *Fuzzy subset logic* is the foundation for *fuzzy set theory*, just as *classical logic* is the foundation for *classical set theory*. Just as *classical logic* and the *algebras of sets* constructed upon it are *Boolean algebra*s (and hence have all the "nice" *Boolean* properties such as *idempotency, excluded middle, non-contradiction, distributivity*, etc.), one might very much prefer to have a *fuzzy subset logic* and resulting *fuzzy subset algebras* that are *Boolean* as well.[1] But it has been found that under what would seem to be very "reasonable" conditions, this is simply not possible. In particular, we have the following crippling constraints:

- 📖 Bellman and Giertz in 1973 demonstrated that under very "reasonable" conditions, if we want a *fuzzy subset logic* that is *distributive*, then it also must be *idempotent*.[2]
- 📖 Dubois and Padre in 1980 demonstrated that under very "reasonable" conditions, if we want a *fuzzy subset logic* that has the *non-contradiction* and *excluded middle* properties, then that logic is *not idempotent*…and therefore not only fails to be a *Boolean algebra*, but also is not even a *lattice*.[3]

---

[1] *excluded middle*: $x \vee \neg x = 1$. *non-contradiction*: $x \wedge \neg x = 0$. *idempotency*: $x \vee x = x$ and $x \wedge x = x$. *distributivity*: $x \vee (y \wedge z) = (x \vee y) \wedge (x \vee z)$ and $x \wedge (y \vee z) = (x \wedge y) \vee (x \wedge z)$. *classic Boolean properties*: Theorem A.42 page 30.

[2] see *fuzzy operators idempotency theorem* (Theorem 1.25 page 12)

[3] *Dubois-Padre 1980 result*: see *fuzzy negation idempotency theorem* (Theorem 1.28 page 14) and *Dubois-Padre 1980 theorem* (Corollary 1.29 page 15). Every *lattice* is a *Boolean algebra*, but not conversely (Definition A.11 page 24, Definition A.41 page 30). A *lattice* is *idempotent*, *commutative*, *associative*, and *absorptive* (Theorem A.14 page 25). A *Boolean algebra* has all these properties but is moreover *bounded*, *distributive*, *complemented*, *de Morgan*, *involutory*, and has *identity* (Theorem A.42 page 30).

Moreover, even if we are willing to give up the *non-contradiction* and *excluded middle* properties and retain *idempotency*, Bellman and Giertz also demonstrated in 1973 that the *only* choice we have for the logic operator pair $(\wedge, \vee)$ is the $(\min, \max)$ operator pair such that $(\wedge, \vee) = (\min, \max)$.[4]

**A solution.**   Section 1 of this paper examines these results in detail, and demonstrates that "hidden" in the hypotheses of these results is the assumption that the operator pair $(\wedge, \vee)$ is *pointwise evaluated*.[5] Section 2 demonstrates that if this constraint is removed, then it is indeed possible to construct *fuzzy subset logic*s that have all the *Boolean* properties, including that of *idempotency*, *non-contradiction*, *excluded middle*, and *distributivity*.

**A solution yielding ortho fuzzy subset logics.**   In this paper, a *logic $L'$* $\triangleq$ $(\ X,\ \vee,\ \wedge,\ \neg,\ 0,\ 1\ ;\ \leq,\ \rightarrow)$ is defined as a *lattice $L$* $\triangleq$ $(\ X,\ \vee,\ \wedge\ ;\ \leq)$ with a *negation function* $\neg$ and *implication function* $\rightarrow$ defined on this lattice. And in this sense, the logic $L'$ is said to be "constructed on" the lattice $L$. This paper demonstrates that it is possible to construct *fuzzy subset logics* on *Boolean lattices* yielding *Boolean fuzzy subset logics*. However, more generally, it is also demonstrated that it is possible to construct *fuzzy subset logics* on *orthocomplemented lattices* yielding *ortho fuzzy subset logics*. The main difference between a *Boolean lattice* and an *orthocomplemented lattice* is that the latter does not in general support *distributivity*.[6] On finite sets, there are significantly more choices of *orthocomplemented lattices* than there are *Boolean lattices*.[7]And so having the option of constructing *ortho fuzzy subset logics* is arguably not without advantage. The disadvantage is that we give up the guarantee of *distributivity*. But some authors[8] have investigated structures without this property anyways. In fact, one could argue that the "crucial" properties that we would really like a logic to have, if possible, are the following:

(1).  *disjunctive idempotence*:   $x \vee x \quad = \quad x$   and
(2).  *conjunctive idempotence*:   $x \wedge x \quad = \quad x$   and
(3).  *excluded middle*:   $x \vee \neg x \quad = \quad 1$   and
(4).  *non-contradiction*:   $x \wedge \neg x \quad = \quad 0$   .

Not all *fuzzy logics* have all the these properties. Of course all *Boolean logics* have them. But more generally than *Boolean logics* and less generally than *fuzzy logics*, all *ortho logics* have them as well.[9]

---

[4]*Bellman-Giertz 1973 result*: see *fuzzy min-max theorem* (Theorem 1.26 page 13) and *Bellman-Giertz 1973 theorem* (Corollary 1.27 page 14). $(\wedge, \vee)$ in an *ordered set*: Definition A.9 page 24 and Definition A.8 page 24. $(\wedge, \vee)$ in a *lattice*: Definition A.11 page 24. $(\wedge, \vee)$ in a *logic*: Definition C.5 page 50. $(\min, \max)$: Definition 1.15 page 7.

[5]*pointwise evaluated*: (Definition 1.12 page 7)

[6] *logic*: Definition C.5 page 50. *lattice*: Definition A.11 page 24. *negation function*: Definition B.2 page 35. *implication function*: Definition C.1 page 45 *Boolean lattice*: Definition A.41 page 30. *orthocomplemented lattice*: Definition A.44 page 31. *ortho negation*: Definition B.3 page 35. *ortho+distributivity=Boolean*: Proposition A.50 page 33

[7]There are a total of 5 *orthocomplemented lattices* with 8 elements; of these 5, only 1 is *Boolean*. There are a total of 10 *orthocomplemented lattices* with 8 elements or less; of these 10, only 4 are *Boolean*. For further details, see Example A.46 page 31.

[8] 📖 [Alsina et al.(1980)Alsina, Trillas, and Valverde], 📖 [Hamacher(1976)] ⟨referenced by 📖 [Alsina et al.(1983)Alsina, Trillas, and Valverde]⟩

[9] properties of *fuzzy negations* and hence also *fuzzy logics*: Theorem B.11 page 36. properties of *ortho negations* and hence also *ortho logics*: Theorem B.15 page 37. relationships between logics: Figure 13 page 50.

**Negation functions.**    There are several types of *negation functions* and information about these functions is scattered about in the literature.  APPENDIX B introduces several types of negation, describes some of their properties, and shows where *fuzzy negation*, *ortho negation*, and *Boolean negation* "fit" into the larger structure of *negations* in general.

**Implication functions.**    Defining an *implication* function for a logic constructed on a *Boolean lattice* is straightforward because we can simply use the *classical implication* $x \stackrel{c}{\to} y \triangleq \neg x \vee y$. However, defining an *implication* function for a *non-Boolean* logic is more difficult.  APPENDIX C addresses the problem of defining implication functions on *lattices*, including lattices that are *non-Boolean*.

# 1   Fuzzy subset operators

A *fuzzy subset* is often specified in terms of a *membership function*.  A *fuzzy subset logic* is a *lattice* of *membership function*s together with a *fuzzy negation* function and an *implication* function. Although its definition is simple and straightforward, *fuzzy subset logic* has some notorious problems attempting to provide some very standard Boolean properties.[10]

## 1.1   Indicator functions

In *classical subset theory*, a subset $A$ of a set $X$ can be specified using an *indicator function* $\mathbb{1}_A(x)$ (next definition).  An indicator function specifies concretely whether or not an element is a member of $A$. That is, it is a convenient "indicator" of whether or not a particular element is in a subset. A subset that can be defined using an indicator function is a *crisp subset* (next definition).

**Definition 1.1**   [11]   Let $2^X$ be the *power set* of a set $X$. Let $Y^X$ be the *set of all functions* mapping from $X$ to a set $Y$. The **indicator function** $\mathbb{1}_A \in \{0,1\}^X$ is defined as
$$\mathbb{1}_A(x) = \left\{ \begin{array}{ll} 1 & \text{if } x \in A \\ 0 & \text{if } x \notin A \end{array} \right\} \quad \forall x \in X, \, A \in 2^X.$$
The parameter $A$ of $\mathbb{1}_A$ is a **crisp subset** of $X$ if $\mathbb{1}_A(x)$ is an *indicator function* on $X$.

Every set $X$ has at least one crisp subset (itself). A set of subsets, together with the relation $\subseteq$, form an *ordered set*, and in some cases also form a *lattice*.  Common set structures include the *power set* $2^X$, *topologies*, *rings of sets* and *algebras of sets*. A set structure may be represented in terms of subsets, or equivalently, in terms of *set indicator functions*.[12]

---

[10] *fuzzy subset*: Definition 1.7 page 6, *fuzzy subset logic*: Definition 1.11 page 6, *membership function*: Definition 1.7 page 6, *lattice*: Definition A.11 page 24, *fuzzy negation*: Definition B.2 page 35, *implication*: Definition C.1 page 45 and Definition C.5 page 50; problems: Theorem 1.26 page 13 and Theorem 1.28 page 14.

[11] [Feller(1971)], page 104 ⟨1. Baire Functions⟩, ✎ [Aliprantis and Burkinshaw(1998)], page 126, ✎ [Hausdorff(1937)], page 22, ✎ [de la Vallée-Poussin(1915)] page 440

[12] *ordered set*: Definition A.1 page 22, *lattice* (Definition A.11 page 24), *set indicator function*: Definition 1.1 page 4, *topologies*: Example 1.3 page 5 and Example 1.4 page 5; examples of *set structures*: Example 1.3 page 5 and Example 1.5 page 5.

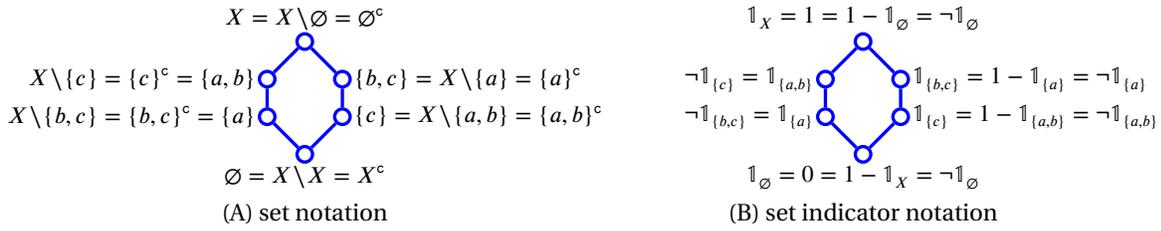

Figure 1: set structures on $O_6$ (see Example 1.5 page 5)

**Remark 1.2** Often set structures are defined in terms of set operators like *intersection* $\cap$, *union* $\cup$, and *set complement* c. The set operators $(\cap, \cup, c, \Rightarrow, \varnothing, X)$ in turn can be defined in terms of arithmetic operators $\big(\min, \max, f(x) \triangleq 1 - x, g(x,y) \triangleq y - xy, 0, 1\big)$ on the *set indicator* function[13] or in terms of classic logic operators $(\wedge, \vee, \neg, \to, 0, 1)$ like this:

$$
\begin{aligned}
0 &\triangleq \mathbb{1}_\varnothing &&= 0 \\
1 &\triangleq \mathbb{1}_X &&= 1 \\
\mathbb{1}_A \vee \mathbb{1}_B &\triangleq \mathbb{1}_{A \cup B} &&= \max(\mathbb{1}_A, \mathbb{1}_B) \\
\mathbb{1}_A \wedge \mathbb{1}_B &\triangleq \mathbb{1}_{A \cap B} &&= \min(\mathbb{1}_A, \mathbb{1}_B) \\
\neg \mathbb{1}_A &\triangleq \mathbb{1}_{A^c} &&= 1 - \mathbb{1}_A \\
\mathbb{1}_A \to \mathbb{1}_B &\triangleq \mathbb{1}_{A \Rightarrow B} &&= \max(1 - \mathbb{1}_A, \mathbb{1}_B)
\end{aligned}
$$

where $A \Rightarrow B \triangleq A^c \cup B$ is the *set implication* from $A$ to $B$.[14]

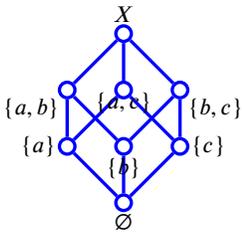

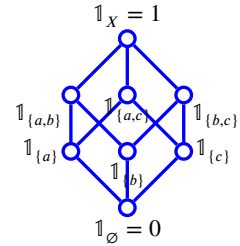

**Example 1.3** The set structures illustrated to the left and right are the **power set** of the set $X \triangleq \{a, b, c\}$. A power set is a special case of an *algebra of sets* and also a *topology*. The lattice to the left uses set notation; the one to the right uses set indicators.

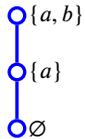

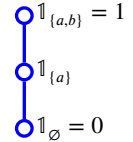

**Example 1.4** The set structures illustrated to the left and right are a **topology** on the set $X \triangleq \{a, b\}$. The lattice to the left uses set notation; the one to the right uses set indicators.

**Example 1.5** The set structures illustrated in Figure 1 (page 5) are not *topologies* (or *algebras of sets* or *power sets*), but are *set structures* none the less. The negation function in the structure is an *ortho negation* (Definition B.3 page 35). The lattice in (A) uses set notation; the one in (B) uses set indicators.

**Definition 1.6** Let $\mathbb{1}^X$ be the *set of all indicator functions* on a set $X$. Let a *logic* be defined as in Definition C.5 (page 50). A **crisp subset logic** is a *logic*
$$\big(\mathbb{1}^X, \vee, \wedge, \neg, \mathbb{1}_\varnothing, \mathbb{1}_X; \leqslant, \to\big).$$

---

[13] 📖 [Aliprantis and Burkinshaw(1998)], page 126, 📖 [Hausdorff(1937)], pages 22–23

[14] 📄 [Ellerman(2010)] ⟨§1.7; $A \Rightarrow B = (A^c \cup B)^*$ where $C^*$ is the *interior* of a set $C$ in a *topological space*⟩

## 1.2   Membership functions

In a crisp subset $A$ of a crisp set $X$ $(A \subseteq X)$, an element $x \in X$ has only two possible "degrees of membership" in $A$: Either $x$ *is* in $A$ or $x$ is *not* in $A$. Said another way, either $x$ has "full membership" in $A$, or $x$ has "absolute non-membership" in $A$. And this "degree of membership" is specified by an *indicator function* (Definition 1.1 page 4) $\mathbb{1}_A(x)$ which maps from $X$ to the 2-valued set $\{0, 1\}$, where 0 represents "absolute non-membership" and 1 represents "full membership".

In a fuzzy subset $B$ of a crisp set $X$ $(B \subseteq X)$, an element $x \in X$ has a range of possible degrees of membership in $B$. And this membership is specified by a *membership function* (next definition) $\mathfrak{m}_B(x)$ which maps from $X$ to the infinte set $[0 : 1]$.

**Definition 1.7** [15] Let $[0 : 1]$ be the *closed interval* on $\mathbb{R}$ such that $[0 : 1] \triangleq \{x \in \mathbb{R} \, | 0 \le x \le 1\}$. Let $X$ be a set. A function $\mathfrak{m}_A(x)$ is a **membership function** on $X$ if $\mathfrak{m}_A \in [0 : 1]^X$. The parameter $A$ is called a **fuzzy subset** of $X$. For any value $x \in X$, $\mathfrak{m}_A(x) \in [0 : 1]$ represents the "degree of membership" of $x$ in $A$. The condition $\mathfrak{m}_A(x) = 1$ indicates that $x$ has "full membership" in $A$, and the condition $\mathfrak{m}_A(x) = 0$ indicates that $x$ has "absolute non-membership" in $A$.

**Remark 1.8** [16] What is typically called a "fuzzy set" arguably should more accurately be called a "fuzzy subset" because every element $x$ at any "degree of membership" in a fuzzy subset $A$ has absolute full membership in some universal crisp set $X$. And thus $A$ is a subset of the crisp set $X$ $(A \subseteq X)$.

**Remark 1.9** In a crisp set $X$, a *fuzzy subset* $A \subseteq X$ should not be confused with a *random subset* $B \subseteq X$. In the fuzzy subset $A$, an element $x \in X$ has a "degree of membership" in $A$ that specifies "to what extent" $x$ can be considered a member of $A$. In the random subset $B$, the element $x \in X$ has a "degree of likelihood" that $x$ is in $B$ and that specifies the probability that $x$ is a member of $B$. Alternatively, a fuzzy subset is a result of "inference under vagueness", while a random subset is a result of "inference under randomness".[17]

**Example 1.10** Let $A$ be the set of all people who are "young" with *membership function* $\mathfrak{m}_A(x)$. Let $B$ be the set of all people who are "middle age" with *membership function* $\mathfrak{m}_B(x)$. Let $C$ be the set of all people who are "old" with *membership function* $\mathfrak{m}_C(x)$. Of course all these are vague, or "fuzzy", concepts; but the following figure illustrates what the *membership function*s (Definition 1.7 page 6) for these sets **might** look like.

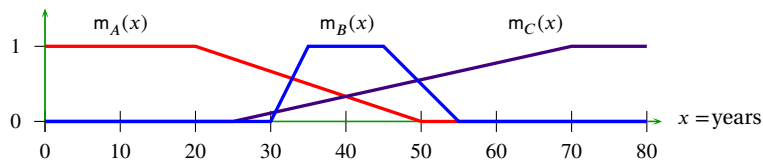

**Definition 1.11** Let $\mathbb{M}$ be a set of *membership function*s (Definition 1.7 page 6).
The structure $\boldsymbol{L} \triangleq (\mathbb{M}, \vee, \wedge, \neg, 0, 1 ; \le)$ is a **fuzzy subset logic** if $\boldsymbol{L}$ is a *fuzzy logic* (Definition C.5 page 50).

---

[15] [Hájek(2011)], page 68 ⟨"absolutely true", "absolutely false⟩, [Dubois(1980)] page 10, [Dubois et al.(2000)Dubois, Ostasiewicz, and Padre], page 42 ⟨"full membership", "absolute non-membership"⟩, [Zadeh(1965)] page 339 ⟨"grade of membership"⟩
[16] [Dubois(1980)] page 10 ⟨Remarks 1⟩, [Kaufmann(1975)]
[17] [Hájek(2011)], page 67 ⟨5.1 Introduction⟩

## 1.3   Operators on membership functions

The *meet-join* operator pair $(\wedge, \vee)$ on a set of indicator functions $\mathbb{1}^X$ induces an ordering relation on $\mathbb{1}^X$.[18] So the operator pairs $(\wedge, \vee)$ can be defined on sets of membership functions to form lattices. But while lattices of set indicators effectively have just one choice for $(\wedge, \vee)$, membership function lattices have many choices.

In this paper, the operators $(\wedge, \vee)$ are called *pointwise evaluated* if at each single value $x$, the functions $[\mathbb{m} \wedge \mathbb{n}](x)$ and $[\mathbb{m} \vee \mathbb{n}](x)$ depend only on the values of $\mathbb{m}(x)$ ($\mathbb{m}$ evaluated at the single value $x$) and $\mathbb{n}(x)$ (next definition).

**Definition 1.12**  [19] Let $L \triangleq (\mathbb{M}, \wedge, \vee)$, where $\mathbb{M}$ is a set of *membership function*s (Definition 1.7 page 6) with operators $(\wedge, \vee)$.  $L$ is **pointwise evaluated**, or said to have **pointwise evaluation**, if there exists $\mathsf{f}, \mathsf{g} \in [0:1]^{[0:1]^2}$ such that

  1. $[\mathbb{m} \wedge \mathbb{n}](x) = \mathsf{f}[\mathbb{m}(x), \mathbb{n}(x)]$  ∀x∈ℝ, and ∀m, n∈𝕄   and
  2. $[\mathbb{m} \vee \mathbb{n}](x) = \mathsf{g}[\mathbb{m}(x), \mathbb{n}(x)]$  ∀x∈ℝ, and ∀m, n∈𝕄

**Example 1.13**

  1. The function $\underline{\wedge}$ defined as  $[\mathbb{m} \,\underline{\wedge}\, \mathbb{n}](x) \triangleq \mathbb{m}(x) + \mathbb{n}(x)$     is *pointwise evaluated* .
  2. The function $\underline{\wedge}$ defined as  $\left[\mathbb{m} \,\underline{\wedge}\, \mathbb{n}\right](x) \triangleq \underbrace{\int_{-\infty}^{x} \mathbb{m}(u)\mathbb{n}(x-u)\,\mathrm{d}u}_{\text{"convolution"}}$  is *not pointwise evaluated*.

**Example 1.14**  Examples of operators that *are* pointwise evaluated include the *min-max operators* (next definition), the *product and probabilistic sum operators* (Definition 1.17 page 8), and the *Łukasiewicz t-norm and t-conorm* (Definition 1.18 page 9).

One of the most common fuzzy logic operator pairs is the *min-max operator pair* (next). As will be demonstrated by the *fuzzy min-max theorem* (Theorem 1.26 page 13), under fairly "reasonable" conditions the *min-max operators* are the *only* choice available for a *fuzzy subset logic*.

**Definition 1.15**  [20] Let $\mathbb{M}$ be a *set of membership functions* on a set $X$. Let $\mathsf{f}(x)$ and $\mathsf{g}(x)$ be functions both with *domain* $X$. Let $\min(\mathsf{f}(x), \mathsf{g}(x))$ and $\max(\mathsf{f}(x), \mathsf{g}(x))$ be the *pointwise minimum* and *pointwise maximum*, respectively, of $\mathsf{f}(x)$ and $\mathsf{g}(x)$ over $X$.  The **min-max operators** $(\wedge, \vee)$ for $L$ are defined as

$$\begin{aligned}
\left[\mathbb{m}_A \vee \mathbb{m}_B\right](x) &\triangleq \max\left[\mathbb{m}_A(x), \mathbb{m}_B(x)\right] \quad \forall m∈𝕄, \ x∈X \\
\left[\mathbb{m}_A \wedge \mathbb{m}_B\right](x) &\triangleq \min\left[\mathbb{m}_A(x), \mathbb{m}_B(x)\right] \quad \forall m∈𝕄, \ x∈X
\end{aligned}$$

**Proposition 1.16**  *Let* $\mathbb{M}$, max, *and* min *defined as in Definition 1.15. Let* $L \triangleq (\mathbb{M}, \vee, \wedge\,; \leq)$ *be an algebraic structure with* $x \leq y \iff x \wedge y = x$.

  $(\wedge, \vee) = (\min, \max) \implies$  *$L$ is a* LATTICE *(Definition A.11 page 24)*.

---

[18] see Remark 1.2 page 5, Proposition A.10 page 24, and Example 1.3 page 5–Example 1.5 page 5
[19] [Dubois(1980)] page 11 ⟨B.a(i)⟩
[20] [Fodor and Yager(2000)], page 133, [Zadeh(1965)] pages 340–341 ⟨(3),(5)⟩; *pointwise ordering*: Definition A.7 page 23

✎PROOF:   To be a lattice, **L** must be *commutative*, *associative*, and *absorptive* (Theorem A.18 page 25).

$$\mathfrak{m} \vee \mathfrak{n} = \max(\mathfrak{m}, \mathfrak{n})$$     by left hypothesis

$$= \max(\mathfrak{n}, \mathfrak{m})$$     by definition of min

$$= \mathfrak{n} \vee \mathfrak{m}$$     by left hypothesis

$$\implies \vee \text{ is } commutative$$

$$\mathfrak{m} \wedge \mathfrak{n} = \min(\mathfrak{m}, \mathfrak{n})$$     by left hypothesis

$$= \min(\mathfrak{n}, \mathfrak{m})$$     by definition of min

$$= \mathfrak{n} \wedge \mathfrak{m}$$     by left hypothesis

$$\implies \wedge \text{ is } commutative$$

$$\mathfrak{m} \vee (\mathfrak{n} \vee \mathfrak{p}) = \max[\mathfrak{m}, \max(\mathfrak{n}, \mathfrak{p})]$$     by left hypothesis

$$= \max[\max(\mathfrak{m}, \mathfrak{n}), \mathfrak{p}]$$     by definition of max

$$= (\mathfrak{m} \vee \mathfrak{n}) \vee \mathfrak{p}$$     by left hypothesis

$$\implies \vee \text{ is } associative$$

$$\mathfrak{m} \wedge (\mathfrak{n} \wedge \mathfrak{p}) = \min[\mathfrak{m}, \min(\mathfrak{n}, \mathfrak{p})]$$     by left hypothesis

$$= \min[\min(\mathfrak{m}, \mathfrak{n}), \mathfrak{p}]$$     by definition of min

$$= (\mathfrak{m} \wedge \mathfrak{n}) \wedge \mathfrak{p}$$     by left hypothesis

$$\implies \wedge \text{ is } associative$$

$$\mathfrak{m} \vee (\mathfrak{m} \wedge \mathfrak{n}) = \max[\mathfrak{m}, \min(\mathfrak{m}, \mathfrak{n})]$$     by left hypothesis

$$= \begin{cases} \max(\mathfrak{m}, \mathfrak{m}) & \text{if } \mathfrak{m}(x) \le \mathfrak{n}(x) \quad \forall x \in X \\ \max(\mathfrak{m}, \mathfrak{n}) & \text{otherwise} \end{cases}$$

$$= \begin{cases} \mathfrak{m} & \text{if } \mathfrak{m}(x) \le \mathfrak{n}(x) \quad \forall x \in X \\ \mathfrak{m} & \text{otherwise} \end{cases}$$

$$= \mathfrak{m}$$

$$\mathfrak{m} \wedge (\mathfrak{m} \vee \mathfrak{n}) = \min[\mathfrak{m}, \max(\mathfrak{m}, \mathfrak{n})]$$     by left hypothesis

$$= \begin{cases} \min(\mathfrak{m}, \mathfrak{n}) & \text{if } \mathfrak{m}(x) \le \mathfrak{n}(x) \quad \forall x \in X \\ \min(\mathfrak{m}, \mathfrak{m}) & \text{otherwise} \end{cases}$$

$$= \begin{cases} \mathfrak{m} & \text{if } \mathfrak{m}(x) \le \mathfrak{n}(x) \quad \forall x \in X \\ \mathfrak{m} & \text{otherwise} \end{cases}$$

$$= \mathfrak{m}$$

$$\implies (\wedge, \vee) \text{ is } absorptive$$

🖝

**Definition 1.17** [21] Let $\mathbb{M}$ be defined as in Definition 1.15. Then for all $\mathfrak{m} \in \mathbb{M}$,
the **probabilistic sum operator** $\vee$ on $\mathbb{M}$ is defined as $\quad [\mathfrak{m}_A \vee \mathfrak{m}_B](x) \triangleq \mathfrak{m}_A(x) + \mathfrak{m}_B(x) - \mathfrak{m}_A(x)\mathfrak{m}_B(x)$
and the **product sum operator** $\wedge$ on $\mathbb{M}$ is defined as $\quad [\mathfrak{m}_A \wedge \mathfrak{m}_B](x) \triangleq \mathfrak{m}_A(x)\mathfrak{m}_B(x)$

Note that the *product and probabilistic sum operators* (previous definition) do *not* in general form a lattice because, for example, they are not in general *idempotent* (a necessary condition for being a lattice—Theorem A.14 page 25). Suppose for example $\mathfrak{m}(p) = \frac{1}{2}$ at some point $p$. Then at that point $p$

$$\mathfrak{m} \vee \mathfrak{m} \triangleq \mathfrak{m} + \mathfrak{m} - \mathfrak{m}\mathfrak{m} = \tfrac{1}{2} + \tfrac{1}{2} - \tfrac{1}{2} \cdot \tfrac{1}{2} = \tfrac{3}{4} \ne \mathfrak{m} \implies \vee \text{ is } non\text{-}idempotent$$

$$\mathfrak{m} \wedge \mathfrak{m} \triangleq \mathfrak{m}\mathfrak{m} = \tfrac{1}{2} \cdot \tfrac{1}{2} = \tfrac{1}{4} \ne \mathfrak{m} \implies \wedge \text{ is } non\text{-}idempotent$$

---

[21] 🕮 [Fodor and Yager(2000)], page 133

**Definition 1.18** [22] Let $L$, $D$, min and max be defined as in Definition 1.15. Then for all $m \in \mathbb{M}$,

the **Łukasiewicz t-conorm** $\vee$ is defined as $\quad [m_A \wedge m_B](x) \triangleq \max\big[0, m_A(x) + m_B(x) - 1\big] \quad$ $\forall m \in \mathbb{M},\, x \in X$ and the **Łukasiewicz t-norm** $\wedge$ is defined as $\quad [m_A \vee m_B](x) \triangleq \min\big[1, m_A(x) + m_B(x)\big] \quad$ $\forall m \in \mathbb{M},\, x \in X$

The *Łukasiewicz t-conorm* is also called the **bold sum**, and the *Łukasiewicz t-norm* is also called the **bold intersection**.

Note that the *Łukasiewicz operators* (previous definition) do *not* in general form a lattice because, for example, they are not in general *idempotent*. Suppose for example $m(p) = \tfrac{1}{2}$ at some point $p$. Then

$$m \vee m \triangleq \min(1, m + m) = \min(1, \tfrac{1}{2} + \tfrac{1}{2}) = 1 \neq m$$
$$m \wedge m \triangleq \max(0, m + m - 1) = \max(0, \tfrac{1}{2} + \tfrac{1}{2} - 1) = 0 \neq m$$

There are several choices for *negations* in a *fuzzy subset logic*. Arguably the "simplest" is the *discrete negation* (Example B.16 page 38). Perhaps the most "common" is the *standard negation* (next definition). More generally there is the *$\lambda$-negation* (Definition 1.20 page 9) which reduces to the standard negation at $\lambda = 0$ and approaches the discrete negation as $\lambda \to \infty$. Alternatively there is also the *Yager negation* (Definition 1.21 page 9) which reduces to the standard negation at $p = 1$.

**Definition 1.19** [23] The function $\neg m(x)$ is the **standard negation** (or **Łukasiewicz negation**) of $m$ if

$$\neg m(x) \triangleq 1 - m(x) \qquad \forall x \in \mathbb{R}.$$

**Definition 1.20** [24] The function $\neg m(x)$ is the **$\lambda$-negation** of a function $m(x)$ if

$$\neg m(x) \triangleq \frac{1 - m(x)}{1 + \lambda m(x)} \qquad \forall \lambda \in (-1 : \infty).$$

**Definition 1.21** [25] The function $\neg m(x)$ is the **Yager negation** of a function $m(x)$ if

$$\neg m(x) \triangleq (1 - m^p)^{1/p} \qquad \forall p \in (0 : \infty).$$

If $\neg m$ is a *$\lambda$-negation*, then the function $\neg$ in a *fuzzy subset lattice* $L$ is a *de Morgan negation* (Definition B.3 page 35) and thus the *de Morgan* properties hold in $L$ (Theorem B.14 page 37). The *standard negation* (Definition 1.19 page 9) is a *$\lambda$-negation* (at $\lambda = 0$) and so the standard negation is also *de Morgan*.

**Theorem 1.22** *Let* $L \triangleq (\mathbb{M}, \vee, \wedge, \neg, 0, 1 ; \leq)$ *be a* LATTICE WITH NEGATION *(Definition B.5 page 35).*

$$\left\{ \begin{array}{c} \neg m(x) \text{ is a } \lambda\text{-NEGATION} \quad \forall m \in \mathbb{M} \\ \textit{(Definition 1.7 page 6)} \end{array} \right\} \implies \left\{ \begin{array}{c} \neg \text{ is a } \text{DE MORGAN NEGATION } \textit{on } L \\ \textit{(Definition B.3 page 35)} \end{array} \right\}$$

✎PROOF:   To be a *de Morgan negation*, $\neg\mathrm{m}_A(x)$ must be *antitone* and *involutory* (Definition B.3 page 35).

$$\mathrm{m}_A(x) \leq \mathrm{m}_B(x) \implies -\mathrm{m}_B(x) \leq -\mathrm{m}_A(x) \qquad \text{by property of } real\ numbers\ \mathbb{R}$$

$$\implies 1 - \mathrm{m}_B(x) \leq 1 - \mathrm{m}_A(x) \qquad \text{by property of } real\ numbers\ \mathbb{R}$$

$$\implies \frac{1 - \mathrm{m}_B(x)}{1 + \lambda\mathrm{m}_B(x)} \leq \frac{1 - \mathrm{m}_A(x)}{1 + \lambda\mathrm{m}_A(x)} \qquad \text{because } 1 + \lambda\mathrm{m} > 0$$

$$\implies \neg\mathrm{m}_B(x) \leq \neg\mathrm{m}_A(x) \qquad \text{by definition of } \lambda\text{-}negation \text{ (Definition 1.20 page 9)}$$

$$\implies \mathrm{m} \text{ is } antitone$$

$$\neg\neg\mathrm{m}_A(x) \triangleq \neg\left(\frac{1 - \mathrm{m}_A(x)}{1 + \lambda\mathrm{m}_A(x)}\right) \qquad \text{by definition of } \lambda\text{-}negation \text{ (Definition 1.20 page 9)}$$

$$\triangleq \frac{1 - \frac{1-\mathrm{m}_A(x)}{1+\lambda\mathrm{m}_A(x)}}{1 + \lambda\frac{1-\mathrm{m}_A(x)}{1+\lambda\mathrm{m}_A(x)}} \qquad \text{by definition of } \lambda\text{-}negation \text{ (Definition 1.20 page 9)}$$

$$= \frac{\left(1 + \lambda\mathrm{m}_A(x)\right) - \left(1 - \mathrm{m}_A(x)\right)}{\left(1 + \lambda\mathrm{m}_A(x)\right) + \lambda\left(1 - \mathrm{m}_A(x)\right)}$$

$$= \frac{(1 + \lambda)\mathrm{m}_A(x)}{1 + \lambda}$$

$$= \mathrm{m}_A(x)$$

$$\implies \neg\mathrm{m} \text{ is } involutory$$

👉

**Corollary 1.23**  *Let* $\mathbf{L} \triangleq (\mathbb{M}, \vee, \wedge, \neg, 0, 1; \leq)$ *be a* LATTICE WITH NEGATION *(Definition B.5 page 35).*

$$\left\{\begin{array}{l} A. \quad \neg\mathrm{m}(x) \text{ is a } \lambda\text{-NEGATION} \quad \forall\mathrm{m}\in\mathbb{M} \quad and \\ B. \quad \neg\mathrm{m}_1 = \mathrm{m}_0 \end{array}\right\} \implies \left\{\begin{array}{l} 1. \quad \neg \text{ is a DE MORGAN NEGATION on } \mathbf{L} \quad and \\ 2. \quad \neg \text{ is a FUZZY NEGATION on } \mathbf{L} \end{array}\right\}$$

✎PROOF:

(1)  Proof for (1): by Theorem 1.22 (page 9)

(2)  Proof for (2): To be a *fuzzy negation*, $\neg\mathrm{m}_A(x)$ must be *antitone*, have *weak double negation*, and have *boundary condition* $\neg\mathrm{m}_1(x) = \mathrm{m}_0(x)$ (Definition B.2 page 35).

  (a)  Proof that $\neg$ is *antitone*: by Theorem 1.22 (page 9).

  (b)  Proof that $\neg$ has *weak double negation*: by Theorem 1.22 (page 9), $\neg$ is *involutory*, which implies $\neg$ has *weak double negation*.

  (c)  Proof that $\neg\mathrm{m}_1(x) = \mathrm{m}_0(x)$: by left hypothesis (B).

👉

We can now define fuzzy subset operators $(\cap, \cup, \mathrm{c})$ in terms of the fuzzy logic operators $(\wedge, \vee, \neg)$ like this (cross reference Remark 1.2 page 5):

$$\begin{array}{llll} \mathrm{m}_\varnothing & \triangleq & 0 & \text{(Definition A.19 page 25)} \\ \mathrm{m}_X & \triangleq & 1 & \text{(Definition A.19 page 25)} \\ \mathrm{m}_{A\cup B} & \triangleq & \mathrm{m}_A \vee \mathrm{m}_B & \text{(Section 1.3 page 7)} \\ \mathrm{m}_{A\cap B} & \triangleq & \mathrm{m}_A \wedge \mathrm{m}_B & \text{(Section 1.3 page 7)} \\ \mathrm{m}_{A^c} & \triangleq & \neg\mathrm{m}_A & \text{(Section 1.3 page 9)} \end{array}$$

In the case of *set indicator functions*, defining $(\wedge, \vee)$ is straightforward. But again here in *fuzzy subset logics*, it is not.

## 1.4    Key theorems

This section contains the following key theorems which under very "reasonable" conditions say very roughly the following about the *fuzzy subset logic* operator pair $(\wedge, \vee)$:

*fuzzy operators idempotency theorem* (Theorem 1.25 page 12):

     *distributive* $\implies$ *idempotent*    and conversely
     *non-idempotent* $\implies$ *non-distributive*

*fuzzy negation idempotency theorem* (Theorem 1.28 page 14):

     *excluded middle* or *non-contradiction* $\implies$ *non-idempotent*    and conversely
     *idempotent* $\implies$ *excluded middle* or *non-contradiction* or both fails

*fuzzy min-max theorem* (Theorem 1.26 page 13):

     *idempotent* $\implies$ $(\wedge, \vee) = (\min, \max)$    and conversely
     $(\wedge, \vee) \neq (\min, \max)$ $\implies$ *non-idempotent*

The *fuzzy min-max boundary theorem* (next theorem) shows that under three pairs of arguably "reasonable" conditions (including *pointwise evaluation*), the functions $\min(\mathfrak{m}, \mathfrak{n})$ and $\max(\mathfrak{m}, \mathfrak{n})$ act as bounds for any possible operators $(\wedge, \vee)$.

**Theorem 1.24**   (fuzzy min-max boundary theorem)   [26] *Let* $\mathbb{M}$ *be a set of* MEMBERSHIP FUNCTIONS *(Definition 1.7 page 6).*

$$\left\{\begin{array}{llllll}
1. & \exists \mathrm{f} \in [0:1]^{[0:1]^2} & \text{such that} & [\mathfrak{m} \wedge \mathfrak{n}](x) = \mathrm{f}[\mathfrak{m}(x), \mathfrak{n}(x)] & \forall \mathfrak{m},\mathfrak{n} \in \mathbb{M} & \text{(POINTWISE EVALUATED)} \quad and \\
2. & \exists \mathrm{g} \in [0:1]^{[0:1]^2} & \text{such that} & [\mathfrak{m} \vee \mathfrak{n}](x) = \mathrm{g}[\mathfrak{m}(x), \mathfrak{n}(x)] & \forall \mathfrak{m},\mathfrak{n} \in \mathbb{M} & \text{(POINTWISE EVALUATED)} \quad and \\
3. & \mathfrak{m} \vee 0 = \mathfrak{m} & & 0 \vee \mathfrak{m} = \mathfrak{m} & \forall \mathfrak{m} \in \mathbb{M} & \text{(DISJUNCTIVE IDENTITY)} \quad and \\
4. & \mathfrak{m} \wedge 1 = \mathfrak{m} & & 1 \wedge \mathfrak{m} = \mathfrak{m} & \forall \mathfrak{m} \in \mathbb{M} & \text{(CONJUNCTIVE IDENTITY)} \quad and \\
5. & \mathfrak{n} \le \mathfrak{p} \implies \mathfrak{m} \vee \mathfrak{n} \le \mathfrak{m} \vee \mathfrak{p} & and & \mathfrak{n} \vee \mathfrak{m} \le \mathfrak{p} \vee \mathfrak{m} & \forall \mathfrak{m},\mathfrak{n},\mathfrak{p} \in \mathbb{M} & \text{(DISJUNCTIVE ISOTONE)} \quad and \\
6. & \mathfrak{n} \le \mathfrak{p} \implies \mathfrak{m} \wedge \mathfrak{n} \le \mathfrak{m} \wedge \mathfrak{p} & and & \mathfrak{n} \wedge \mathfrak{m} \le \mathfrak{p} \wedge \mathfrak{m} & \forall \mathfrak{m},\mathfrak{n},\mathfrak{p} \in \mathbb{M} & \text{(CONJUNCTIVE ISOTONE)} 
\end{array}\right\}$$

$$\implies \left\{\mathfrak{m} \wedge \mathfrak{n} \le \min(\mathfrak{m}, \mathfrak{n}) \quad and \quad \max(\mathfrak{m}, \mathfrak{n}) \le \mathfrak{m} \vee \mathfrak{n} \quad \forall \mathfrak{m} \in \mathbb{M}\right\}$$

PROOF:

$$
\begin{aligned}
\boxed{\max(\mathfrak{m}, \mathfrak{n})} &= \max([\mathfrak{m} \vee 0], [0 \vee \mathfrak{n}]) && \text{by } \textit{disjunctive identity} \text{ property} \\
&\le \max(\mathfrak{m} \vee \mathfrak{n}, 0 \vee \mathfrak{n}) && \text{by } \textit{disjunctive isotone} \text{ property: } 0 \le \mathfrak{n} \implies \mathfrak{m} \vee 0 \le \mathfrak{m} \vee \mathfrak{n} \\
&\le \max(\mathfrak{m} \vee \mathfrak{n}, \mathfrak{m} \vee \mathfrak{n}) && \text{by } \textit{disjunctive isotone} \text{ property: } 0 \le \mathfrak{m} \implies 0 \vee \mathfrak{n} \le \mathfrak{m} \vee \mathfrak{n} \\
&= \boxed{\mathfrak{m} \vee \mathfrak{n}} && \text{by definition of } \max(\cdot, \cdot) \\
\boxed{\mathfrak{m} \wedge \mathfrak{n}} &= \min(\mathfrak{m} \wedge \mathfrak{n}, \mathfrak{m} \wedge \mathfrak{n}) && \text{by definition of } \min(\cdot, \cdot) \\
&\le \min([\mathfrak{m} \wedge 1], [\mathfrak{m} \wedge \mathfrak{n}]) && \text{by } \textit{conjunctive isotone} \text{ property: } \mathfrak{n} \le 1 \implies \mathfrak{m} \wedge \mathfrak{n} \le \mathfrak{m} \wedge 1 \\
&\le \min([\mathfrak{m} \wedge 1], [1 \wedge \mathfrak{n}]) && \text{by } \textit{conjunctive isotone} \text{ property: } \mathfrak{m} \le 1 \implies \mathfrak{m} \wedge \mathfrak{n} \le 1 \wedge \mathfrak{n} \\
&= \boxed{\min(\mathfrak{m}, \mathfrak{n})} && \text{by } \textit{conjunctive identity} \text{ property}
\end{aligned}
$$

How reasonable are the "reasonable conditions" of Theorem 1.24? Let's discuss them briefly:

The strength of the **pointwise evaluation** condition is perhaps more in its simplicity than in it's reasonableness. In mathematics in general, functions are often mapped to other functions in blatant disregard to this property or one like it. Often such a mapping is referred to as an "operator".

---

[26]  [Alsina et al.(1983)Alsina, Trillas, and Valverde] page 16 ⟨§1⟩

🐚 In fuzzy logic, the **identity** properties are "reasonable" because if either the "degree of membership" of $x$ is $\mathfrak{m}_{(x)}$ **or** $x$ has "full membership", then arguably the "degree of membership" of $x$ is $\mathfrak{m}_{(x)}$. Likewise, if both the "degree of membership" of $x$ is $\mathfrak{m}_{(x)}$ **and** $x$ has "absolute non-membership", then arguably the "degree of membership" of $x$ is $\mathfrak{m}_{(x)}$. In order theory, $x \vee 0 = x$ and $x \wedge 1 = x$ are true of any *bounded lattice* (Proposition A.21 page 26). Their commuted counterparts follow from a weakened form of the *commutative* property. Note that all lattices are *commutative* (Theorem A.14 page 25).

🐚 The **isotone** properties are a natural requirement of fuzzy logic—if the "degree of membership" $\mathfrak{m}_{(x)}$ increases, then we might expect that the "degrees of membership" $[\mathfrak{m} \vee \mathfrak{n}]_{(x)}$, $[\mathfrak{n} \vee \mathfrak{m}]_{(x)}$, $[\mathfrak{m} \wedge \mathfrak{n}]_{(x)}$, and $[\mathfrak{n} \wedge \mathfrak{m}]_{(x)}$ to also increase. In order theory, the *isotone* properties hold for all lattices (Proposition A.15 page 25).

The *fuzzy operators idempotency theorem* (next theorem) shows that under a handful of additional arguably "somewhat reasonable" conditions (including the rather "strong" *distributivity* property), the functions $\wedge$ and $\vee$ are both *idempotent*.

**Theorem 1.25**  (fuzzy operators idempotency theorem)  [27] *Let* $\mathbb{M}$ *be a set of* MEMBERSHIP FUNCTION*s* (Definition 1.7 page 6).

$$
\left\{
\begin{array}{llll}
1. & \exists f \in [0:1]^{[0:1]^2} \text{ such that } & [\mathfrak{m} \wedge \mathfrak{n}](x) = f[\mathfrak{m}(x), \mathfrak{n}(x)] \quad \forall \mathfrak{m}, \mathfrak{n} \in \mathbb{M} & \text{(POINTWISE EVALUATED)} \quad and \\
2. & \exists g \in [0:1]^{[0:1]^2} \text{ such that } & [\mathfrak{m} \vee \mathfrak{n}](x) = g[\mathfrak{m}(x), \mathfrak{n}(x)] \quad \forall \mathfrak{m}, \mathfrak{n} \in \mathbb{M} & \text{(POINTWISE EVALUATED)} \quad and \\
3. & 0 \wedge 0 = 0 \quad \bigg| \quad 1 \vee 1 = 1 & & \text{(BOUNDARY CONDITION)} \quad and \\
4. & \mathfrak{m} \vee 0 = \mathfrak{m} \quad \bigg| \quad 0 \vee \mathfrak{m} = \mathfrak{m} & \forall \mathfrak{m} \in \mathbb{M} & \text{(DISJUNCTIVE IDENTITY)} \quad and \\
5. & \mathfrak{m} \wedge 1 = \mathfrak{m} \quad \bigg| \quad 1 \wedge \mathfrak{m} = \mathfrak{m} & \forall \mathfrak{m} \in \mathbb{M} & \text{(CONJUNCTIVE IDENTITY)} \quad and \\
6. & \mathfrak{m} \wedge (\mathfrak{n} \vee \mathfrak{p}) = (\mathfrak{m} \wedge \mathfrak{n}) \vee (\mathfrak{m} \wedge \mathfrak{p}) & \forall \mathfrak{m}, \mathfrak{n}, \mathfrak{p} \in \mathbb{M} & \text{(DISJUNCTIVE DISTRIBUTIVE)} \quad and \\
7. & \mathfrak{m} \vee (\mathfrak{n} \wedge \mathfrak{p}) = (\mathfrak{m} \vee \mathfrak{n}) \wedge (\mathfrak{m} \vee \mathfrak{p}) & \forall \mathfrak{m}, \mathfrak{n}, \mathfrak{p} \in \mathbb{M} & \text{(CONJUNCTIVE DISTRIBUTIVE)}
\end{array}
\right\}
$$

$$
\implies \left\{
\begin{array}{lll}
1. & \mathfrak{m} = \mathfrak{m} \vee \mathfrak{m} \quad \forall \mathfrak{m} \in \mathbb{M} & \text{(DISJUNCTIVE IDEMPOTENT)} \quad and \\
2. & \mathfrak{m} = \mathfrak{m} \wedge \mathfrak{m} \quad \forall \mathfrak{m} \in \mathbb{M} & \text{(CONJUNCTIVE IDEMPOTENT)}
\end{array}
\right\}
$$

✎PROOF:

$$
\begin{array}{ll}
\mathfrak{m} = \mathfrak{m} \wedge 1 & \text{by } \textit{conjunctive identity} \text{ property} \\
\quad = \mathfrak{m} \wedge (1 \vee 1) & \text{by } \textit{boundary condition} \\
\quad = (\mathfrak{m} \wedge 1) \vee (\mathfrak{m} \wedge 1) & \text{by } \textit{conjunctive distributive} \text{ property} \\
\quad = \mathfrak{m} \vee \mathfrak{m} & \text{by } \textit{conjunctive identity} \text{ property} \\
\mathfrak{m} = \mathfrak{m} \vee 0 & \text{by } \textit{disjunctive identity} \text{ property} \\
\quad = \mathfrak{m} \vee (0 \wedge 0) & \text{by } \textit{boundary condition} \\
\quad = (\mathfrak{m} \vee 0) \wedge (\mathfrak{m} \vee 0) & \text{by } \textit{disjunctive distributive} \text{ property} \\
\quad = \mathfrak{m} \wedge \mathfrak{m} & \text{by } \textit{disjunctive identity} \text{ property}
\end{array}
$$

How reasonable are the "reasonable conditions" of Theorem 1.25? Let's discuss them briefly:

🐚 In fuzzy logic, the **boundary conditions** are "reasonable" because if $x$ has both "absolute non-membership" **and** "absolute non-membership", then arguably $x$ has "absolute non-membership".

---

[27] 📖 [Bellman and Giertz(1973)] page 154 $\langle a \vee a = a \wedge a = a \dots (10) \rangle$, 📖 [Alsina et al.(1983)Alsina, Trillas, and Valverde] page 15 $\langle x = G(x, x) \rangle$

Likewise, if $x$ has either " full membership" **or** "full membership", then arguably $x$ has "full membership". In order theory, the boundary conditions are simply a weakened form of the *idempotent* property, which holds for all lattices (Theorem A.14 page 25).

The **distributive** properties hold in *classical logic* (2-valued logic) and more generally in any *Boolean logic*, but not necessarily in any other form of *logic* (Definition C.5 page 50). In order theory, a comparatively small but important class of lattices are *distributive*. But note that in any lattice, the *distributive inequalities* always hold (Theorem A.16 page 25); and if *one* of the distributive properties hold, then they *both* hold (Theorem A.28 page 27).

The *fuzzy min-max theorem* (next theorem) shows that under the *identity* and *isotone* conditions (Theorem 1.24 page 11) and the additional condition of *weak idempotency*, the **only** functions for $(\wedge, \vee)$ are $(\wedge, \vee) = (\min, \max)$….

**Theorem 1.26**   (fuzzy min-max theorem)    [28] *Let* $\mathbb{M}$ *be a set of* MEMBERSHIP FUNCTION$s$ *(Definition 1.7 page 6).*

$$
\left\{
\begin{array}{llll}
1. & \exists f \in [0:1]^{[0:1]^2} \text{ such that } & [\mathtt{m} \wedge \mathtt{n}](x) = f[\mathtt{m}(x), \mathtt{n}(x)] & \forall \mathtt{m},\mathtt{n} \in \mathbb{M} \quad \text{(POINTWISE EVALUATED)} \quad and \\
2. & \exists g \in [0:1]^{[0:1]^2} \text{ such that } & [\mathtt{m} \vee \mathtt{n}](x) = g[\mathtt{m}(x), \mathtt{n}(x)] & \forall \mathtt{m},\mathtt{n} \in \mathbb{M} \quad \text{(POINTWISE EVALUATED)} \quad and \\
3. & \mathtt{m} \vee 0 = \mathtt{m} & 0 \vee \mathtt{m} = \mathtt{m} & \forall \mathtt{m} \in \mathbb{M} \quad \text{(DISJUNCTIVE IDENTITY)} \quad and \\
4. & \mathtt{m} \wedge 1 = \mathtt{m} & 1 \wedge \mathtt{m} = \mathtt{m} & \forall \mathtt{m} \in \mathbb{M} \quad \text{(CONJUNCTIVE IDENTITY)} \quad and \\
5. & \mathtt{m} \wedge \mathtt{m} \geq \mathtt{m} & \mathtt{m} \vee \mathtt{m} \leq \mathtt{m} & \forall \mathtt{m} \in \mathbb{M} \quad \text{(WEAK IDEMPOTENT)} \quad and \\
6. & \mathtt{n} \leq \mathtt{p} \implies \mathtt{m} \vee \mathtt{n} \leq \mathtt{m} \vee \mathtt{p} \text{ and} & \mathtt{n} \vee \mathtt{m} \leq \mathtt{p} \vee \mathtt{m} & \forall \mathtt{m},\mathtt{n},\mathtt{p} \in \mathbb{M} \quad \text{(DISJUNCTIVE ISOTONE)} \quad and \\
7. & \mathtt{n} \leq \mathtt{p} \implies \mathtt{m} \wedge \mathtt{n} \leq \mathtt{m} \wedge \mathtt{p} \text{ and} & \mathtt{n} \wedge \mathtt{m} \leq \mathtt{p} \wedge \mathtt{m} & \forall \mathtt{m},\mathtt{n},\mathtt{p} \in \mathbb{M} \quad \text{(CONJUNCTIVE ISOTONE)}
\end{array}
\right\}
$$

$$
\implies \left\{
\begin{array}{lll}
1. & \mathtt{m} \vee \mathtt{n} = \max(\mathtt{m}, \mathtt{n}) & \forall \mathtt{m},\mathtt{n} \in \mathbb{M} \quad and \\
2. & \mathtt{m} \wedge \mathtt{n} = \min(\mathtt{m}, \mathtt{n}) & \forall \mathtt{m},\mathtt{n} \in \mathbb{M}
\end{array}
\right\}
$$

✎PROOF:

$$
\begin{array}{lll}
\boxed{\max(\mathtt{m}, \mathtt{n})} \leq \boxed{\mathtt{m} \vee \mathtt{n}} & & \text{by } \textit{fuzzy min-max boundary theorem (Theorem 1.24 page 11)} \\
\leq \max(\mathtt{m}, \mathtt{n}) \vee \mathtt{n} & & \text{by } \textit{disjunctive isotone property: } \mathtt{m} \leq \max(\mathtt{m}, \mathtt{n}) \\
\leq \max(\mathtt{m}, \mathtt{n}) \vee \max(\mathtt{m}, \mathtt{n}) & & \text{by } \textit{disjunctive isotone property: } \mathtt{n} \leq \max(\mathtt{m}, \mathtt{n}) \\
\leq \boxed{\max(\mathtt{m}, \mathtt{n})} & & \text{by } \textit{weak idempotent property} \\
\boxed{\min(\mathtt{m}, \mathtt{n})} \leq \min(\mathtt{m}, \mathtt{n}) \wedge \min(\mathtt{m}, \mathtt{n}) & & \text{by } \textit{weak idempotent property} \\
\leq \mathtt{m} \wedge \min(\mathtt{m}, \mathtt{n}) & & \text{by } \textit{isotone property of } \wedge: \min(\mathtt{m}, \mathtt{n}) \leq \mathtt{m} \\
\leq \boxed{\mathtt{m} \wedge \mathtt{n}} & & \text{by } \textit{isotone property of } \wedge: \min(\mathtt{m}, \mathtt{n}) \leq \mathtt{n} \\
\leq \boxed{\min(\mathtt{m}, \mathtt{n})} & & \text{by } \textit{fuzzy min-max boundary theorem (Theorem 1.24 page 11)}
\end{array}
$$

How reasonable are the "reasonable conditions" of Theorem 1.26? Let's discuss them briefly:

One way to get the *weak idempotent* property or even the stronger *idempotent* property is to force $(\min, \max)$ to have the *boundary* and *distributive* properties (Theorem 1.25 page 12). However, this is arguably a kind of sledge hammer approach and is not really necessary.

---

[28] This result is very similar to the celebrated result of Bellman and Giertz (1973): 📖 [Bellman and Giertz(1973)] pages 153–154 ⟨§4⟩

In fuzzy logic, even the stronger *idempotent* property is arguably "reasonable" because if an element $x$ both has a "degree of membership" $\mathbb{m}(x)$ **and** a "degree of membership" $\mathbb{m}(x)$, then arguably $x$ has a "degree of membership" $\mathbb{m}(x)$. Likewise, if $x$ either has a "degree of member-ship" $\mathbb{m}(x)$ **or** a "degree of membership" $\mathbb{m}(x)$, then arguably $x$ has a "degree of membership" $\mathbb{m}(x)$. In order theory, all lattices are *idempotent* (Theorem A.14 page 25). But, again, here we only require *weak idempotency*, not *idempotency*.

**Corollary 1.27** (Bellman-Giertz 1973 theorem) [29] *Let* $\mathbb{M}$ *be a set of* MEMBERSHIP FUNCTION*s (Definition 1.7 page 6).*

$$\left\{\begin{array}{llll}
1. & \exists f \in [0:1]^{[0:1]^2} \quad \text{such that} \quad [\mathbb{m} \wedge \mathbb{n}](x) = f[\mathbb{m}(x), \mathbb{n}(x)] & \forall \mathbb{m}, \mathbb{n} \in \mathbb{M} & \text{(POINTWISE EVALUATED)} \quad \text{and} \\
2. & \exists g \in [0:1]^{[0:1]^2} \quad \text{such that} \quad [\mathbb{m} \vee \mathbb{n}](x) = g[\mathbb{m}(x), \mathbb{n}(x)] & \forall \mathbb{m}, \mathbb{n} \in \mathbb{M} & \text{(POINTWISE EVALUATED)} \quad \text{and} \\
3. & \mathbb{m} \vee 0 = \mathbb{m} \qquad\qquad\qquad 0 \vee \mathbb{m} = \mathbb{m} & \forall \mathbb{m} \in \mathbb{M} & \text{(DISJUNCTIVE IDENTITY)} \quad \text{and} \\
4. & \mathbb{m} \wedge 1 = \mathbb{m} \qquad\qquad\qquad 1 \wedge \mathbb{m} = \mathbb{m} & \forall \mathbb{m} \in \mathbb{M} & \text{(CONJUNCTIVE IDENTITY)} \quad \text{and} \\
5. & \mathbb{m} \wedge \mathbb{m} = \mathbb{m} \qquad\qquad\qquad \mathbb{m} \vee \mathbb{m} = \mathbb{m} & \forall \mathbb{m} \in \mathbb{M} & \text{(IDEMPOTENT)} \quad \text{and} \\
6. & \mathbb{n} \le \mathbb{p} \implies \mathbb{m} \vee \mathbb{n} \le \mathbb{m} \vee \mathbb{p} \quad and \quad \mathbb{n} \vee \mathbb{m} \le \mathbb{p} \vee \mathbb{m} & \forall \mathbb{m}, \mathbb{n}, \mathbb{p} \in \mathbb{M} & \text{(DISJUNCTIVE ISOTONE)} \quad \text{and} \\
7. & \mathbb{n} \le \mathbb{p} \implies \mathbb{m} \wedge \mathbb{n} \le \mathbb{m} \wedge \mathbb{p} \quad and \quad \mathbb{n} \wedge \mathbb{m} \le \mathbb{p} \wedge \mathbb{m} & \forall \mathbb{m}, \mathbb{n}, \mathbb{p} \in \mathbb{M} & \text{(CONJUNCTIVE ISOTONE)}
\end{array}\right.$$

$$\implies \left\{\begin{array}{lll}
1. & \mathbb{m} \vee \mathbb{n} = \max(\mathbb{m}, \mathbb{n}) & \forall \mathbb{m}, \mathbb{n} \in \mathbb{M} \quad and \\
2. & \mathbb{m} \wedge \mathbb{n} = \min(\mathbb{m}, \mathbb{n}) & \forall \mathbb{m}, \mathbb{n} \in \mathbb{M}
\end{array}\right\}$$

✎PROOF:    This follows directly from Theorem 1.26 (page 13).          ✏

One big difficulty in *fuzzy subset logic* (Definition 1.11 page 6) is that under "reasonable" conditions, if the fuzzy subset logic is required to have either the *excluded middle* property *or* the *non-contradiction* property (Boolean algebras have both), then the fuzzy subset logic cannot be *idempotent* (next theorem). Furthermore, if a structure is not idempotent, then it is *not a lattice* (Theorem A.14 page 25).

**Theorem 1.28** (fuzzy negation idempotency theorem)   *Let* $\boldsymbol{L} \triangleq (\mathbb{M}, \vee, \wedge, \neg, 0, 1; \le)$ *be a* FUZZY SUB-SET LOGIC *(Definition 1.11 page 6). Let* $(\wedge, \vee)$ *be* POINTWISE EVALUATED *(Definition 1.12 page 7). If there exists* $p$ *such that* $\neg\mathbb{m}(p) = \mathbb{m}(p) \in (0:1)$ *then*

$$\begin{array}{llll}
& & \underbrace{\qquad\qquad\qquad}_{\text{(FIXED POINT CONDITION)}} & \\
(A). & \mathbb{m} \vee \neg\mathbb{m} = 1 \quad \forall \mathbb{m} \in \mathbb{M} & \text{(EXCLUDED MIDDLE)} & \implies \mathbb{m} \vee \mathbb{m} \ne \mathbb{m} \quad \text{(NON-IDEMPOTENT)} \\
(B). & \mathbb{m} \wedge \neg\mathbb{m} = 0 \quad \forall \mathbb{m} \in \mathbb{M} & \text{(NON-CONTRADICTION)} & \implies \mathbb{m} \wedge \mathbb{m} \ne \mathbb{m} \quad \text{(NON-IDEMPOTENT)}
\end{array}$$

✎PROOF:

$$\begin{array}{ll}
1 = \mathbb{m}(p) \vee \neg\mathbb{m}(p) & \text{by } \textit{excluded middle} \text{ hypothesis (A)} \\
\phantom{1} = \mathbb{m}(p) \vee \mathbb{m}(p) & \text{by fixed point hypothesis} \\
\phantom{1} = \mathbb{m}(p) & \textbf{if } \vee \text{ is } \textit{idempotent} \\
\implies \neg\mathbb{m}(p) = 0 & \text{because } \neg\mathbb{m}(p) = \neg 1 = 0 \\
\implies \mathbb{m}(p) = 0 & \text{by fixed point hypothesis} \\
\implies \text{contradiction} & \text{because } \mathbb{m}(p) = 1 \ne 0 = \mathbb{m}(p) \text{ is a contradiction} \\
\implies \vee \text{ is } \textit{non-idempotent} & \\
\\
0 = \mathbb{m}(p) \wedge \neg\mathbb{m}(p) & \text{by } \textit{non-contradiction} \text{ hypothesis (B)}
\end{array}$$

---

[29] 📖 [Bellman and Giertz(1973)] pages 153–154 ⟨§4⟩

| | | |
|---|---|---|
| $= \mathtt{m}(p) \wedge \mathtt{m}(p)$ | by fixed point hypothesis |
| $= \mathtt{m}(p)$ | **if** $\wedge$ is *idempotent* |
| $\implies$ | $\neg\mathtt{m}(p) = 1$ | because $\neg\mathtt{m}(p) = \neg 0 = 1$ |
| $\implies$ | $\mathtt{m}(p) = 0$ | by fixed point hypothesis |
| $\implies$ | contradiction | because $\mathtt{m}(p) = 0 \neq 1 = \mathtt{m}(p)$ is a contradiction |
| $\implies$ | $\wedge$ is *non-idempotent* |

How reasonable are the "reasonable conditions" of Theorem 1.28? Let's discuss them briefly: One of these "reasonable conditions" is that at some point $p$, $\neg\mathtt{m}(p) = \mathtt{m}(p) \in (0 : 1)$. Because fuzzy negations are *antitone*, in some cases this is arguably a "reasonable" assumption, especially if $\mathtt{m}(x)$ is *continuous* and *strictly antitone*. However, be warned that it is not always the case that there is such a point $p$ in a fuzzy subset logic ($\mathbb{M}$, $\vee$, $\wedge$, $\neg$, 0, 1 ; $\leq$). For example, under standard negation, and if the universal set is finite, then it is certainly possible that $p$ does not exist, as in the example illustrated to the right with $X \triangleq \{a, b, c, d\}$ :

| $x$ | $\mathtt{m}(x)$ | $\neg\mathtt{m}(x)$ |
|---|---|---|
| $d$ | 1 | 0 |
| $c$ | ¾ | ¼ |
| $b$ | ¼ | ¾ |
| $a$ | 0 | 1 |

**Corollary 1.29**  (Dubois-Padre 1980 theorem)  [30] *Let* $\boldsymbol{L} \triangleq$ ($\mathbb{M}$, $\vee$, $\wedge$, $\neg$, 0, 1 ; $\leq$) *be a* FUZZY SUBSET LOGIC *(Definition 1.11 page 6). Let* ($\wedge$, $\vee$) *be* POINTWISE EVALUATED *(Definition 1.12 page 7).*
*If* $\neg(x)$ *is* CONTINUOUS *and* STRICTLY ANTITONE *then*

*(A).*  $\mathtt{m} \vee \neg\mathtt{m} = 1$  $\forall_{\mathtt{m} \in \mathbb{M}}$  *(EXCLUDED MIDDLE)*  $\implies$  $\mathtt{m} \vee \mathtt{m} \neq \mathtt{m}$  *(NON-IDEMPOTENT)*
*(B).*  $\mathtt{m} \wedge \neg\mathtt{m} = 0$  $\forall_{\mathtt{m} \in \mathbb{M}}$  *(NON-CONTRADICTION)*  $\implies$  $\mathtt{m} \wedge \mathtt{m} \neq \mathtt{m}$  *(NON-IDEMPOTENT)*

✎PROOF:   This follows directly from Theorem 1.28 (page 14).

## 1.5   Examples of non-ortho and non-Boolean fuzzy subset

This section presents some examples of *fuzzy subset logic*s. They all have "problems". The problem of the first example is just that is a kind of trivial fuzzy subset logic in that it is 2-valued and equivalent to the classical subset logic. In all the other examples, the "problem" involves not having one or more of the following four properties:

(1).  *disjunctive idempotence*:   $x \vee x = x$   and
(2).  *conjunctive idempotence*:   $x \wedge x = x$   and
(3).  *excluded middle*:   $x \vee \neg x = 1$   and
(4).  *non-contradiction*:   $x \wedge \neg x = 0$

Actually, this is a problem only as far as not having an *ortho* or *Boolean* logic is a problem—because all *ortho logics* and all *Boolean logics* have these properties. And so if even one is missing, the logic is neither an *ortho logic* nor a *Boolean logic*. Also note that if a logic does not have both (1) and (2), then it cannot even be constructed on a *lattice* at all…and as defined in this paper, is not even a *logic*.

**Example 1.30**  Consider the structure $\boldsymbol{L} \triangleq$ ($\mathbb{M}$, $\vee$, $\wedge$, $\neg$, 0, 1 ; $\leq$) in Figure 2 page 16 (A).

1.  $\boldsymbol{L}$ is a *Boolean lattice* (Definition A.41 page 30).

---

[30] 🖙 [Dubois and Padre(1980)], page 62 ⟨P1, requires $\neg(x)$ be to *continuous* and *strictly antitone*⟩
🖙 [Fodor and Yager(2000)], pages 130–131 ⟨Theorem 2, reference to previous without proof⟩

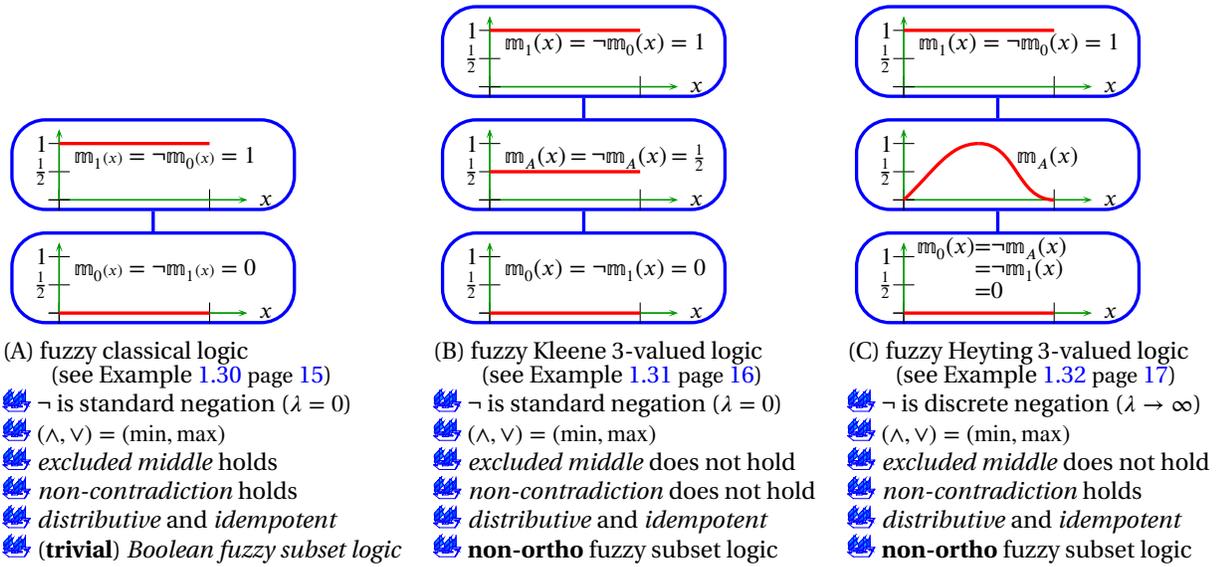

(A) fuzzy classical logic
(see Example 1.30 page 15)
🔥 ¬ is standard negation ($\lambda = 0$)
🔥 ($\wedge$, $\vee$) = (min, max)
🔥 *excluded middle* holds
🔥 *non-contradiction* holds
🔥 *distributive* and *idempotent*
🔥 (**trivial**) *Boolean fuzzy subset logic*

(B) fuzzy Kleene 3-valued logic
(see Example 1.31 page 16)
🔥 ¬ is standard negation ($\lambda = 0$)
🔥 ($\wedge$, $\vee$) = (min, max)
🔥 *excluded middle* does not hold
🔥 *non-contradiction* does not hold
🔥 *distributive* and *idempotent*
🔥 **non-ortho** fuzzy subset logic

(C) fuzzy Heyting 3-valued logic
(see Example 1.32 page 17)
🔥 ¬ is discrete negation ($\lambda \to \infty$)
🔥 ($\wedge$, $\vee$) = (min, max)
🔥 *excluded middle* does not hold
🔥 *non-contradiction* holds
🔥 *distributive* and *idempotent*
🔥 **non-ortho** fuzzy subset logic

Figure 2:  *fuzzy logic*s (Definition C.5 page 50) on *linear lattice*s (Definition A.11 page 24) $\mathbf{L}_2$ and $\mathbf{L}_3$

2. The function ¬ is an *ortho negation* (Definition B.3 page 35) (and hence also is a *fuzzy negation* Definition B.2 page 35, Figure 9 page 34).

3. The negation ¬$\mathbb{m}$ of each *membership function* $\mathbb{m}$ (Definition 1.7 page 6) is the *standard negation* (Definition 1.19 page 9).

4. $\mathbf{L}$ together with the *classical implication* (Example C.4 page 46) is the *classical logic* (Example C.6 page 50) and is also a *fuzzy logic* (Definition C.5 page 50).

5. Because the membership functions $\mathbb{m}(x)$ equal 0 or 1 only, the fuzzy subsets are equivalent to crisp sets.

6. $\mathbf{L}$ is *linear* (Definition A.11 page 24) and therefore *distributive* (Theorem A.30 page 27); and therefore ($\wedge$, $\vee$) are *idempotent* (Theorem 1.25 page 12).

7. The *excluded middle* and *non-contradiction* properties hold in $\mathbf{L}$, but $\mathbf{L}$ is also *idempotent*. This does not contradict Theorem 1.28 (page 14), because ¬ does not satisfy the *fixed point condition* (there is no point $p$ such that ¬$\mathbb{m}(p) = \mathbb{m}(p) \in (0 : 1)$).

**Example 1.31**  Consider the structure $\mathbf{L} \triangleq (\mathbb{M}, \vee, \wedge, \neg, 0, 1 ; \le)$ in Figure 2 page 16 (B).

1. The function ¬ is a *Kleene negation* (Definition B.3 page 35) (and hence a *de Morgan negation*), and is also a *fuzzy negation* (Example B.25 page 41).

2. The negation ¬$\mathbb{m}$ of each *membership function* $\mathbb{m}$ is the *standard negation* because for example $\mathbb{m}_A(x) \triangleq \frac{1}{2} = 1 - \frac{1}{2} = 1 - \mathbb{m}_A(x) \triangleq \neg\mathbb{m}_A(x)$.

3. $\mathbf{L}$ is *linear* (Definition A.11 page 24) and therefore *distributive* (Definition A.27 page 27, Theorem A.30 page 27); and therefore ($\wedge$, $\vee$) are *idempotent* (Theorem 1.25 page 12).

4. $\mathbf{L}$ does not have the *excluded middle* property because
$\mathbb{m}_A \vee \neg\mathbb{m}_A = \mathbb{m}_A \vee \mathbb{m}_A = \mathbb{m}_A \triangleq \frac{1}{2} \ne 1$.

5. $\mathbf{L}$ does not have the *non-contradiction* property because
$\mathbb{m}_A \wedge \neg\mathbb{m}_A \mathbb{m}_A \wedge \mathbb{m}_A = \mathbb{m}_A \triangleq \frac{1}{2} \ne 0$.

6. ($\wedge$, $\vee$) = (min, max) (Definition 1.15 page 7), which together with the *idempotence* property agrees with Theorem 1.26 (page 13).

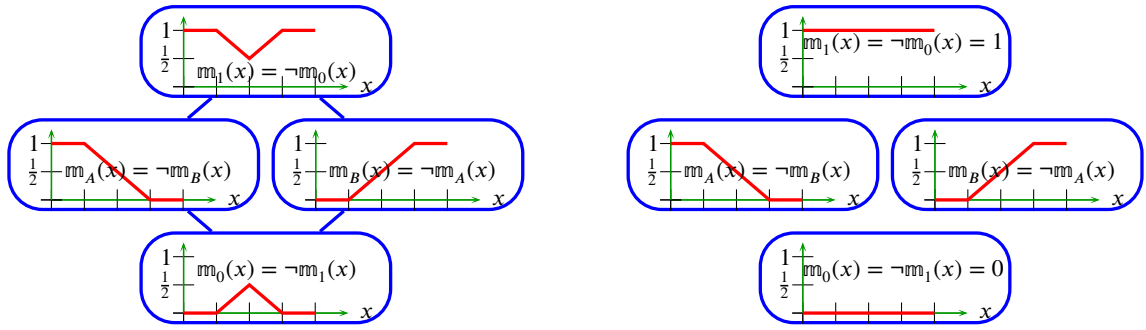

(A) min-max operators (see Example 1.33 page 17)
- ✍ ¬m(x) is standard negation (λ = 0)
- ✍ *excluded middle* and *non-contradiction* do not hold
- ✍ *distributive* and *idempotent*
- ✍ **non-ortho** fuzzy subset logic

(B) Łukasiewicz operators (see Example 1.34 page 18)
- ✍ ¬m(x) is standard negation (λ = 0)
- ✍ *excluded middle* and *non-contradiction* hold
- ✍ *non-idempotent*
- ✍ **not a logic**

Figure 3: *fuzzy logic* on $\boldsymbol{M_2}$ lattice

7. **L** together with the *classical implication* (Example C.4 page 46) is a *Kleene 3-valued logic* (Example C.7 page 51) and also a *fuzzy logic* (Definition C.5 page 50).

**Example 1.32** Consider the structure $\boldsymbol{L} \triangleq (\,\mathbb{M},\ \vee,\ \wedge,\ \neg,\ 0,\ 1\,;\ \leq)$ in Figure 2 page 16 (C).

1. The function ¬ is an *intuitionistic negation* (Definition B.3 page 35) (and hence also a *fuzzy negation* Example B.26 page 41).
2. The negation ¬m of each *membership function* m is the *discrete negation* (Example B.16 page 38).
3. **L** does **not** have the *excluded middle* property because $m_A \vee \neg m_A \neq 1$.
4. **L** does have the *non-contradiction* property.
5. **L** is *linear* (Definition A.11 page 24) and therefore *distributive* (Definition A.27 page 27, Theorem A.30 page 27); and therefore (∧, ∨) are *idempotent* (Theorem 1.25 page 12).
6. Note that having both *non-contradiction* and *idempotency* does not conflict with Theorem 1.28 (page 14) because it does not satisfy the *fixed point condition*.
7. (∧, ∨) = (min, max) (Definition 1.15 page 7), which together with the *idempotence* property agrees with (Theorem 1.26 page 13).
8. **L** together with the *classical implication* (Example C.4 page 46) is a *Heyting 3-valued logic* (Example C.10 page 52) and also a *fuzzy logic* (Definition C.5 page 50).

**Example 1.33** Consider the structure **L** illustrated in Figure 3 page 17 (A).

1. The function ¬ is a *Kleene negation* (Definition B.3 page 35) and also a *fuzzy negation* (Definition B.2 page 35).
2. The negation ¬m of each membership function m is the *standard negation* (Definition 1.19 page 9).
3. The ∧ and ∨ operators are the *min-max operators* (Definition 1.15 page 7).
4. Because (∧, ∨) = (min, max), **L** is a lattice (Proposition 1.16 page 7).
5. Because **L** is a lattice, **L** is *idempotent* (Theorem A.14 page 25). Conversely, *idempotence* and (min, max) are in agreement with Theorem 1.26 (page 13).
6. **L** does *not* have the *excluded middle* property because $m_A \vee \neg m_A = m_1 \neq 1$.
7. **L** does *not* have the *non-contradiction* property because $m_A \wedge \neg m_A = m_0 \neq 0$.
8. The *idempotence* property is *not* in disagreement with Theorem 1.28 (page 14) because **L** does not have the *excluded middle* or *non-contradiction* properties.

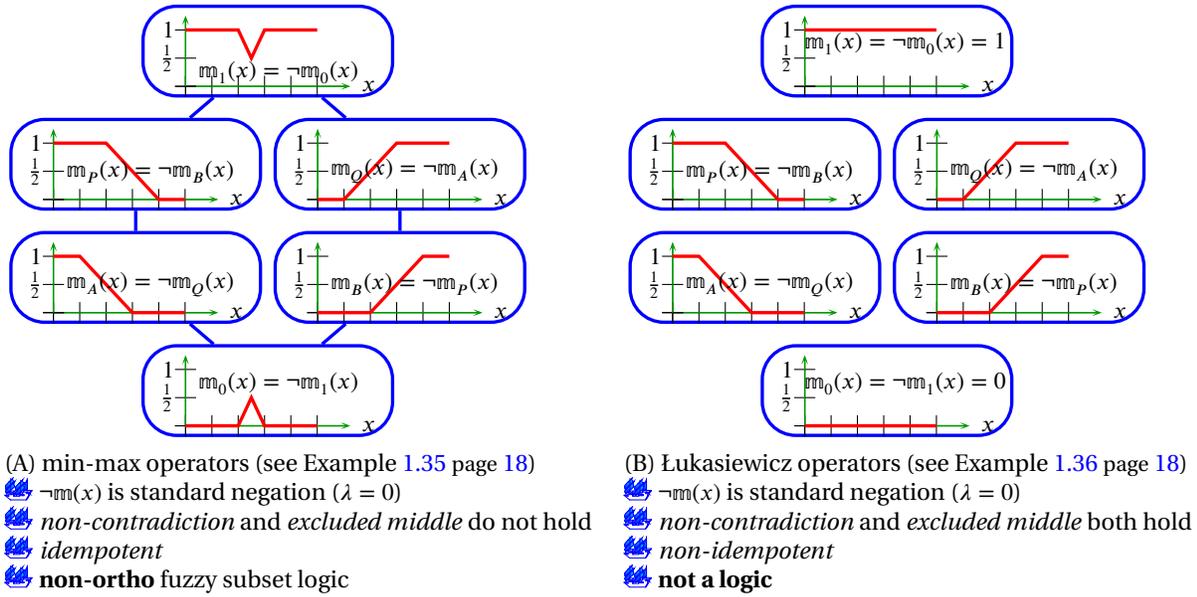

(A) min-max operators (see Example 1.35 page 18)
🜨 ¬𝕞(x) is standard negation (λ = 0)
🜨 *non-contradiction* and *excluded middle* do not hold
🜨 *idempotent*
🜨 **non-ortho** fuzzy subset logic

(B) Łukasiewicz operators (see Example 1.36 page 18)
🜨 ¬𝕞(x) is standard negation (λ = 0)
🜨 *non-contradiction* and *excluded middle* both hold
🜨 *non-idempotent*
🜨 **not a logic**

Figure 4:  *fuzzy logic* on $O_6$ *lattice*

9. *L* together with any of the six *implication* functions listed in Example C.4 (page 46) is a *fuzzy subset logic* (Definition 1.11 page 6).

**Example 1.34**   Consider the structure *L* illustrated in Figure 3 page 17 (B).

1. The function ¬ is an *ortho negation* (Definition B.3 page 35) (and thus also a *fuzzy negation*).
2. The negation ¬𝕞 of each membership function 𝕞 is the *standard negation* (Definition 1.19 page 9).
3. The ∧ and ∨ operators are the *Łukasiewicz operators* (Definition 1.18 page 9). Under these operators, *L* has the *non-contradiction* and *excluded middle* properties, but *L* is *not idempotent* (e.g. $𝕞_A ∨ 𝕞_A ≠ 𝕞_A$), and so *L* is not a lattice (Theorem 1.28 page 14, Theorem A.14 page 25).

**Example 1.35**   Consider the structure *L* illustrated in Figure 4 page 18 (A).

1. The function ¬ is an *ortho negation* (Definition B.3 page 35) (and thus also a *fuzzy negation*).
2. The negation ¬𝕞 of each membership function 𝕞 is the *standard negation* (Definition 1.19 page 9).
3. The ∧ and ∨ operators are the *min-max operators* (Definition 1.15 page 7).
4. Because (∧, ∨) = (min, max), *L* is a lattice (Proposition 1.16 page 7).
5. Because *L* is a lattice, *L* is *idempotent* (Theorem A.14 page 25). Conversely, *idempotence* and (min, max) are in agreement with Theorem 1.26 (page 13).
6. *L* does *not* have the *excluded middle* property because $𝕞_A ∨ ¬𝕞_A = 𝕞_1 ≠ 1$.
7. *L* does *not* have the *non-contradiction* property because $𝕞_A ∧ ¬𝕞_A = 𝕞_0 ≠ 0$.
8. *L* together with any of the six *implication* functions listed in Example C.4 (page 46) is a *fuzzy subset logic* (Definition 1.11 page 6).

**Example 1.36**   Consider the structure *L* illustrated in Figure 4 page 18 (B).

1. The function ¬ is an *ortho negation* (Definition B.3 page 35) (and thus also a *fuzzy negation*).
2. The negation ¬𝕞 of each membership function 𝕞 is the *standard negation* (Definition 1.19 page 9).

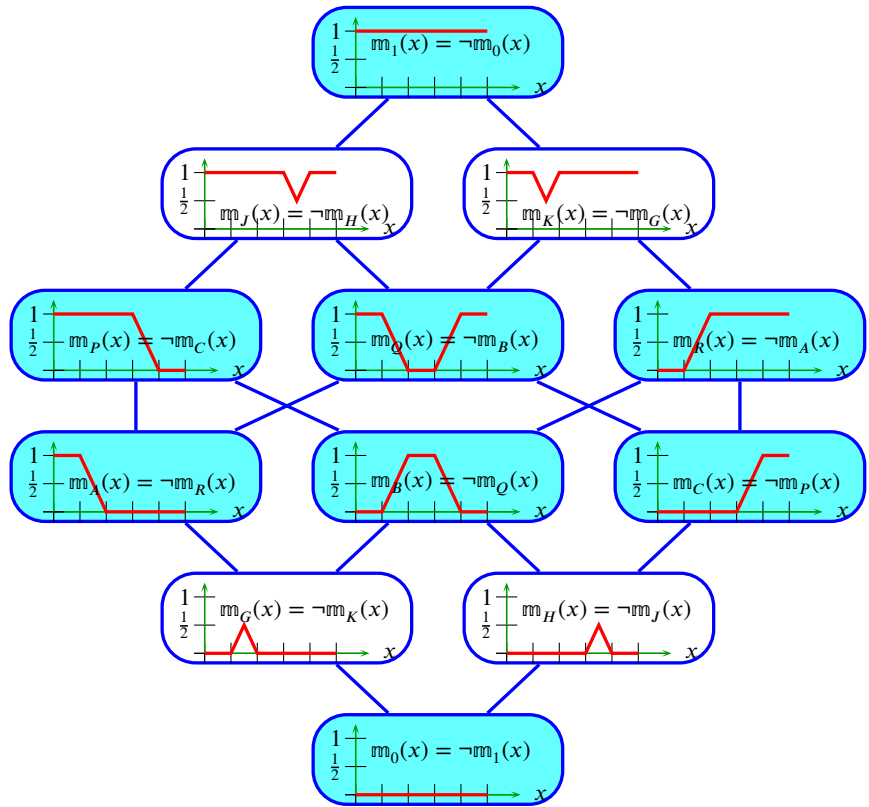

🚢 min-max operators (see Example 1.37 page 19)
🚢 ¬𝕞($x$) is standard negation ($\lambda = 0$)
🚢 *exluded middle* and *non-contradiction* do not hold
🚢 *distributive* and *idempotent*
🚢 **non-ortho** fuzzy subset logic

Figure 5: *fuzzy logic* on lattice with $\mathbf{L}_2^3$ *sublattice*

3. The $\wedge$ and $\vee$ operators are the *Łukasiewicz operators* (Definition 1.18 page 9). Under these operators, **L** has the *non-contradiction* and *excluded middle* properties, but **L** is *not idempotent*, and so **L** is not a lattice (Theorem 1.28 page 14).

**Example 1.37**  Consider the structure **L** illustrated in Figure 5 (page 19).

1. The function ¬ is an *ortho negation* (Definition B.3 page 35) (and thus also a *fuzzy negation*).
2. The negation ¬𝕞 of each membership function 𝕞 is the *standard negation* (Definition 1.19 page 9).
3. The $\wedge$ and $\vee$ operators are the *min-max operators* (Definition 1.15 page 7).
4. Because $(\wedge, \vee) = (\min, \max)$, **L** is a lattice (Proposition 1.16 page 7).
5. Because **L** is a lattice, **L** is *idempotent* (Theorem A.14 page 25). Conversely, *idempotence* and $(\min, \max)$ are in agreement with Theorem 1.26 (page 13).
6. **L** does *not* have the *excluded middle* property because for example
   $$\mathbb{m}_A \vee \neg\mathbb{m}_A = \mathbb{m}_A \vee \mathbb{m}_R = \mathbb{m}_K \neq 1.$$
7. **L** does *not* have the *non-contradiction* property because for example
   $$\mathbb{m}_A \wedge \neg\mathbb{m}_A = \mathbb{m}_A \wedge \mathbb{m}_R = \mathbb{m}_G \neq 0.$$

8. *L* does not contain $M_3$ or $N_5$ and so is *distibutive* (Theorem A.30 page 27). (also cross reference Theorem 1.25 page 12 and Theorem 1.28 page 14).

9. *L* is *non-Boolean*, but has an $L_2^3$ Boolean sublattice ( shaded in Figure 5).

## 2   Boolean and ortho fuzzy subset logics

The Introduction described the problem of constructing *Boolean fuzzy subet logics* and more generally *ortho fuzzy subet logics*. It also briefly described a "solution". This section presents this solution in more detail.

Simply put, a solution is available if we are willing to give up the *pointwise evaluation* condition (Definition 1.12 page 7). In particular, we can proceed as follows:

(1) We give up the *pointwise evaluation* condition.

(2) We define the *ordering relation* (Definition A.1 page 22) $\leq$ in the *fuzzy subset logic* $L \triangleq (\mathbb{M}, \vee, \wedge, \neg, 0, 1 ; \leq)$ to be the *pointwise ordering relation* (Definition A.7 page 23).

(3) In a *lattice* (Definition A.11 page 24), the definitions of the ordering relation $\leq$ and operators $(\wedge, \vee)$ are not independent—the ordering relation defines the operators (Definition A.9 page 24, Definition A.8 page 24) and the operators define the ordering relation (Proposition A.10 page 24).

(4) Traditionally in fuzzy logic literature, we first define a *pointwise evaluated* (Definition 1.12 page 7) pair of operators $(\wedge, \vee)$, and then define the ordering relation $\leq$ in terms of $(\wedge, \vee)$. For example, if $(\wedge, \vee) = (\min, \max)$, then

$$x \leq y \overset{\text{def}}{\iff} \max(x, y) = y$$
$$x \leq y \overset{\text{def}}{\iff} \min(x, y) = x$$

(5) However, here we take a kind of converse approach: We first define a *pointwise* ordering relation $\leq$ (Definition A.7 page 23), and then define the operators $(\wedge, \vee)$ in terms of $\leq$. In doing so, $(\wedge, \vee)$ may possibly no longer satisfy the *pointwise evaluation* condition.

(6) By carefully constructing a set of *membership function*s (Definition 1.7 page 6) $\mathbb{M}$, we can construct *fuzzy subset logic*s (Definition 1.11 page 6) on Boolean and other types of lattice structures.

(7) A fuzzy subset logic then inherits the properties of the lattice it is constructed on. So, for example, if a fuzzy subset logic is constructed on a Boolean lattice, then that fuzzy subset logic is also *Boolean* with all the properties of a Boolean algebra (Theorem A.42 page 30) including the *non-contradiction*, *excluded middle*, *idempotent*, and *distributive* properties.

(8) Despite Theorem 1.26 page 13 and Theorem 1.28, this is all possible because $(\wedge, \vee)$ is no longer *pointwise evaluated* (Definition 1.12 page 7). The result of, say, $[\mathtt{m} \vee \mathtt{n}](x)$ at the point $x$ is no longer necessarily the result of the two values $\mathtt{m}(x)$ and $\mathtt{n}(x)$ alone, but instead $[\mathtt{m} \vee \mathtt{n}](x)$ at the point $x$ may be the result of entire membership functions in the structure or even the position of $\mathtt{m}$ and $\mathtt{n}$ in the structure.

(9) Examples follow.

**Example 2.1**  Consider the structure $L \triangleq (\mathbb{M}, \vee, \wedge, \neg, 0, 1 ; \leq)$ with $M \triangleq \{\mathtt{m}_0, \mathtt{m}_A, \mathtt{m}_B, \mathtt{m}_1\}$ illustrated in Figure 6 page 21 (A).

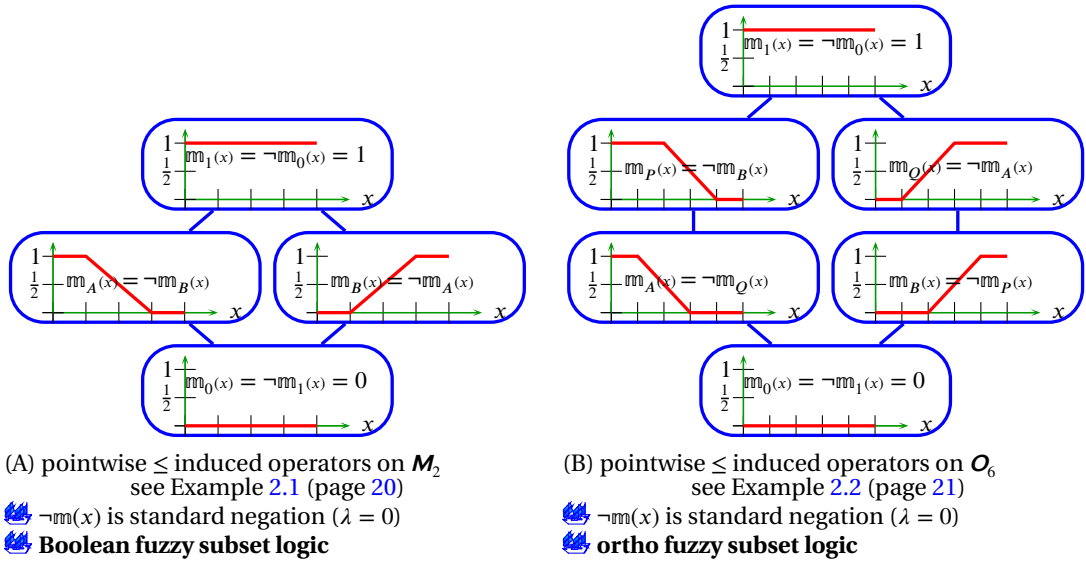

(A) pointwise $\le$ induced operators on $\boldsymbol{M}_2$
see Example 2.1 (page 20)
📖 $\neg m(x)$ is standard negation ($\lambda = 0$)
📖 **Boolean fuzzy subset logic**

(B) pointwise $\le$ induced operators on $\boldsymbol{O}_6$
see Example 2.2 (page 21)
📖 $\neg m(x)$ is standard negation ($\lambda = 0$)
📖 **ortho fuzzy subset logic**

Figure 6: *fuzzy logic* on $\boldsymbol{M}_2$ *lattice*

1. The function $\neg$ is an *ortho negation* (Definition B.3 page 35) (and thus also a *fuzzy negation*).
2. The negation $\neg m$ of each membership function $m$ is the *standard negation* (Definition 1.19 page 9).
3. $\boldsymbol{L}$ is very similar to the structure in Example 1.34 (page 18), which fails to even be a logic.
4. However the structure of this example has a valid ordering relation $\le$ (*pointwise ordering relation*), has valid operators ($\wedge$, $\vee$) defined in terms of $\le$ (Definition A.9 page 24, Definition A.8 page 24), and is a *Boolean lattice* with all the accompanying Boolean properties including the *non-contradiction, excluded middle, idempotency,* and *distributivity*.
5. In this example, the operators are no longer *Łukasiewicz operators* (as in Example 1.34), but some other operators (not explicitly given in terms of a function of the form given in Theorem 1.26 (page 13)).
6. This Boolean lattice together with the *classical implication* (Example C.4 page 46) is an *ortho logic* (and thus also a *fuzzy subset logic*—Definition 1.11 page 6).

**Example 2.2** Consider the structure $\boldsymbol{L} \triangleq (\mathbb{M}, \vee, \wedge, \neg, 0, 1 ; \le)$ with $M \triangleq \{ m_0, m_A, m_B, m_P, m_Q, m_1 \}$ illustrated in Figure 6 page 21 (B).

1. The function $\neg$ is an *ortho negation* (Definition B.3 page 35) (and thus also a *fuzzy negation*).
2. The negation $\neg m$ of each membership function $m$ is the *standard negation* (Definition 1.19 page 9).
3. $\boldsymbol{L}$ is very similar to the structure in Example 1.36 (page 18), which fails to be a logic.
4. However the structure of this example has a valid ordering relation $\le$, has valid operators ($\wedge$, $\vee$) defined in terms of $\le$, and is an *orthocomplemented lattice* (Definition A.44 page 31) with all the accompanying properties of an orthocomplemented lattice including the *non-contradiction, excluded middle* and *idempotency* properties (Theorem A.14 page 25, Definition A.44 page 31, Theorem A.47 page 32).
5. In this example, the operators are no longer *Łukasiewicz operators* (as in Example 1.36, but some other operators.
6. This orthocomplemented lattice together with any one of the *implication*s given in Example C.4 (page 46) is an *ortho logic* (and thus also a *fuzzy subset logic*—Definition 1.11 page 6).

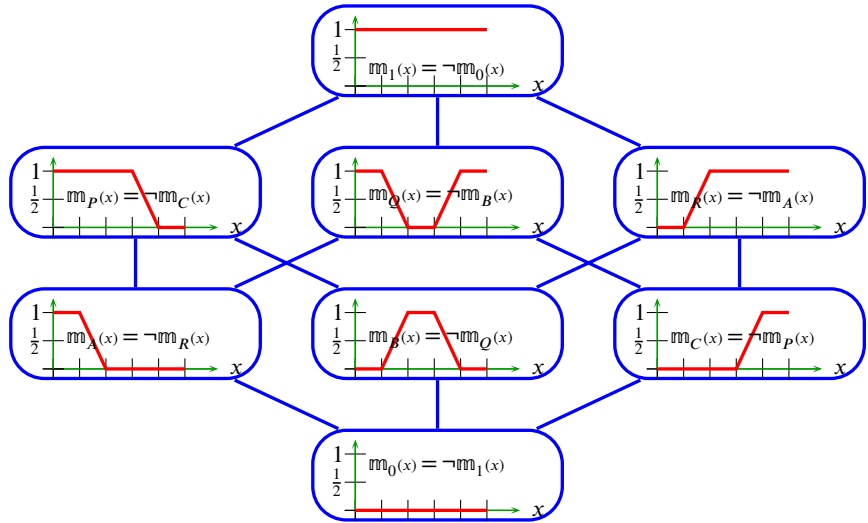

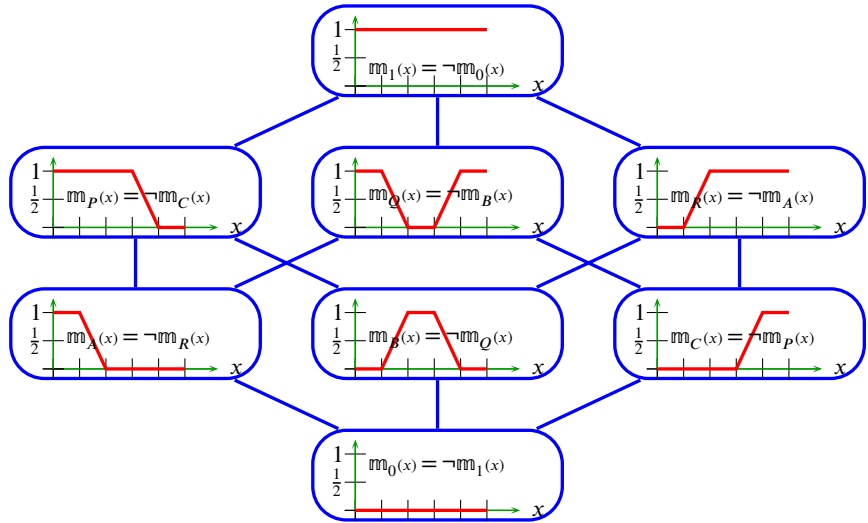

 pointwise $\leq$ induced operators on $\boldsymbol{L}_2^3$ (see Example 2.3 page 22)
 $\neg\mathbb{m}(x)$ is standard negation ($\lambda = 0$)
 **Boolean fuzzy subset logic**

Figure 7:  *fuzzy logic* on *Boolean lattice* $\boldsymbol{L}_2^3$

**Example 2.3**  Consider the structure $\boldsymbol{L} \triangleq (\mathbb{M}, \vee, \wedge, \neg, 0, 1 ; \leq)$ illustrated in Figure 7 (page 22).

1. The function $\neg$ is an *ortho negation* (Definition B.3 page 35) (and thus also a *fuzzy negation*).
2. The negation $\neg\mathbb{m}$ of each membership function $\mathbb{m}$ is the *standard negation* (Definition 1.19 page 9).
3. $\boldsymbol{L}$ is somewhat similar to the fuzzy subset logic of Example 1.37 (page 19), which fails to be *Boolean*.
4. However the structure of this example has a valid ordering relation $\leq$, has valid operators ($\wedge, \vee$) defined in terms of $\leq$, and is a *Boolean lattice* with the accompanying Boolean properties including the *non-contradiction, excluded middle, idempotent,* and *distributivity* properties (Theorem A.42 page 30).
5. In this example, the operators are no longer *min-max operators* (as in Example 1.37), but some other operators.
6. This Boolean lattice together with the *classical implication* (Example C.4 page 46) is an *ortho logic* (and thus also a *fuzzy logic*).

# Appendix A    Background: Order

## A.1    Ordered sets

**Definition A.1**   [31]   Let $2^{XX}$ be the set of all *relation*s on a set $X$.
 A relation $\leq$ is an **order relation** in $2^{XX}$ if

| | | | |
|---|---|---|---|
| 1. $x \le x$ | $\forall x \in X$ | (*reflexive*) | and |
| 2. $x \le y$ and $y \le z \implies x \le z$ | $\forall x,y,z \in X$ | (*transitive*) | and |
| 3. $x \le y$ and $y \le x \implies x = y$ | $\forall x,y \in X$ | (*anti-symmetric*) | |

} *preorder*

The pair $(X, \le)$ is an **ordered set** if $\le$ is an *order relation* on a set $X$. If $x \le y$ or $y \le x$, then elements $x$ and $y$ are said to be **comparable**, denoted $x \sim y$. Otherwise they are **incomparable**, denoted $x || y$.

**Definition A.2**  [32]  Let $(X, \le)$ be an *ordered set*. Let $2^{XX}$ be the set of all relations on $X$. The relations $\ge, <, > \in 2^{XX}$ are defined as follows:

$$x \ge y \overset{\text{def}}{\iff} y \le x \qquad \forall x,y \in X$$

$$x < y \overset{\text{def}}{\iff} x \le y \quad \text{and} \quad x \ne y \quad \forall x,y \in X$$

$$x > y \overset{\text{def}}{\iff} x \ge y \quad \text{and} \quad x \ne y \quad \forall x,y \in X$$

**Definition A.3**  [33]  An *ordered set* $(X, \le)$ (Definition A.1 page 22) is **linear**, or is a **linearly ordered set**, if
$$x \le y \quad \text{or} \quad y \le x \qquad \forall x,y \in X \qquad \text{(*comparable*)}$$
A *linearly ordered set* is also called a **totally ordered set**, a **fully ordered set**, and a **chain**.

**Definition A.4**  [34]  $y$ **covers** $x$, denoted $x \prec y$, in the ordered set $(X, \le)$ if
1. $x \le y$  (*$y$ is greater than $x$*)  and
2. $(x \le z \le y) \implies (z = x \text{ or } z = y)$  (there is no element between $x$ and $y$).

An ordered set can be represented graphically by a *Hasse diagram* (next definition).

**Definition A.5**  Let $(X, \le)$ be an ordered pair. A diagram is a **Hasse diagram** of $(X, \le)$ if
1. Each element in $X$ is represented by a dot or small circle and
2. for each $x, y \in X$, if $x \prec y$, then $y$ appears at a higher position than $x$ and a line connects $x$ and $y$.

**Example A.6**  Here are three ways of representing the ordered set $(2^{\{x,y\}}, \subseteq)$:

(1) *Hasse diagram*:  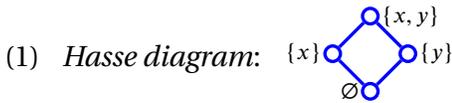

(2) Sets of ordered pairs specifying *order relation*s:
$$\subseteq = \left\{ \begin{array}{llll} (\varnothing, \varnothing), & (\{x\}, \{x\}), & (\{y\}, \{y\}), & (\{x,y\}, \{x,y\}), \\ (\varnothing, \{x\}), & (\varnothing, \{y\}), & (\varnothing, \{x,y\}), & (\{x\}, \{x,y\}), (\{y\}, \{x,y\}) \end{array} \right\}$$

(3) Sets of ordered pairs specifying *covering relation*s:
$$\prec = \left\{ (\varnothing, \{x\}), \quad (\varnothing, \{y\}), \quad (\{x\}, \{x,y\}), (\{y\}, \{x,y\}) \right\}$$

**Definition A.7**  Let $Y^X$ be the set of all functions that map from a set $X$ to a set $Y$. Let $(Y, \lessgtr)$ be an *ordered set*. The relation $\le$ is a *pointwise ordering relation* on $Y^X$ with respect to $\lessgtr$ if for all $f, g \in Y^X$
$$f \le g \implies \{f(x) \lessgtr g(x) \quad \forall x \in X\}$$

---

**Definition A.8**   Let $(X, \leq)$ be an ordered set and $2^X$ the power set of $X$.
 For any set $A \in 2^X$, $c$ is an **upper bound** of $A$ in $(X, \leq)$ if

    1.  $x \in A \implies x \leq c$ .

An element $b$ is the **least upper bound**, or **l.u.b.**, of $A$ in $(X, \leq)$ if

    2.  $b$ and $c$ are *upper bound*s of $A \implies b \leq c$ .

The *least upper bound* of the set $A$ is denoted $\bigvee A$. It is also called the **supremum** of $A$, which is denoted $\sup A$. The **join** $x \vee y$ of $x$ and $y$ is defined as $x \vee y \triangleq \bigvee \{x, y\}$.

**Definition A.9**   Let $(X, \leq)$ be an ordered set and $2^X$ the power set of $X$.   For any set $A \in 2^X$, $p$ is a **lower bound** of $A$ in $(X, \leq)$ if

    1.  $p \leq x \quad \forall x \in A$.

An element $a$ is the **greatest lower bound**, or **glb**, of $A$ in $(X, \leq)$ if

    2.  $a$ and $p$ are *lower bound*s of $A \implies p \leq a$.

The *greatest lower bound* of the set $A$ is denoted $\bigwedge A$. It is also called the **infimum** of $A$, which is denoted $\inf A$. The **meet** $x \wedge y$ of $x$ and $y$ is defined as $x \wedge y \triangleq \bigwedge \{x, y\}$.

**Proposition A.10**

$$x \leq y \iff \left\{ \begin{array}{llll} 1. & x \wedge y & = & x \quad and \\ 2. & x \vee y & = & y \end{array} \right\} \quad \forall x, y \in X$$

## A.2  Lattices

### A.2.1   General lattices

**Definition A.11**   [35]  An algebraic structure $\boldsymbol{L} \triangleq (X, \vee, \wedge ; \leq)$ is a **lattice** if

    1.  $(X, \leq)$ is an ordered set      <span style="font-size:smaller">($(X, \leq)$ is a partially or totally ordered set)</span>      and

    2.  $x, y \in X \implies \exists (x \vee y) \in X$   <span style="font-size:smaller">(every pair of elements in $X$ has a *least upper bound* in $X$)</span>      and

    3.  $x, y \in X \implies \exists (x \wedge y) \in X$   <span style="font-size:smaller">(every pair of elements in $X$ has a *greatest lower bound* in $X$).</span>

The *lattice* $\boldsymbol{L}$ is *linear* if $(X, \leq)$ is a *linearly ordered set* <span style="font-size:smaller">(Definition A.3 page 23)</span>.

**Example A.12**   [36]The *ordered set* $(X, \leq)$ illustrated by the *Hasse diagram* to the right is **not** a *lattice* because, $a$ and $b$ have no *lower bound* in $X$.

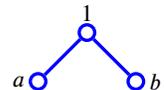

**Example A.13**   [37]  The *ordered set* illustrated by the *Hasse diagram* to the right is **not** a *lattice* because, for example, while $a$ and $b$ have *upper bound*s $c$, $d$, and $1$, still $a$ and $b$ have no *least upper bound*. The element $1$ is not the *least upper bound* because $c \leq 1$ and $d \leq 1$. And neither $c$ nor $d$ is a *least upper bound* because $c \not\leq d$ and $d \not\leq c$; rather, $c$ and $d$ are *incomparable* $(a \| b)$. Note that if we remove either or both of the two lines crossing the center, the *ordered set* becomes a *lattice*.

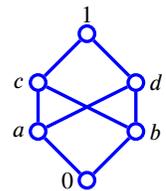

**Theorem A.14** [38]  $(X, \vee, \wedge; \le)$ *is a* LATTICE   $\Longleftrightarrow$

$$
\left\{
\begin{array}{rcl}
x \vee x & = & x \\
x \vee y & = & y \vee x \\
(x \vee y) \vee z & = & x \vee (y \vee z) \\
x \vee (x \wedge y) & = & x
\end{array}
\right.
\left|
\begin{array}{rcl}
x \wedge x & = & x \\
x \wedge y & = & y \wedge x \\
(x \wedge y) \wedge z & = & x \wedge (y \wedge z) \\
x \wedge (x \vee y) & = & x
\end{array}
\right.
\left.
\begin{array}{ll}
\forall x \in X & \text{(IDEMPOTENT)} \quad and \\
\forall x,y \in X & \text{(COMMUTATIVE)} \quad and \\
\forall x,y,z \in X & \text{(ASSOCIATIVE)} \quad and \\
\forall x,y \in X & \text{(ABSORPTIVE)}.
\end{array}
\right\}
$$

**Proposition A.15**   (Monotony laws)  [39] *Let* $(X, \vee, \wedge; \le)$ *be a* LATTICE.

$$
\left\{
\begin{array}{rcll}
a & \le & b & and \\
x & \le & y &
\end{array}
\right\}
\Longrightarrow
\left\{
\begin{array}{rcll}
a \wedge x & \le & b \wedge y & and \\
a \vee x & \le & b \vee y &
\end{array}
\right\}
\quad \forall a,b,x,y \in X
$$

**Theorem A.16**   (distributive inequalities)  [40]   $(X, \vee, \wedge; \le)$ *is a* LATTICE  $\Longrightarrow$

$$
\left\{
\begin{array}{rcl}
x \wedge (y \vee z) & \ge & (x \wedge y) \vee (x \wedge z) \\
x \vee (y \wedge z) & \le & (x \vee y) \wedge (x \vee z) \\
(x \wedge y) \vee (x \wedge z) \vee (y \wedge z) & \le & (x \vee y) \wedge (x \vee z) \wedge (y \vee z)
\end{array}
\right.
\left.
\begin{array}{ll}
\forall x,y,z \in X & \text{(JOIN SUPER-DISTRIBUTIVE)} \quad and \\
\forall x,y,z \in X & \text{(MEET SUB-DISTRIBUTIVE)} \quad and \\
\forall x,y,z \in X & \text{(MEDIAN INEQUALITY)}.
\end{array}
\right\}
$$

**Theorem A.17**   (Modular inequality)  [41]  *Let* $(X, \vee, \wedge; \le)$ *be a* LATTICE.

$$ x \le y \quad \Longrightarrow \quad x \vee (y \wedge z) \le y \wedge (x \vee z) $$

Theorem A.14 (page 25) gives 4 necessary and sufficient pairs of properties for a structure $(X, \vee, \wedge; \le)$ to be a *lattice*. However, these 4 pairs are actually *overly* sufficient (they are not *independent*), as demonstrated next.

**Theorem A.18** [42]
$(X, \vee, \wedge; \le)$ *is a* LATTICE   $\Longleftrightarrow$

$$
\left\{
\begin{array}{rcl}
x \vee y & = & y \vee x \\
(x \vee y) \vee z & = & x \vee (y \vee z) \\
x \vee (x \wedge y) & = & x
\end{array}
\right.
\left|
\begin{array}{rcl}
x \wedge y & = & y \wedge x \\
(x \wedge y) \wedge z & = & x \wedge (y \wedge z) \\
x \wedge (x \vee y) & = & x
\end{array}
\right.
\left.
\begin{array}{ll}
\forall x,y \in X & \text{(COMMUTATIVE)} \quad and \\
\forall x,y,z \in X & \text{(ASSOCIATIVE)} \quad and \\
\forall x,y \in X & \text{(ABSORPTIVE)}
\end{array}
\right\}
$$

### A.2.2   Bounded *lattice*s

**Definition A.19**   Let $L \triangleq (X, \vee, \wedge; \le)$ be a *lattice*. Let $\bigvee X$ be the least upper bound of $(X, \le)$ and let $\bigwedge X$ be the greatest lower bound of $(X, \le)$.   $L$ is **upper bounded** if $(\bigvee X) \in X$. $L$ is **lower bounded** if $(\bigwedge X) \in X$. $L$ is **bounded** if $L$ is both upper and lower bounded. A *bounded lattice* is optionally denoted $(X, \vee, \wedge, 0, 1; \le)$, where $0 \triangleq \bigwedge X$ and $1 \triangleq \bigvee X$.

**Proposition A.20**   *Let* $L \triangleq (X, \vee, \wedge; \le)$ *be a* LATTICE.
$\{L \text{ is FINITE}\} \quad \Longrightarrow \quad \{L \text{ is BOUNDED}\}$

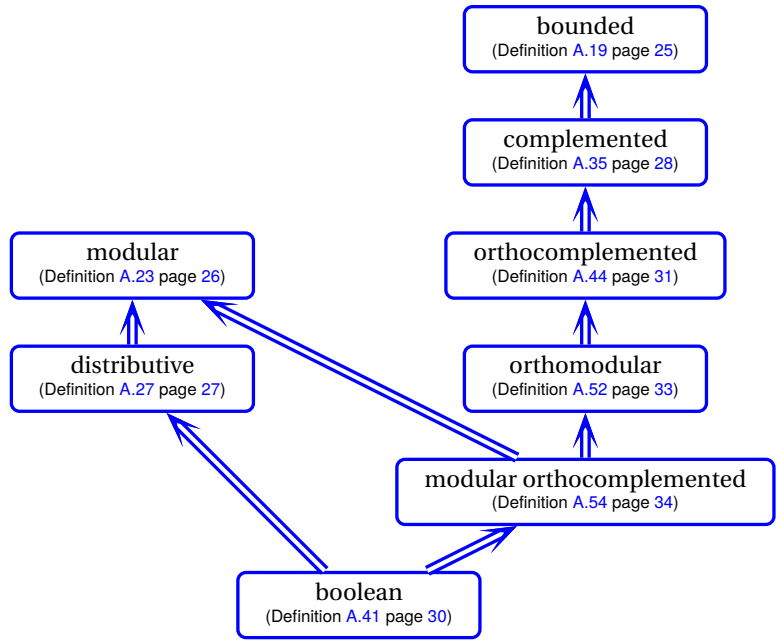

Figure 8: relationships between selected lattice types

**Proposition A.21**  *Let* $L \triangleq (X, \vee, \wedge; \leq)$ *be a* LATTICE *with* $\bigvee X \triangleq 1$ *and* $\bigwedge X \triangleq 0$.

$$\{L \text{ is BOUNDED}\} \implies \left\{ \begin{array}{llll} x \vee 1 & = & 1 & \forall_{x \in X} & \text{(upper bounded)} & \text{and} \\ x \wedge 0 & = & 0 & \forall_{x \in X} & \text{(lower bounded)} & \text{and} \\ x \vee 0 & = & x & \forall_{x \in X} & \text{(join-identity)} & \text{and} \\ x \wedge 1 & = & x & \forall_{x \in X} & \text{(meet-identity)} \end{array} \right\}$$

### A.2.3  Modular lattices

**Definition A.22**  [43] Let $(X, \vee, \wedge; \leq)$ be a lattice. The **modularity** relation $\circledR \in 2^{XX}$ is defined as

$$x \circledR y \quad \overset{\text{def}}{\iff} \quad \left\{ (x, y) \in X^2 \, | a \leq y \quad \implies \quad y \wedge (x \vee a) = (y \wedge x) \vee a \quad \forall a \in X \right\}.$$

Modular lattices are a generalization of *distributive lattice*s in that all distributive lattices are modular, but not all modular lattices are distributive (Example A.33 page 28, Example A.34 page 28).

**Definition A.23**  [44]  A lattice $(X, \vee, \wedge; \leq)$ is **modular** if $\quad x \circledR y \quad \forall x, y \in X$.

**Definition A.24**  (N5 lattice/pentagon)  [45]  The **N5 lattice** is the ordered set $(\{0, a, b, p, 1\}, \leq)$ with cover relation

$\quad \prec = \{(0, a), (a, b), (b, 1), (p, 1), (0, p)\}$.

The N5 lattice is also called the **pentagon**.  The N5 lattice is illustrated by the Hasse diagram to the right.

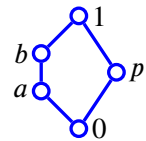

---

[43]   [Stern(1999)] page 11,   [Maeda and Maeda(1970)], page 1 ⟨Definition (1.1)⟩,  ☞ [Maeda(1966)] page 248

[44]  ☞ [Birkhoff(1967)] page 82,  ☞ [Maeda and Maeda(1970)], page 3 ⟨Definition (1.7)⟩

[45]  [Beran(1985)] pages 12–13,  📖 [Dedekind(1900)] pages 391–392 ⟨(44) and (45)⟩

**Theorem A.25** [46] *Let* **L** *be a* LATTICE *(Definition A.11 page 24).*

    **L** *is* MODULAR  $\iff$  **L** *does* NOT *contain the* N5 LATTICE

Examples of *modular lattice*s are provided in Example A.33 (page 28) and Example A.34 (page 28).

### A.2.4   Distributive lattices

**Definition A.26** [47] Let $(X, \vee, \wedge; \leq)$ be a *lattice* (Definition A.11 page 24).
The **distributivity relation** $\circledD \in 2^{XXX}$ and the **dual distributivity relation** $\circledD^* \in 2^{XXX}$ are defined as

    $\circledD \triangleq \left\{ (x, y, z) \in X^3 \,|\, x \wedge (y \vee z) = (x \wedge y) \vee (x \wedge z) \right\}$  (each $(x, y, z)$ is *disjunctive distributive*).

    $\circledD^* \triangleq \left\{ (x, y, z) \in X^3 \,|\, x \vee (y \wedge z) = (x \vee y) \wedge (x \vee z) \right\}$  (each $(x, y, z)$ is *conjunctive distributive*).

A triple $(x, y, z)$ is **distributive** if $(x, y, z) \in \circledD$ and such a triple is alternatively denoted as $(x, y, z)\circledD$.

**Definition A.27** [48]  A lattice $(X, \vee, \wedge; \leq)$ is **distributive** if  $(x, y, z) \in \circledD$  $\forall x, y, z \in X$

Not all lattices are *distributive.* But if a lattice **L** does happen to be distributive—that is all triples in **L** satisfy the *distributive* property—then all triples in **L** also satisfy the *dual distributive* property, as well as another property called the *median property*. The converses also hold (next theorem).

**Theorem A.28** [49] *Let* $\mathbf{L} \triangleq (X, \vee, \wedge; \leq)$ *be a* LATTICE. *The following statements are all equivalent:*

| | | | |
|---|---|---|---|
| *(1).* | **L** *is* DISTRIBUTIVE | | *(Definition A.27 page 27)* |
| $\iff$ *(2).* | $x \wedge (y \vee z) = (x \wedge y) \vee (x \wedge z)$ | $\forall x,y,z \in X$ | (DISJUNCTIVE DISTRIBUTIVE) |
| $\iff$ *(3).* | $x \vee (y \wedge z) = (x \vee y) \wedge (x \vee z)$ | $\forall x,y,z \in X$ | (CONJUNCTIVE DISTRIBUTIVE) |
| $\iff$ *(4).* | $(x \vee y) \wedge (x \vee z) \wedge (y \vee z) = (x \wedge y) \vee (x \wedge z) \vee (y \wedge z)$ | $\forall x,y,z \in X$ | (MEDIAN PROPERTY) |

**Definition A.29** (M3 lattice/diamond) [50] The **M3 lattice** is the ordered set $(\{0, p, q, r, 1\}, \leq)$ with covering relation

    $\lessdot = \{(p, 1), (q, 1), (r, 1), (0, p), (0, q), (0, r)\}$.

The M3 lattice is also called the **diamond**, and is illustrated by the Hasse diagram to the right.

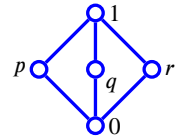

**Theorem A.30** (Birkhoff distributivity criterion) [51] *Let* $\mathbf{L} \triangleq (X, \vee, \wedge; \leq)$ *be a* LATTICE.

    **L** *is* DISTRIBUTIVE  $\iff$  $\begin{cases} \textbf{L} \textit{ does } \textbf{not} \textit{ contain N5 as a sublattice} & \text{and} \\ \textbf{L} \textit{ does } \textbf{not} \textit{ contain M3 as a sublattice} \end{cases}$

Distributive lattices are a special case of modular lattices. That is, all distributive lattices are modular, but not all modular lattices are distributive (next theorem). An example is the *M3 lattice*—it is modular, but yet it is not *distributive*.

**Theorem A.31** [52] *Let* $( X, \vee, \wedge; \leq )$ *be a lattice.*

$( X, \vee, \wedge; \leq )$ *is* DISTRIBUTIVE $\quad\overset{\Longrightarrow}{\Longleftarrow}\quad$ $( X, \vee, \wedge; \leq )$ *is* MODULAR.

**Proposition A.32** [53] *Let* $X_n$ *be a finite set with order* $n = \big| X_n \big|$. *Let* $l_n$ *be the number of unlabeled lattices on* $X_n$, $m_n$ *the number of unlabeled modular lattices on* $X_n$. *and* $d_n$ *the number of unlabeled distributive lattices on* $X_n$.

| $n$ | 0 | 1 | 2 | 3 | 4 | 5 | 6 | 7 | 8 | 9 | 10 | 11 | 12 | 13 | 14 |
|---|---|---|---|---|---|---|---|---|---|---|---|---|---|---|---|
| $l_n$ | 1 | 1 | 1 | 1 | 2 | 5 | 15 | 53 | 222 | 1078 | 5994 | 37622 | 262776 | 2018305 | 16873364 |
| $m_n$ | 1 | 1 | 1 | 1 | 2 | 4 | 8 | 16 | 34 | 72 | 157 | 343 | 766 | 1718 | 3899 |
| $d_n$ | 1 | 1 | 1 | 1 | 2 | 3 | 5 | 8 | 15 | 26 | 47 | 82 | 151 | 269 | 494 |

**Example A.33** [54] There are a total of 5 unlabeled lattices on a five element set. Of these, 3 are *distributive* (Proposition A.32 page 28, and thus also *modular*), one is *modular* but *non-distributive*, and one is *non-distributive* (and *non-modular*).

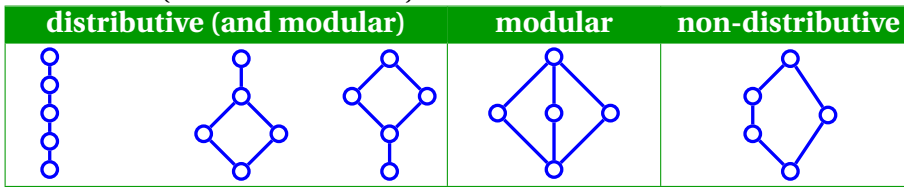

**Example A.34** [55] There are a total of 15 unlabeled lattices on a six element set; and of these 15, five are distributive (Proposition A.32 page 28). The following illustrates the 5 distributive lattices. Note that none of these lattices are *complemented* (none are *Boolean* Definition A.41 page 30).

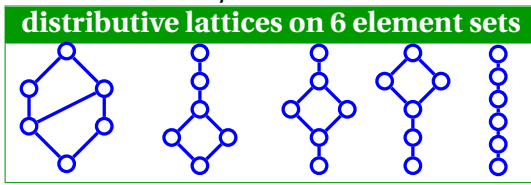

## A.3  Complemented lattices

### A.3.1  Definitions

**Definition A.35** [56] Let $L \triangleq ( X, \vee, \wedge, 0, 1; \leq )$ be a *bounded lattice* (Definition A.19 page 25). An element $x' \in X$ is a **complement** of an element $x$ in $L$ if

---

1. $x \wedge x' = 0$    (*non-contradiction*)    and
2. $x \vee x' = 1$    (*excluded middle*).

An element $x'$ in $\mathbf{L}$ is the *unique complement* of $x$ in $\mathbf{L}$ if $x'$ is a *complement* of $x$ and $y'$ is a *complement* of $x \implies x' = y'$. $\mathbf{L}$ is **complemented** if every element in $X$ has a complement in $X$. $\mathbf{L}$ is **uniquely complemented** if every element in $X$ has a unique complement in $X$. A complemented lattice that is *not* uniquely complemented is **multiply complemented**.

**Example A.36**   Here are some examples:

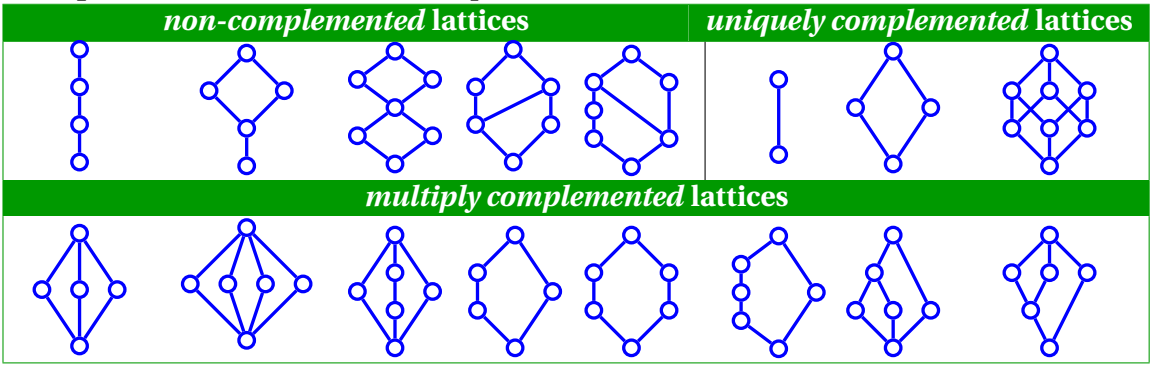

**Example A.37**   Of the 53 unlabeled lattices on a 7 element set, 0 are *uniquely complemented*, 17 are *multiply complemented*, and 36 are *non-complemented*.

Theorem A.38 (next) is a landmark theorem in mathematics.

**Theorem A.38**  [57]  *For every lattice $\mathbf{L}$, there exists a lattice $\mathbf{U}$ such that*
1. $\mathbf{L} \subseteq \mathbf{U}$ *($\mathbf{L}$ is a sublattice of $\mathbf{U}$)*    and
2. $\mathbf{U}$ *is* UNIQUELY COMPLEMENTED.

**Corollary A.39**  [58]  *Let $\mathbf{L} \triangleq (X, \vee, \wedge; \leq)$ be a lattice.*

$$\left\{ \begin{array}{l} \text{1. } \mathbf{L} \text{ is DISTRIBUTIVE} \quad and \\ \text{2. } \mathbf{L} \text{ is COMPLEMENTED} \end{array} \right\} \quad \overset{\implies}{\Longleftarrow} \quad \{ \mathbf{L} \text{ is UNIQUELY COMPLEMENTED} \}$$

**Theorem A.40**   (Huntington properties)  [59]  *Let $\mathbf{L}$ be a lattice.*

$$\left\{ \begin{array}{l} \mathbf{L} \text{ is} \\ \text{UNIQUELY} \\ \text{COMPLEMENTED} \end{array} \right\} \text{ and } \underbrace{\left\{ \begin{array}{ll} \mathbf{L} \text{ is MODULAR} & or \\ \mathbf{L} \text{ is ATOMIC} & or \\ \mathbf{L} \text{ is ORTHOCOMPLEMENTED} & or \\ \mathbf{L} \text{ has FINITE WIDTH} & or \\ \mathbf{L} \text{ is DE MORGAN} \end{array} \right\}}_{\text{HUNTINGTON PROPERTIES}} \implies \left\{ \begin{array}{l} \mathbf{L} \text{ is} \\ \text{DISTRIBUTIVE} \end{array} \right\}$$

### A.3.2  Boolean lattices

**Definition A.41** [60]  A *lattice* (Definition A.11 page 24) $L$ is **Boolean** if
1.  $L$ is *bounded*        (Definition A.19 page 25)    and
2.  $L$ is *distributive*        (Definition A.27 page 27)    and
3.  $L$ is *complemented*        (Definition A.35 page 28)    .

In this case, $L$ is a **Boolean algebra** or a **Boolean lattice**. In this paper, a *Boolean lattice* is denoted $( X, \vee, \wedge, 0, 1 \,; \leq)$, and a *Boolean lattice* with $2^N$ elements is sometimes denoted $L_2^N$.

**Theorem A.42**  (classic 10 Boolean properties) [61]  *Let* $A \triangleq ( X, \vee, \wedge, 0, 1 \,; \leq)$ *be an algebraic structure. In the event that* $A$ *is a* BOUNDED LATTICE (Definition A.19 page 25), *let* $x'$ *represent a* COMPLEMENT (Definition A.35 page 28) *of an element* $x$ *in* $A$.

$A$ *is a Boolean algebra*    $\implies$    $\forall x, y, z \in X$

| disjunctive properties | | conjunctive properties | | property name | |
|---|---|---|---|---|---|
| $x \vee x$ | $= x$ | $x \wedge x$ | $= x$ | (IDEMPOTENT) | *and* |
| $x \vee y$ | $= y \vee x$ | $x \wedge y$ | $= y \wedge x$ | (COMMUTATIVE) | *and* |
| $x \vee (y \vee z)$ | $= (x \vee y) \vee z$ | $x \wedge (y \wedge z)$ | $= (x \wedge y) \wedge z$ | (ASSOCIATIVE) | *and* |
| $x \vee (x \wedge y)$ | $= x$ | $x \wedge (x \vee y)$ | $= x$ | (ABSORPTIVE) | *and* |
| $x \vee 1$ | $= 1$ | $x \wedge 0$ | $= 0$ | (BOUNDED) | *and* |
| $x \vee 0$ | $= x$ | $x \wedge 1$ | $= x$ | (IDENTITY) | *and* |
| $x \vee (y \wedge z)$ | $= (x \vee y) \wedge (x \vee z)$ | $x \wedge (y \vee z)$ | $= (x \wedge y) \vee (x \wedge z)$ | (DISTRIBUTIVE) | *and* |
| $x \vee x'$ | $= 1$ | $x \wedge x'$ | $= 0$ | (COMPLEMENTED) | *and* |
| $(x \vee y)'$ | $= x' \wedge y'$ | $(x \wedge y)'$ | $= x' \vee y'$ | (DE MORGAN) | *and* |
| | | $(x')' = x$ | | (INVOLUTORY) | |

**Lemma A.43**
$\left.\begin{array}{l}( X, \vee, \wedge, 0, 1 \,; \leq) \\ \textit{is a } \textsc{Boolean algebra}\end{array}\right\} \implies \left\{\begin{array}{ll}\textit{1. } & x' \vee (x \wedge y) = x' \vee y \quad \forall x, y \in X \quad (\textsc{Sasaki hook}) \quad \textit{and} \\ \textit{2. } & x \vee (x' \wedge y) = x \vee y \quad \forall x, y \in X\end{array}\right.$

PROOF:

$x' \vee (x \wedge y) = \underbrace{x' \vee (x' \wedge y)}_{x'} \vee (x \wedge y)$    by *absorption* property (Theorem A.42 page 30)

$= x' \vee \big[(x' \vee x) \wedge y\big]$    by *associative* and *distributive* properties (Theorem A.42 page 30)

$= x' \vee [1 \wedge y]$    by *excluded middle* property (Theorem A.42 page 30)

$= x' \vee y$    by definition of 1 (Definition A.8 page 24)

$x \vee (x' \wedge y) = \underbrace{x \vee (x \wedge y)}_{x} \vee (x \wedge y)$    by *absorption* property (Theorem A.42 page 30)

$= x \vee \big[(x \vee x') \wedge y\big]$    by *associative* and *distributive* properties (Theorem A.42 page 30)

$= x \vee [1 \wedge y]$    by *excluded middle* property (Theorem A.42 page 30)

$= x \vee y$    by definition of 1 (Definition A.8 page 24)

---

### A.3.3  Orthocomplemented Lattices

**Definition A.44**  [62]  Let $L \triangleq ( X, \vee, \wedge, 0, 1 ; \le )$ be a *bounded lattice* (Definition A.19 page 25).   An element $x^{\perp} \in X$ is an **orthocomplement** of an element $x \in X$ if

1. $x^{\perp\perp} = x$           $\forall x \in X$   (*involutory*)   and
2. $x \wedge x^{\perp} = 0$           $\forall x \in X$   (*non-contradiction*)   and
3. $x \le y \implies y^{\perp} \le x^{\perp}$   $\forall x, y \in X$   (*antitone*).

The lattice $L$ is **orthocomplemented** ($L$ is an **orthocomplemented lattice**) if every element $x$ in $X$ has an *orthocomplement*.

**Definition A.45**  [63]  The $\mathbf{O_6}$ **lattice** is the ordered set $\left( \{ 0, p, q, p^{\perp}, q^{\perp}, 1 \}, \le \right)$ with cover relation
$$\lessdot = \left\{ (0, p), (0, q), \left( p, q^{\perp} \right), \left( q, p^{\perp} \right), \left( p^{\perp}, 1 \right), \left( q^{\perp}, 1 \right) \right\}.$$
The $O_6$ lattice is illustrated by the Hasse diagram to the right.

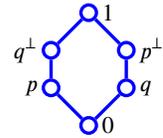

**Example A.46**  [64]  There are a total of 10 **orthocomplemented lattices** with 8 elements or less.  These along with some other orthocomplemented lattices are illustrated next:[65]

Lattices that are **orthocomplemented** but *non-orthomodular* and hence also *non-modular-orthocomplemented* and *non-Boolean*:

1.  $O_6$ *lattice*          2.  $O_8$ *lattice*          3.          4.

5.          6.          7.

Lattices that are **orthocomplemented** and **orthomodular** but *non-modular-orthocomplemented* and hence also *non-Boolean*:

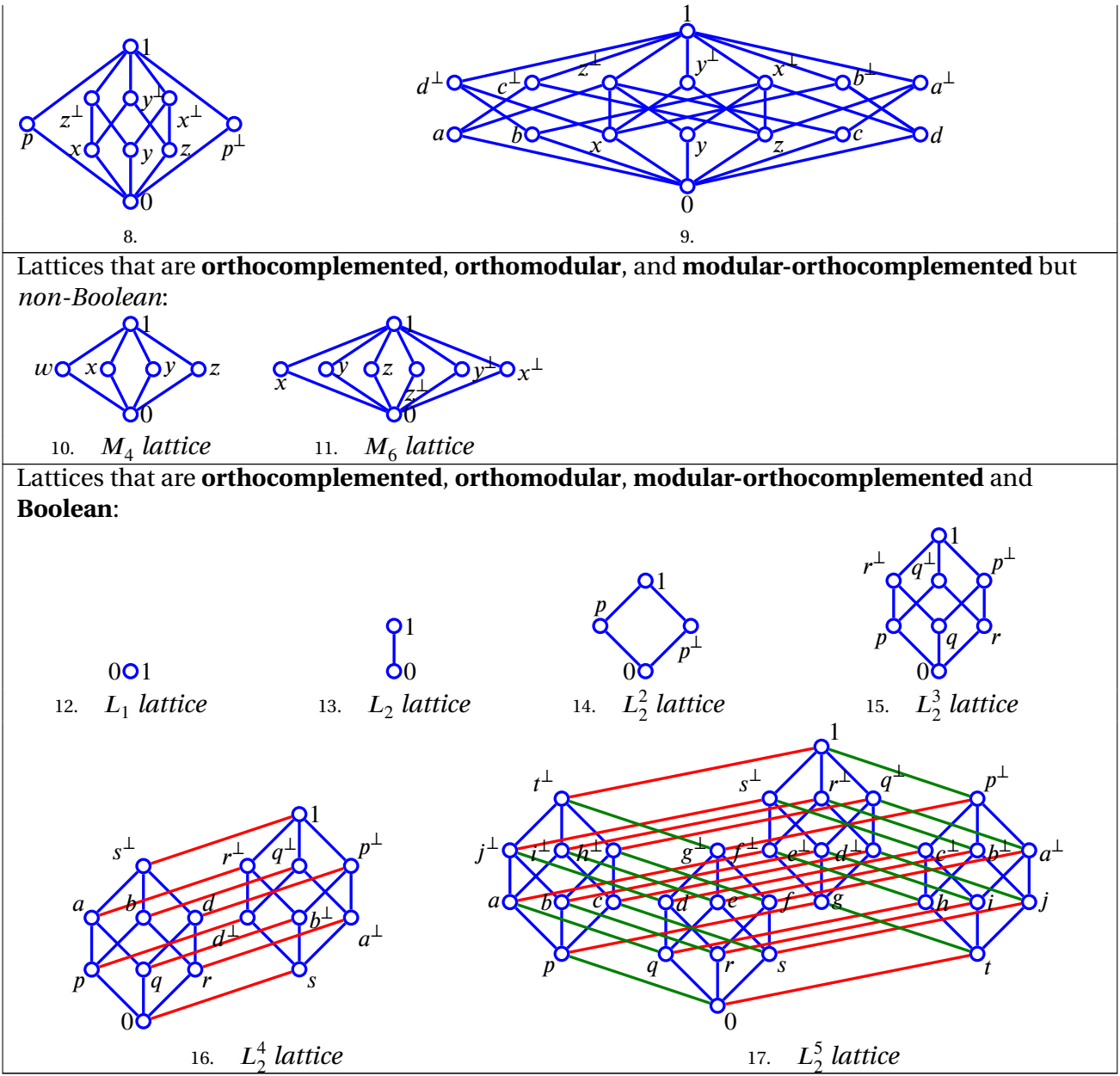

8.                                                9.

Lattices that are **orthocomplemented**, **orthomodular**, and **modular-orthocomplemented** but *non-Boolean*:

10.   $M_4$ *lattice*          11.   $M_6$ *lattice*

Lattices that are **orthocomplemented**, **orthomodular**, **modular-orthocomplemented** and **Boolean**:

12.   $L_1$ *lattice*          13.   $L_2$ *lattice*          14.   $L_2^2$ *lattice*          15.   $L_2^3$ *lattice*

16.   $L_2^4$ *lattice*                          17.   $L_2^5$ *lattice*

**Theorem A.47**  [66]  *Let $x^\perp$ be the* ORTHOCOMPLEMENT *of an element $x$ in a* BOUNDED LATTICE $\boldsymbol{L} \triangleq (X, \vee, \wedge, 0, 1 ; \leq)$.

$$\left\{ \begin{array}{l} \boldsymbol{L} \text{ is} \\ \text{ORTHOCOMPLEMENTED} \end{array} \right\} \implies \left\{ \begin{array}{llll} (1). & 0^\perp = 1 & & \text{(BOUNDARY CONDITION)} \quad \text{and} \\ (2). & 1^\perp = 0 & & \text{(BOUNDARY CONDITION)} \quad \text{and} \\ (3). & (x \vee y)^\perp = x^\perp \wedge y^\perp & \forall x,y \in X & \text{(DISJUNCTIVE DE MORGAN)} \quad \text{and} \\ (4). & (x \wedge y)^\perp = x^\perp \vee y^\perp & \forall x,y \in X & \text{(CONJUNCTIVE DE MORGAN)} \quad \text{and} \\ (5). & x \vee x^\perp = 1 & \forall x \in X & \text{(EXCLUDED MIDDLE)}. \end{array} \right\}$$

✎PROOF:      Let $x^\perp \triangleq \neg x$, where $\neg$ is an *ortho negation* function (Definition B.3 page 35). Then this theorem follows

───────────────

[66] ☞ [Beran(1985)] pages 30–31, 📖 [Birkhoff and Neumann(1936)] page 830 ⟨L74⟩, ☞ [Cohen(1989)] page 37
⟨3B.13. Theorem⟩

directly from Theorem B.15 (page 37).

**Corollary A.48**  *Let* $\mathbf{L} \triangleq ( X, \vee, \wedge, 0, 1 ; \leq )$ *be a* LATTICE *(Definition A.11 page 24).*

$$\left\{ \begin{array}{c} \mathbf{L} \text{ \textit{is} \textbf{orthocomplemented}} \\ \textit{(Definition A.44 page 31)} \end{array} \right\} \implies \left\{ \begin{array}{c} \mathbf{L} \text{ \textit{is} \textbf{complemented}} \\ \textit{(Definition A.35 page 28)} \end{array} \right\}$$

✎PROOF:   This follows directly from the definition of *orthocomplemented lattice*s (Definition A.44 page 31) and *complemented lattice*s (Definition A.35 page 28).

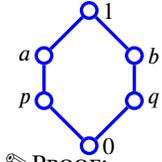

**Example A.49**  The $O_6$ *lattice* (Definition A.45 page 31) illustrated to the left is both **orthocomplemented** (Definition A.44 page 31) and **multiply complemented** (Definition A.35 page 28). The lattice illustrated to the right is **multiply complemented**, but is **non-orthocomplemented**.

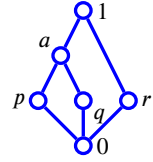

✎PROOF:

(1) Proof that $O_6$ *lattice* is multiply complemented: $b$ and $q$ are both *complements* of $p$.

(2) Proof that the right side lattice is multiply complemented: $a$, $p$, and $q$ are all *complements* of $r$.

**Proposition A.50**  [67] *Let* $\mathbf{L} = ( X, \vee, \wedge, 0, 1 ; \leq )$ *be a* BOUNDED LATTICE *(Definition A.19 page 25).*

$$\left\{ \begin{array}{ll} \textit{1.} & \mathbf{L} \text{ \textit{is} \textbf{orthocomplemented}} \quad \textit{(Definition A.44 page 31)} \quad \textit{and} \\ \textit{2.} & \mathbf{L} \text{ \textit{is} \textbf{distributive}} \qquad\qquad \textit{(Definition A.27 page 27)} \end{array} \right\} \implies \left\{ \begin{array}{c} \mathbf{L} \text{ \textit{is} \textbf{Boolean}} \\ \textit{(Definition A.41 page 30)} \end{array} \right\}$$

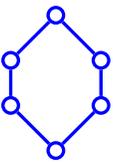

**Example A.51**  The $O_6$ *lattice* (Definition A.45 page 31) illustrated to the left is **orthocomplemented** (Definition A.44 page 31) but **non-join-distributive** (Definition A.27 page 27),and hence *non-Boolean*. The lattice illustrated to the right is **orthocomplemented** *and* **distributive** and hence also **Boolean** (Proposition A.50 page 33).

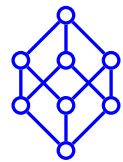

### A.3.4  Orthomodular lattices

**Definition A.52**  [68] Let $\mathbf{L} \triangleq ( X, \vee, \wedge, 0, 1 ; \leq )$ be a *bounded lattice* (Definition A.19 page 25).
$\mathbf{L}$ is **orthomodular** if

1. $\mathbf{L}$ is *orthocomplemented*                    and
2. $x \leq y \implies x \vee \left( x^{\perp} \wedge y \right) = y \quad \forall x,y \in X \quad$ *(orthomodular identity)*

**Theorem A.53**  [69] *Let* $\mathbf{L} = ( X, \vee, \wedge, 0, 1 ; \leq )$ *be an algebraic structure.*

$$\left\{ \begin{array}{c} \mathbf{L} \text{ \textit{is an} \textbf{\textit{orthomodular lattice}}} \quad \textit{and} \\ \underbrace{\left( x \wedge y^{\perp} \right)^{\perp} = y \vee \left( x^{\perp} \wedge y^{\perp} \right)}_{\text{ELKAN'S LAW}} \quad \forall x, y \in X \end{array} \right\} \implies \left\{ \begin{array}{c} \mathbf{L} \text{ \textit{is a}} \\ \textbf{\textit{Boolean algebra}} \\ \textit{(Definition A.41 page 30)} \end{array} \right\}$$

---

**Definition A.54**  Let $L \triangleq (\,X,\,\vee,\,\wedge,\,0,\,1\,;\,\leq)$ be a *bounded lattice* (Definition A.19 page 25).
$L$ is a **modular orthocomplemeted lattice** if
  1. $L$ is **orthocomplemented**  (Definition A.44 page 31)  and
  2. $L$ is **modular**  (Definition A.23 page 26)

**Theorem A.55**  [70] *Let* $L$ *be a lattice.*
$\{L$ *is* BOOLEAN$\}$  $\implies$  $\{L$ *is* MODULAR ORTHOCOMPLEMENTED  *(Definition A.54 page 34)*$\}$
  $\implies$  $\{L$ *is* ORTHOMODULAR  *(Definition A.52 page 33)*$\}$
  $\implies$  $\{L$ *is* ORTHOCOMPLEMENTED  *(Definition A.44 page 31)*$\}$

# Appendix B    Background: Negation

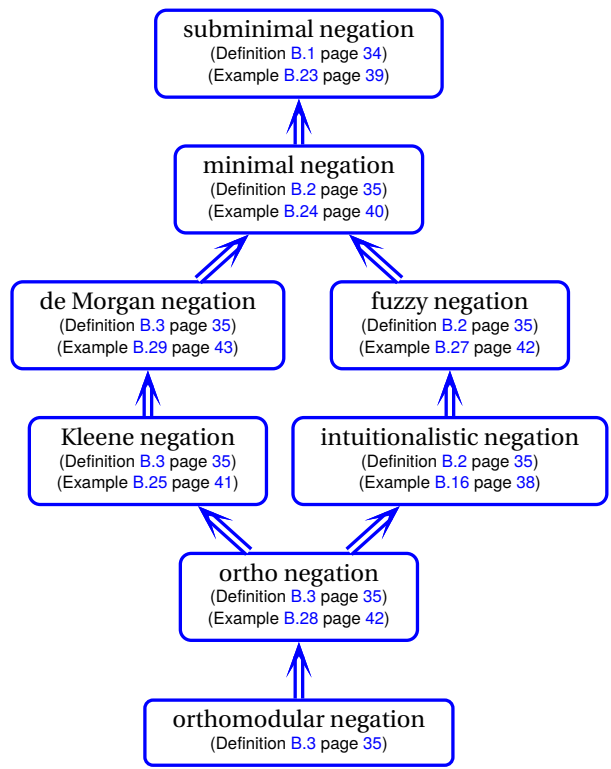

Figure 9: lattice of negations

## B.1   Definitions

**Definition B.1**  [71] Let $L \triangleq (\,X,\,\vee,\,\wedge,\,0,\,1\,;\,\leq)$ be a *bounded lattice* (Definition A.19 page 25).
  A *function* $\neg \in X^X$ is a **subminimal negation** on $L$ if [72]
    $x \leq y \implies \neg y \leq \neg x \quad \forall x,y \in X \quad$ *(antitone)*.

---

**Definition B.2**  [73]  Let $\boldsymbol{L} \triangleq (X, \vee, \wedge, 0, 1 ; \leq)$ be a *bounded lattice* (Definition A.19 page 25).
A *function* $\neg \in X^X$ is a **negation**, or **minimal negation**, on $\boldsymbol{L}$ if

1. $x \leq y \implies \neg y \leq \neg x$   $\forall x,y \in X$   *(antitone)*          and
2. $x \leq \neg\neg x$   $\forall x \in X$   *(weak double negation)*.

A *minimal negation* $\neg$ is an **intuitionistic negation** on $\boldsymbol{L}$ if

3. $x \wedge \neg x = 0$   $\forall x,y \in X$   *(non-contradiction)*.

A *minimal negation* $\neg$ is a **fuzzy negation** on $\boldsymbol{L}$ if

4. $\neg 1 = 0$   *(boundary condition)*.

**Definition B.3**  [74]  Let $\boldsymbol{L} \triangleq (X, \vee, \wedge, 0, 1 ; \leq)$ be a *bounded lattice* (Definition A.19 page 25).
A *minimal negation* $\neg$ is a **de Morgan negation** on $\boldsymbol{L}$ if

5. $x = \neg\neg x$   $\forall x \in X$   *(involutory)*.

A *de Morgan negation* $\neg$ is a **Kleene negation** on $\boldsymbol{L}$ if

6. $x \wedge \neg x \leq y \vee \neg y$   $\forall x,y \in X$   *(Kleene condition)*.

A *de Morgan negation* $\neg$ is an **ortho negation** on $\boldsymbol{L}$ if

7. $x \wedge \neg x = 0$   $\forall x,y \in X$   *(non-contradiction)*.

A *de Morgan negation* $\neg$ is an **orthomodular negation** on $\boldsymbol{L}$ if

8. $x \wedge \neg x = 0$   $\forall x,y \in X$   *(non-contradiction)*   and
9. $x \leq y \implies x \vee \left(x^{\perp} \wedge y\right) = y$   $\forall x,y \in X$   *(orthomodular)*.

**Remark B.4**  [75]  The *Kleene condition* is a weakened form of the *non-contradiction* and *excluded middle* properties in the sense $\underbrace{x \wedge \neg x = 0}_{non\text{-}contradiction} \leq \underbrace{1 = y \vee \neg y}_{excluded\ middle}$.

**Definition B.5**  Let $\boldsymbol{L} \triangleq (X, \vee, \wedge, \neg, 0, 1 ; \leq)$ be a *bounded lattice* (Definition A.19 page 25) with a function $\neg \in X^X$.  If $\neg$ is a *negation* (Definition B.2 page 35), then $\boldsymbol{L}$ is a **lattice with negation**.

## B.2  Properties of negations

**Lemma B.6**  [76]  *Let* $\neg \in X^X$ *be a function on a* LATTICE $\boldsymbol{L} \triangleq (X, \vee, \wedge ; \leq)$ *(Definition A.11 page 24).*

$\left.\underbrace{\begin{aligned}x \leq y \implies \\ \neg y \leq \neg x\end{aligned}}_{\text{ANTITONE}}\right\} \implies \begin{cases} \neg x \vee \neg y \leq \neg(x \wedge y) & \forall x,y \in X & \text{(CONJUNCTIVE DE MORGAN INEQUALITY)} & and \\ \neg(x \vee y) \leq \neg x \wedge \neg y & \forall x,y \in X & \text{(DISJUNCTIVE DE MORGAN INEQUALITY)} \end{cases}$

---

[72] In the context of natural language, D. Devidi has argued that, *subminimal negation* (Definition B.1 page 34) is "difficult to take seriously as" a negation.  For further details see 📖 [Devidi(2010)], page 511, 📖 [Devidi(2006)], page 568

[73] 📖 [Dunn(1996)] pages 4–6, 📖 [Dunn(1999)] pages 24–26 ⟨2 THE KITE OF NEGATIONS⟩, 📖 [TROELSTRA AND VAN DALEN(1988)] PAGE 4 ⟨1.6 Intuitionism. (B)⟩, 📖 [DE VRIES(2007)] page 11 ⟨Definition 16⟩, 📖 [GOTTWALD(1999)] PAGE 21 ⟨DEFINITION 3.3⟩, 📖 [NOVÁK ET AL.(1999)Novák, Perfilieva, and Močkoř] PAGE 50 ⟨Definition 2.26⟩, 📖 [NGUYEN AND WALKER(2006)] PAGES 98–99 ⟨5.4 NEGATIONS⟩, 📖 [BELLMAN AND GIERTZ(1973)] PAGES 155–156 ⟨(N1) ¬0 = 1 AND ¬1 = 0, (N3) ¬¬x = x⟩

[74] 📖 [Dunn(1999)] pages 24–26 ⟨2 THE KITE OF NEGATIONS⟩, 📖 [JENEI(2003)] PAGE 283, 📖 [KALMBACH(1983)] PAGE 22, 📖 [LIDL AND PILZ(1998)] PAGE 90, 📖 [HUSIMI(1937)]

[75] 📖 [Cattaneo and Ciucci(2009)] page 78

[76] 📖 [Beran(1985)] page 31 ⟨Theorem 1.2 Proof⟩, 📖 [Fáy(1967)] page 268 ⟨Lemma 1 Proof⟩, 📖 [de Vries(2007)] page 12 ⟨Theorem 18⟩

**Lemma B.7** [77] *Let* $\neg \in X^X$ *be a function on a* LATTICE $\boldsymbol{L} \triangleq ( X, \vee, \wedge\,; \leq)$ *(Definition A.11 page 24).*
*If* $x = (\neg\neg x)$ *for all* $x \in X$ (INVOLUTORY), *then*

$$\underbrace{x \leq y \implies \neg y \leq \neg x}_{\text{ANTITONE}} \implies \underbrace{\left\{ \begin{array}{lll} \neg(x \vee y) & = & \neg x \wedge \neg y \quad \forall x, y \in X \quad \text{(DISJUNCTIVE DE MORGAN)} \quad and \\ \neg(x \wedge y) & = & \neg x \vee \neg y \quad \forall x, y \in X \quad \text{(CONJUNCTIVE DE MORGAN)} \end{array} \right.}_{\text{DE MORGAN}}$$

**Lemma B.8** *Let* $\neg \in X^X$ *be a function on a* BOUNDED LATTICE $\boldsymbol{L} \triangleq ( X, \vee, \wedge, 0, 1\,; \leq)$.

$$\left\{ \begin{array}{ll} 1. & x \leq \neg\neg x \quad \forall x \in X \quad \text{(WEAK DOUBLE NEGATION)} \quad and \\ 2. & \neg 1 = 0 \qquad\qquad\quad \text{(BOUNDARY CONDITION)} \end{array} \right\} \implies \left\{ \; \neg 0 \; = \; 1 \quad \text{(BOUNDARY CONDITION)} \; \right\}$$

PROOF:

$$
\begin{array}{ll}
\neg 0 = \neg\neg 1 & \text{by } \textit{boundary condition} \text{ hypothesis (2)} \\
\qquad \geq 1 & \text{by } \textit{weak double negation} \text{ hypothesis (1)} \\
\implies \neg 0 = 1 & \text{by } \textit{upper bound} \text{ property (Definition A.19 page 25)}
\end{array}
$$

**Lemma B.9** *Let* $\neg \in X^X$ *be a function on a* BOUNDED LATTICE $\boldsymbol{L} \triangleq ( X, \vee, \wedge, 0, 1\,; \leq)$.

$$\left\{ \; x \wedge \neg x \; = \; 0 \quad \forall x \in X \quad \text{(NON-CONTRADICTION)} \; \right\} \implies \left\{ \; \neg 1 \; = \; 0 \quad \text{(BOUNDARY CONDITION)} \; \right\}$$

PROOF:

$$
\begin{array}{ll}
0 = 1 \wedge \neg 1 & \text{by } \textit{non-contradiction} \text{ hypothesis} \\
\; = \neg 1 & \text{by definition of g.u.b. 1 and } \wedge
\end{array}
$$

**Lemma B.10** [78] *Let* $\neg \in X^X$ *be a function on a* BOUNDED LATTICE $\boldsymbol{L} \triangleq ( X, \vee, \wedge, 0, 1\,; \leq)$.

$$\left\{ \begin{array}{ll} (A). & \neg \textit{ is } \text{BIJECTIVE} \qquad\qquad and \\ (B). & x \leq y \implies \neg y \leq \neg x \quad \forall x, y \in X \quad \text{(ANTITONE)} \end{array} \right\} \implies \underbrace{\left\{ \begin{array}{ll} (1). & \neg 0 \; = \; 1 \quad and \\ (2). & \neg 1 \; = \; 0 \end{array} \right\}}_{\text{BOUNDARY CONDITIONS}}$$

**Theorem B.11** *Let* $\neg \in X^X$ *be a function on a* BOUNDED LATTICE $\boldsymbol{L} \triangleq ( X, \vee, \wedge, 0, 1\,; \leq)$.

$$\left\{ \begin{array}{l} \neg \textit{ is a} \\ \text{FUZZY NEGATION} \end{array} \right\} \implies \left\{ \; \neg 0 \; = \; 1 \quad \text{(BOUNDARY CONDITION)} \; \right\}$$

PROOF:   This follows directly from Definition B.2 (page 35) and Lemma B.8 (page 36).

**Theorem B.12** *Let* $\neg \in X^X$ *be a function on a* BOUNDED LATTICE $\boldsymbol{L} \triangleq ( X, \vee, \wedge, 0, 1\,; \leq)$.

$$\left\{ \begin{array}{l} \neg \textit{ is an} \\ \text{INTUITIONISTIC NEGATION} \end{array} \right\} \implies \left\{ \begin{array}{lll} (a) & \neg 1 \; = \; 0 \quad \text{(BOUNDARY CONDITION)} \quad and \\ (b) & \neg 0 \; = \; 1 \quad \text{(BOUNDARY CONDITION)} \quad and \\ (c) & \neg \textit{ is a } \text{FUZZY NEGATION} \end{array} \right\}$$

---

✎Proof:

$$\neg \text{ is an } \textit{intuitionistic negation} \implies x \wedge \neg x = 0 \qquad \text{by Definition B.2 page 35}$$

$$\implies \boxed{\neg 1 = 0} \qquad \text{by Lemma B.9 page 36}$$

$$\implies \neg \text{ is a } \textit{fuzzy negation} \qquad \text{by Definition B.2 page 35}$$

$$\implies \boxed{\neg 0 = 1} \qquad \text{by Theorem B.11 page 36}$$

**Theorem B.13**   *Let* $\neg \in X^X$ *be a function on a* BOUNDED LATTICE $\boldsymbol{L} \triangleq (\, X, \, \vee, \, \wedge, \, 0, \, 1 \,;\, \leq).$

$$\left\{ \begin{array}{l} \neg \text{ is a} \\ \textit{minimal} \\ \textit{negation} \end{array} \right\} \implies \left\{ \begin{array}{llll} \neg x \vee \neg y & \leq & \neg(x \wedge y) & \forall_{x,y \in X} & \text{(CONJUNCTIVE DE MORGAN INEQUALITY)} \quad and \\ \neg(x \vee y) & \leq & \neg x \wedge \neg y & \forall_{x,y \in X} & \text{(DISJUNCTIVE DE MORGAN INEQUALITY)} \end{array} \right\}$$

✎Proof:   This follows directly from Definition B.5 (page 35) and Lemma B.6 (page 35).

**Theorem B.14**   *Let* $\neg \in X^X$ *be a function on a* BOUNDED LATTICE $\boldsymbol{L} \triangleq (\, X, \, \vee, \, \wedge, \, 0, \, 1 \,;\, \leq).$

$$\left. \begin{array}{l} \neg \text{ is a} \\ \textit{de Morgan negation} \end{array} \right\} \implies \left\{ \begin{array}{llll} \neg(x \vee y) & = & \neg x \wedge \neg y & \forall_{x,y \in X} & \text{(DISJUNCTIVE DE MORGAN)} \quad and \\ \neg(x \wedge y) & = & \neg x \vee \neg y & \forall_{x,y \in X} & \text{(CONJUNCTIVE DE MORGAN)} \end{array} \right\}$$

✎Proof:   This follows directly from Definition B.5 (page 35) and Lemma B.7 (page 36).

**Theorem B.15**   *Let* $\neg \in X^X$ *be a function on a* BOUNDED LATTICE $\boldsymbol{L} \triangleq (\, X, \, \vee, \, \wedge, \, 0, \, 1 \,;\, \leq).$

$$\left\{ \begin{array}{l} \neg \text{ is an} \\ \textit{ortho} \\ \textit{negation} \end{array} \right\} \implies \left\{ \begin{array}{llllll} 1. & \neg 0 & = & 1 & & \text{(BOUNDARY CONDITION)} \quad and \\ 2. & \neg 1 & = & 0 & & \text{(BOUNDARY CONDITION)} \quad and \\ 3. & \neg(x \vee y) & = & \neg x \wedge \neg y & \forall_{x,y \in X} & \text{(DISJUNCTIVE DE MORGAN)} \quad and \\ 4. & \neg(x \wedge y) & = & \neg x \vee \neg y & \forall_{x,y \in X} & \text{(CONJUNCTIVE DE MORGAN)} \quad and \\ 5. & x \vee \neg x & = & 1 & \forall_{x \in X} & \text{(EXCLUDED MIDDLE)} \quad and \\ 6. & x \wedge \neg x & \leq & y \vee \neg y & \forall_{x,y \in X} & \text{(KLEENE CONDITION)}. \end{array} \right\}$$

✎Proof:

(1)  Proof for $0 = \neg 1$ boundary condition: by Lemma B.9 (page 36)

(2)  Proof for boundary conditions:

$$1 = \neg\neg 1 \qquad \text{by } \textit{involutory} \text{ property}$$

$$= \neg 0 \qquad \text{by previous result}$$

(3)  Proof for *de Morgan* properties:

    (a)  By Definition B.5 (page 35), *ortho negation* is *involutory* and *antitone*.

    (b)  Therefore by Lemma B.7 (page 36), *de Morgan* properties hold.

(4)  Proof for *excluded middle* property:

$$x \vee \neg x = \neg\neg(x \vee \neg x) \qquad \text{by } \textit{involutory} \text{ property of } \textit{ortho negation} \text{ (Definition B.5 35)}$$

$$= \neg(x\neg \wedge [\neg\neg x]) \qquad \text{by } \textit{disjunctive de Morgan} \text{ property}$$

$$= \neg(\neg x \wedge x) \qquad \text{by } \textit{involutory} \text{ property of } \textit{ortho negation} \text{ (Definition B.5 page 35)}$$

$$= \neg(x \wedge \neg x) \qquad \text{by } \textit{commutative} \text{ property of } \textit{lattices} \text{ (Definition A.11 page 24)}$$

$$= \neg 0 \qquad \text{by } \textit{non-contradiction} \text{ property of } \textit{ortho negation} \text{ (Definition B.5 35)}$$

$$= 1 \qquad \text{by } \textit{boundary condition} \text{ (item (2) page 37) of } \textit{minimal negation}$$

(5)  Proof for *Kleene condition*:

$$x \wedge \neg x = 0 \qquad\qquad \text{by } \textit{non-contradiction} \text{ property (Definition B.5 page 35)}$$
$$\leq 1 \qquad\qquad\qquad \text{by definition of 0 and 1}$$
$$= y \vee \neg y \qquad\qquad \text{by } \textit{excluded middle} \text{ property (item (4) page 37)}$$

☞

## B.3   Examples

**Example B.16**  (discrete negation)  [79]  Let $L \triangleq (\, X, \vee, \wedge, \neg, 0, 1\,; \leq\,)$ be a *bounded lattice* (Definition A.19 page 25) with a function $\neg \in X^X$.  The function $\neg x$ defined as

$$\neg x \triangleq \begin{cases} 1 & \text{for } x = 0 \\ 0 & \text{otherwise} \end{cases}$$

is an **intuitionistic negation** (Definition B.2 page 35, and a *fuzzy negation*).

✎Proof:      To be an *intuitionistic negation*, $\neg x$ must be *antitone*, have *weak double negation*, and have the *non-contradiction* property (Definition B.2 page 35).

$$\left. \begin{array}{lll} \neg y \leq \neg x & \iff & 1 \leq 1 \quad \text{for } 0 = x = y \\ \neg y \leq \neg x & \iff & 0 \leq 1 \quad \text{for } 0 = x \leq y \\ \neg y \leq \neg x & \iff & 0 \leq 0 \quad \text{for } 0 \neq x \leq y \end{array} \right\} \implies \neg x \text{ is } \textit{antitone}$$

$$\left. \begin{array}{lllllll} \neg\neg x & = & \neg 1 & = & 0 & \geq & 0 & = & x & \text{for } x = 0 \\ \neg\neg x & \neq & \neg 0 & = & 1 & \geq & x & = & x & \text{for } x \neq 0 \end{array} \right\} \implies \neg x \text{ has } \textit{weak double negation}$$

$$\left. \begin{array}{lllllll} x \wedge \neg x & = & x \wedge 1 & = & 0 \wedge 0 & = & 0 & \text{for } x = 0 \\ x \wedge \neg x & = & x \wedge 0 & = & x \wedge 0 & = & 0 & \text{for } x \neq 0 \end{array} \right\} \implies \neg x \text{ has } \textit{non-contradiction} \text{ property}$$

☞

**Example B.17**  (dual discrete negation)  [80]  Let $L \triangleq (\, X, \vee, \wedge, \neg, 0, 1\,; \leq\,)$ be a *bounded lattice* (Definition A.19 page 25) with a function $\neg \in X^X$.  The function $\neg x$ defined as

$$\neg x \triangleq \begin{cases} 0 & \text{for } x = 1 \\ 1 & \text{otherwise} \end{cases}$$

is a *subminimal negation* (Definition B.1 page 34) but it is *not* a *minimal negation* (Definition B.2 page 35) (and not any other negation defined here).

✎Proof:      To be a *subminimal negation*, $\neg x$ must be *antitone* (Definition B.1 page 34).  To be a *minimal negation*, $\neg x$ must be *antitone* and have *weak double negation* (Definition B.2 page 35).

$$\left. \begin{array}{lll} \neg y \leq \neg x & \iff & 0 \leq 0 \quad \text{for } x = y = 1 \\ \neg y \leq \neg x & \iff & 0 \leq 1 \quad \text{for } x \leq y = 1 \\ \neg y \leq \neg x & \iff & 1 \leq 1 \quad \text{for } x \leq y \neq 1 \end{array} \right\} \implies \neg x \text{ is } \textit{antitone}$$

$$\left. \begin{array}{lllll} \neg\neg x & = & \neg 0 & = & 1 & \geq & x & \text{for } x = 1 \\ \neg\neg x & = & \neg 1 & = & 0 & \leq & x & \text{for } x \neq 1 \end{array} \right\} \implies \neg x \text{ does } \textit{not} \text{ have } \textit{weak double negation}$$

☞

---

**Example B.18**

The function ¬ illustrated to the right is an *ortho negation* (Definition B.3 page 35).

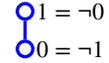

✎Proof:

(1) Proof that ¬ is *antitone*:
$0 \leq 1 \implies \neg 1 = 0 \leq x = \neg 0 \implies \neg$ is *antitone* over $(0, 1)$

(2) Proof that ¬ is *involutory*: $1 = \neg 0 = \neg\neg 1$

(3) Proof that ¬ has the *non-contradiction* property:   $1 \ \wedge \ \neg 1 \ = \ 1 \ \wedge \ 0 \ = \ 0$
$0 \ \wedge \ \neg 0 \ = \ 0 \ \wedge \ 1 \ = \ 0$

**Example B.19**

The functions ¬ illustrated to the right are *not* any negation defined here. In particular, none of them is *antitone*.

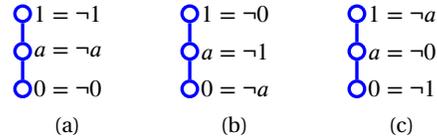

✎Proof:

1. Proof that (a) is *not antitone*: $a \leq 1 \implies \neg 1 = 1 \nleq a = \neg a$
2. Proof that (b) is *not antitone*: $a \leq 1 \implies \neg 1 = a \nleq 0 = \neg a$
3. Proof that (c) is *not antitone*: $0 \leq a \implies \neg a = 1 \nleq a = \neg 0$

**Example B.20**   The function ¬ as illustrated to the right is *not* a *sub-minimal negation* (it is *not antitone*) and so is *not* any negation defined here. Note however that the problem is *not* the $O_6$ *lattice*—it is possible to define a negation on an $O_6$ *lattice* (Example B.31 page 43).
✎Proof:   Proof that ¬ is *not antitone*: $a \leq c \implies \neg c = d \nleq b = \neg a$

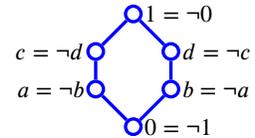

**Remark B.21**   The concept of a *complement* (Definition A.35 page 28) and the concept of a *negation* are fundamentally different. A *complement* is a *relation* on a lattice **L** and a *negation* is a *function*. In Example B.20 (page 39), $b$ and $d$ are both complements of $a$ (and so the lattice is *multiply complemented*), but yet ¬ is *not* a negation. In the right side lattice of Example B.31 (page 43), both $b$ and $d$ are complements of $a$, but yet only $d$ is equal to the negation of $a$ ($d = \neg a$). It can also be said that complementation is a property *of* a lattice, whereas negation is a function defined *on* a lattice.

**Remark B.22**   If a lattice is *complemented*, then by definition each element $x$ in the lattice has a *complement* $x'$ such that $x \wedge x' = 0$ (*non-contradiction* property) and $x \vee x' = 1$ (*excluded middle* property). If a lattice **L** is both *distributive* and *complemented*, then **L** is *uniquely complemented* (Definition A.41 page 30, Theorem A.42 page 30). If **L** is *uniquely complemented* and satisfies any one of *Huntington's properties* (**L** is *modular*, *atomic*, *ortho-complemented*, has *finite width*, or *de Morgan*), then **L** is *distributive* (Theorem A.40 page 29).

**Example B.23**   Each of the functions ¬ illustrated in Figure 10 (page 40) is a *subminimal negation* (Definition B.1 page 34); *none* of them is a *minimal negation* (each fails to have *weak double negation*).

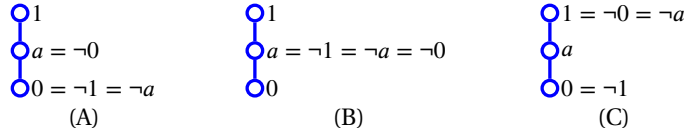

Figure 10: *subminimal negations* on $\mathbf{L}_3$ (Example B.23 page 39)

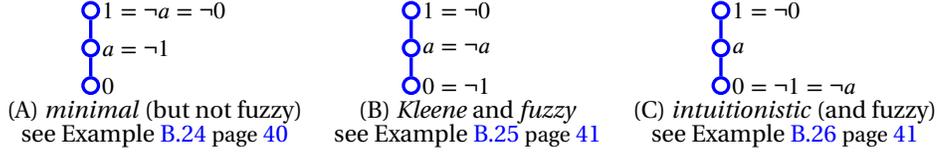

Figure 11: negations on $\mathbf{L}_3$

✎ PROOF:

(1) Proof that (A) ¬ is *antitone*:

$$a \le 1 \implies \neg 1 = 0 \le 0 = \neg a \implies \neg \text{ is } antitone \text{ over } (a,1)$$
$$0 \le 1 \implies \neg 1 = 0 \le a = \neg 0 \implies \neg \text{ is } antitone \text{ over } (0,1)$$
$$0 \le a \implies \neg a = 0 \le a = \neg 0 \implies \neg \text{ is } antitone \text{ over } (0,a)$$

(2) Proof that (A) ¬ *fails* to have *weak double negation*:
$$1 \not\le a = \neg 0 = \neg\neg 1$$

(3) Proof that (B) ¬ is *antitone*:

$$a \le 1 \implies \neg 1 = a \le a = \neg a \implies \neg \text{ is } antitone \text{ over } (a,1)$$
$$0 \le 1 \implies \neg 1 = a \le a = \neg 0 \implies \neg \text{ is } antitone \text{ over } (0,1)$$
$$0 \le a \implies \neg a = a \le a = \neg 0 \implies \neg \text{ is } antitone \text{ over } (0,a)$$

(4) Proof that (B) ¬ *fails* to have *weak double negation*: $1 \not\le a = \neg a = \neg\neg 1$

(5) (C) is a special case of the *dual discrete negation* (Example B.17 page 38).

✎

**Example B.24** Consider the function ¬ on $\mathbf{L}_3$ illustrated in Figure 11 page 40 (A):

1. ¬ is a **minimal negation** (Definition B.2 page 35);
2. ¬ is *not* an *intuitionistic negation* and it is *not* a *de Morgan negation*.

✎ PROOF:

(1) Proof that ¬ is *antitone*:

$$a \le 1 \implies \neg 1 = a \le 1 = \neg a \implies \neg \text{ is } antitone \text{ over } (a,1)$$
$$0 \le 1 \implies \neg 1 = a \le 1 = \neg 0 \implies \neg \text{ is } antitone \text{ over } (0,1)$$
$$0 \le a \implies \neg a = 1 \le 1 = \neg 0 \implies \neg \text{ is } antitone \text{ over } (0,a)$$

(2) Proof that ¬ is a *weak double negation* (and so is a *minimal negation*, but is *not* a *de Morgan negation*):

$$1 = 1 = \neg a = \neg\neg 1 \implies \neg \text{ is } involutory \text{ at } 1$$
$$a = a = \neg 1 = \neg\neg a \implies \neg \text{ is } involutory \text{ at } a$$
$$0 \le a = \neg 1 = \neg\neg 0 \implies \neg \text{ is a } weak \text{ double negation at } 0$$

(3) Proof that ¬ does *not* have the *non-contradiction* property (and so is not an *intuitionistic negation*):
$$1 \wedge \neg 1 = 1 \wedge a = a \ne 0$$

(4) Proof that ¬ is not a *fuzzy negation*: $\neg 1 = a \ne 0$

**Example B.25**   [81] Consider the function ¬ on $L_3$ illustrated in Figure 11 page 40 (B).

1. ¬ is a **Kleene negation** (Definition B.3 page 35) and is also a **fuzzy negation** (Definition B.2 page 35, Example 1.31 page 16).

2. ¬ is *not* an *ortho negation* and is *not* an *intuitionistic negation* (it does not have the *non-contradiction* property).

3. This negation on $L_3$ is used with an *implication* function to construct the *Kleene 3-valued logic* in Example C.7 (page 51), with another *implication* to construct the *Łukasiewicz 3-valued logic* in Example C.8 (page 52), and with yet another *implication* to construct the $RM_3$ *logic* in Example C.9 (page 52).

✎ PROOF:

(1) Proof that ¬ is *antitone*:

$$a \leq 1 \implies \neg 1 = 0 \leq a = \neg a \implies \neg \text{ is } antitone \text{ over } (a, 1)$$
$$0 \leq 1 \implies \neg 1 = 0 \leq 1 = \neg 0 \implies \neg \text{ is } antitone \text{ over } (0, 1)$$
$$0 \leq a \implies \neg a = a \leq 1 = \neg 0 \implies \neg \text{ is } antitone \text{ over } (0, a)$$

(2) Proof that ¬ is *involutory* (and so is a *de Morgan negation*):

$$1 = \neg 0 = \neg\neg 1 \implies \neg \text{ is } involutory \text{ at } 1$$
$$a = \neg a = \neg\neg a \implies \neg \text{ is } involutory \text{ at } a$$
$$0 = \neg 0 = \neg\neg 0 \implies \neg \text{ is } involutory \text{ at } 0$$

(3) Proof that ¬ does *not* have the *non-contradiction* property (and so is not an *ortho negation*):

$$x \wedge \neg x = x \wedge x = x \neq 0$$

(4) Proof that ¬ satisfies the *Kleene condition* (and so is a *Kleene negation*):

$$1 \wedge \neg 1 = 1 \wedge 0 = 0 \leq a = a \vee a = a \vee \neg a$$
$$1 \wedge \neg 1 = 1 \wedge 0 = 0 \leq 1 = 0 \vee 1 = 0 \vee \neg 0$$
$$a \wedge \neg a = 1 \wedge a = a \leq 1 = 1 \vee 0 = 1 \vee \neg 1$$
$$a \wedge \neg a = 1 \wedge a = a \leq 1 = 0 \vee 1 = 0 \vee \neg 0$$
$$0 \wedge \neg 0 = 0 \wedge 1 = 0 \leq 1 = 1 \vee 0 = 1 \vee \neg 1$$
$$0 \wedge \neg 0 = 0 \wedge 1 = 0 \leq a = a \vee a = a \vee \neg a$$

**Example B.26**   (Heyting 3-valued logic/Jaśkowski's first matrix)   [82] Consider the the function ¬ on $L_3$ illustrated in Figure 11 page 40 (C):

1. ¬ is an **intuitionistic negation** (Definition B.2 page 35) (and thus is also a **fuzzy negation**).
2. ¬ is *not* a *de Morgan negation*.
3. This negation on $L_3$ is used with an *implication* function to construct the *Heyting 3-valued logic* in Example C.10 (page 52).

✎ PROOF:   This is simply a special case of the *discrete negation* (Example B.16 page 38).

---

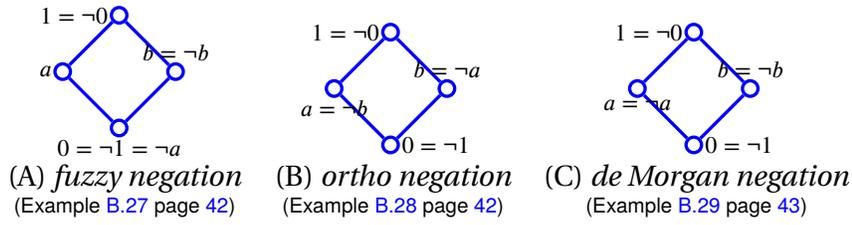

(A) *fuzzy negation*
(Example B.27 page 42)

(B) *ortho negation*
(Example B.28 page 42)

(C) *de Morgan negation*
(Example B.29 page 43)

Figure 12: negations on $M_2$

**Example B.27**   The function $\neg$ illustrated in Figure 12 page 42 (A) is a **fuzzy negation** (Definition B.2 page 35). It is not an *intuitionistic negation* (it does not have the *non-contradiction* property) and it is *not* a *de Morgan negation* (it is not *involutory*).

✎ PROOF:    Note that

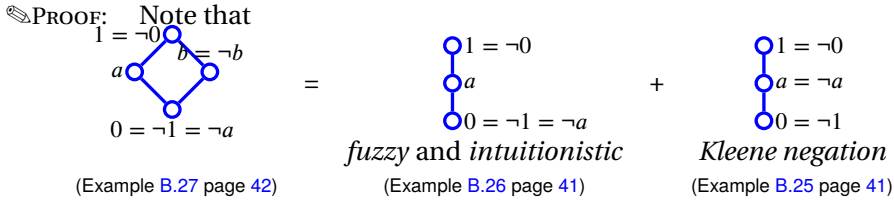

*fuzzy* and *intuitionistic*
(Example B.26 page 41)

*Kleene negation*
(Example B.25 page 41)

(1)  Proof that $\neg$ is *antitone*:

$$a \;\le\; 1 \;\implies\; \neg 1 \;=\; 0 \;\le\; 0 \;=\; \neg a \;\implies\; \neg \text{ is } \textit{antitone at } (a, 1)$$
$$0 \;\le\; 1 \;\implies\; \neg 1 \;=\; 0 \;\le\; 1 \;=\; \neg 0 \;\implies\; \neg \text{ is } \textit{antitone at } (0, 1)$$
$$0 \;\le\; a \;\implies\; \neg a \;=\; 0 \;\le\; 1 \;=\; \neg 0 \;\implies\; \neg \text{ is } \textit{antitone at } (0, a)$$
$$b \;\le\; 1 \;\implies\; \neg 1 \;=\; 0 \;\le\; b \;=\; \neg b \;\implies\; \neg \text{ is } \textit{antitone at } (b, 1)$$
$$0 \;\le\; b \;\implies\; \neg b \;=\; b \;\le\; 1 \;=\; \neg 0 \;\implies\; \neg \text{ is } \textit{antitone at } (0, b)$$

(2)  Proof that $\neg$ has *weak double negation* property (and so is a *minimal negation*, but *not* a *de Morgan negation*):

$$1 \;=\; \neg 0 \;=\; \neg\neg 1 \qquad\qquad \implies \quad \neg \text{ is } \textit{involutory at } 1$$
$$a \;\le\; 1 \;=\; \neg 0 \;=\; \neg\neg a \quad \implies \quad \neg \text{ has } \textit{weak double negation at } a$$
$$0 \;=\; \neg 1 \;=\; \neg\neg 0 \qquad\qquad \implies \quad \neg \text{ is } \textit{involutory at } 0$$
$$b \;=\; \neg b \;=\; \neg\neg b \;=\; \qquad\quad \implies \quad \neg \text{ is } \textit{involutory at } b$$

(3)  Proof that $\neg$ does *not* have the *non-contradiction* property (and so is *not* an *intuitionistic negation*):
$b \wedge \neg b = b \wedge b = b \ne 0$

(4)  Proof that $\neg$ is has *boundary conditions* (and so is a *fuzzy negation*): $\neg 1 = 0, \quad \neg 0 = 1$

**Example B.28**   [83] The function $\neg$ illustrated in Figure 12 page 42 (B) is an **ortho negation** (Definition B.3 page 35).

✎ PROOF:

(1)  Proof that $\neg$ is *antitone*:  $a \;\le\; 1 \;\implies\; \neg 1 \;=\; 0 \;\le\; b \;=\; \neg a$
$$0 \;\le\; 1 \;\implies\; \neg 1 \;=\; 0 \;\le\; 1 \;=\; \neg 0$$
$$0 \;\le\; a \;\implies\; \neg a \;=\; b \;\le\; 1 \;=\; \neg 0$$
$$b \;\le\; 1 \;\implies\; \neg 1 \;=\; 0 \;\le\; a \;=\; \neg b$$
$$0 \;\le\; b \;\implies\; \neg b \;=\; a \;\le\; 1 \;=\; \neg 0$$

---

[83] 📖 [Belnap(1977)] page 13, ✎ [Restall(2000)] page 177 ⟨Example 8.44⟩, 📖 [Pavičić and Megill(2009)] page 28 ⟨Definition 2, *classical implication*⟩

(2) Proof that ¬ is *involutory* (and so is a *de Morgan negation*):  $\begin{aligned} 1 &= \neg 0 = \neg\neg 1 \\ a &= \neg a = \neg\neg a \\ b &= \neg b = \neg\neg b \\ 0 &= \neg 0 = \neg\neg 0 \end{aligned}$

(3) Proof that ¬ is has the *non-contradiction* property (and so is an *ortho negation*):

$$\begin{aligned} 1 \;\wedge\; \neg 1 &= 1 \;\wedge\; 0 = 0 \\ a \;\wedge\; \neg a &= a \;\wedge\; b = 0 \\ b \;\wedge\; \neg b &= b \;\wedge\; a = 0 \\ 0 \;\wedge\; \neg 0 &= 0 \;\wedge\; 1 = 0 \end{aligned}$$

**Example B.29**  (BN$_4$)  [84] The function ¬ illustrated in Figure 12 page 42 (C) is a **de Morgan negation** (Definition B.3 page 35), but it is *not* a *Kleene negation* and not an *ortho negation* (it does *not* satisfy the *Kleene condition*).

✎ PROOF:

(1) Proof that ¬ is *antitone*:  $\begin{aligned} a &\leq 1 \implies \neg 1 = 0 \leq b = \neg a \\ 0 &\leq 1 \implies \neg 1 = 0 \leq 1 = \neg 0 \\ 0 &\leq a \implies \neg a = a \leq 1 = \neg 0 \\ b &\leq 1 \implies \neg 1 = 0 \leq b = \neg b \\ 0 &\leq b \implies \neg b = b \leq 1 = \neg 0 \end{aligned}$

(2) Proof that ¬ is *involutory* (and so is a *de Morgan negation*):  $\begin{aligned} 1 &= \neg 0 = \neg\neg 1 \\ a &= \neg a = \neg\neg a \\ b &= \neg b = \neg\neg b \\ 0 &= \neg 0 = \neg\neg 0 \end{aligned}$

(3) Proof that ¬ does *not* have the *non-contradiction* property (and so is *not* an *ortho negation*):

$$\begin{aligned} a \;\wedge\; \neg a &= a \;\wedge\; a = a \neq 0 \\ b \;\wedge\; \neg b &= b \;\wedge\; b = b \neq 0 \end{aligned}$$

(4) Proof that ¬ does *not* satisfy the *Kleene condition* (and so is a *de Morgan negation*):

$$a \;\wedge\; \neg a = a \;\wedge\; a = a \nleq b \;\wedge\; \neg b = b$$

**Example B.30**

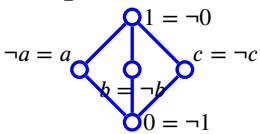

The function ¬ illustrated to the left is a **de Morgan negation**, but it is *not* a *Kleene negation* and not an *ortho negation*. The *negation* illustrated to the right is a **Kleene negation**, but it is *not* an *ortho negation*.

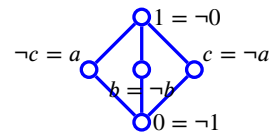

**Example B.31**

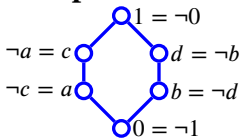

The function ¬ illustrated to the left is a **de Morgan negation** (Definition B.3 page 35); it is *not* a *Kleene negation* (it does not satisfy the Kleene condition). The *negation* illustrated to the right is an **ortho negation** (Definition B.3 page 35).

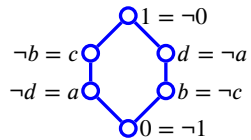

**Example B.32**

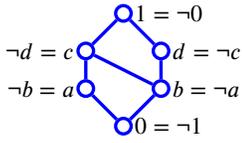

The function ¬ illustrated to the left is **not antitone** and therefore is not a *negation* (Definition B.2 page 35). The function ¬ illustrated to the right is a **Kleene negation** (Definition B.3 page 35); it is *not an ortho negation* (it does not have the *non-contradiction* property).

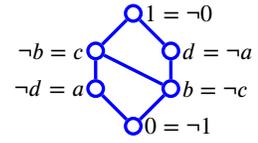

✎ Proof:

(1) Proof that left ¬ is *not antitone*: $a \leq c$ but $\neg c \nleq \neg a$.

(2) Proof that right ¬ satisfies the *Kleene condition*:
$$x \wedge \neg x = \begin{cases} b & \text{for } x = b \\ 0 & \text{otherwise} \end{cases} \quad \forall_{x \in X} \quad \text{and} \quad y \wedge \neg y = \begin{cases} c & \text{for } y = c \\ 0 & \text{otherwise} \end{cases} \quad \forall_{y \in X}$$
$$\implies \quad x \wedge \neg x \leq y \vee \neg y \quad \forall_{x,y \in X}$$

(3) Proof that right ¬ does not have the *non-contradiction* property: $b \wedge \neg b = b \wedge c = b \neq 0$

**Example B.33**

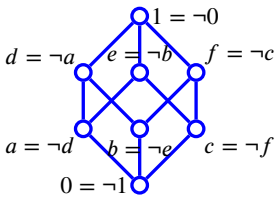

The lattices illustrated to the left and right are *Boolean* (Definition A.41 page 30). The function ¬ illustrated to the left is a **Kleene negation** (Definition B.3 page 35), but it is *not an ortho negation* (it does *not* have the *non-contradiction* property). The *negation* illustrated to the right is an **ortho negation** (Definition B.3 page 35).

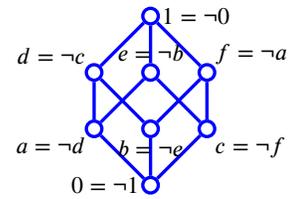

✎ Proof:

(1) Proof that left side negation does *not* have *non-contradiction* property (and so is *not* an *ortho negation*):
$a \wedge \neg a = a \wedge d = a \neq 0$

(2) Proof that left side negation does *not* satisfy *Kleene condition* (and so is *not* a *Kleene negation*):
$a \wedge \neg a = a \wedge d = a \nleq f = c \vee f = c \vee \neg c$

# Appendix C    New implication functions for non-Boolean logics

## C.1    Implication functions

This paper deals with how to construct a *fuzzy subset logic* not only on a *Boolean lattice*, but more generally on other types of *lattice*s as well. However, any logic (fuzzy or otherwise) is arguably not complete without the inclusion of an *implication* function →. If we were only concerned with logics on *Boolean lattice*s, then arguably the *classical implication* $x \xrightarrow{c} y \triangleq \neg x \vee y$ would suffice. However, for some *non-Boolean* lattices, we may do well to have other options. Two common properties of *classical implication* are *entailment* and *modus ponens*. However, these properties do not always support well known logic systems that are constructed on *non-orthocomplemented* (and hence also *non-Boolean*) lattices. For example,

-  the $RM_3$ *logic* does not support the *strong entailment* property,
-  the *Łukasiewicz 3-valued logic* does not support the *strong modus ponens* property, and
-  the *Kleene 3-valued logic* and $BN_4$ *logic* do not support either property.

This section introduces a new definition for an *implication* function with weakened forms of *entailment* and *modus ponens* (herein called *weak entailment* and *weak modus ponens*), and that supports logics constructed on a large class of lattices including *non-orthocomplemented* (and *non-Boolean*) ones.

**Definition C.1** Let $\boldsymbol{L} \triangleq (X, \vee, \wedge, 0, 1; \leq)$ be a *bounded lattice* (Definition A.19 page 25). The function $\to$ in $X^X$ is an **implication** on $\boldsymbol{L}$ if

1. $\{x \leq y\} \implies x \to y \geq x \vee y$    $\forall x,y \in X$    *(weak entailment)*    and
2. $x \wedge (x \to y) \leq \neg x \vee y$    $\forall x,y \in X$    *(weak modus ponens)*

**Proposition C.2** *Let $\to$ be an* IMPLICATION *(Definition C.1 page 45) on a* BOUNDED LATTICE $\boldsymbol{L} \triangleq (X, \vee, \wedge, 0, 1; \leq)$ *(Definition A.19 page 25).*

$\{x \leq y\}$    $\iff$    $\{x \to y \geq x \vee y\}$    $\forall x,y \in X$

✎PROOF:

(1) Proof for $\implies$ case: by *weak entailment* property of *implication*s (Definition C.1 page 45).

(2) Proof for $\impliedby$ case:

| | |
|---|---|
| $y \geq x \wedge (x \to y)$ | by right hypothesis |
| $\geq x \wedge (x \vee y)$ | by *modus ponens* property of $\to$ (Definition C.1 page 45) |
| $= x$ | by *absorptive* property of *lattice*s (Definition A.11 page 24) |

**Remark C.3** Let $\boldsymbol{L} \triangleq (X, \vee, \wedge, 0, 1; \leq)$ be a *bounded lattice* (Definition A.19 page 25). In the context of *orthocomplemented lattice*s, a more common (and stronger) definition of *implication* $\to$ might be[85]

1. $x \leq y \implies x \to y = 1$    $\forall x,y \in X$    *(entailment / strong entailment)*    and
2. $x \wedge (x \to y) \leq y$    $\forall x,y \in X$    *(modus ponens / strong modus ponens)*

This definition yields a result stronger than that of Proposition C.2 (page 45):

$\{x \leq y\}$    $\iff$    $\{x \to y = 1\}$    $\forall x,y \in X$

The *Heyting 3-valued logic* (Example C.10 page 52) and *Boolean 4-valued logic* (Example C.12 page 53) have both *strong entailment* and *strong modus ponens*. However, for non-Boolean logics in general, these two properties seem inappropriate to serve as a definition for *implication*. For example, the *Kleene 3-valued logic* (Example C.7 page 51), $RM_3$ *logic* (Example C.9 page 52), and $BN_4$ *logic* (Example C.13 page 54) do not have the *strong entailment* property; and the *Kleene 3-valued logic*, *Łukasiewicz 3-valued logic* (Example C.8 page 52), and $BN_4$ *logic* do not have the *strong modus ponens* property.

✎PROOF:

(1) Proof for $\implies$ case: by *entailment* property of *implication*s (Definition C.1 page 45).

(2)   Proof for $\Longleftarrow$ case:

$$x \to y = 1 \implies x \land 1 \le y \qquad \text{by } \textit{modus ponens} \text{ property (Definition C.1 page 45)}$$
$$\implies x \le y \qquad\qquad \text{by definition of 1 } (\textit{least upper bound}) \text{ (Definition A.8 page 24)}$$

$\hookleftarrow$

**Example C.4** [86] Let $\textbf{L} \triangleq (X, \lor, \land, \neg, 0, 1 ; \le)$ be a *lattice with negation* (Definition B.5 page 35). If $\textbf{L}$ is an **orthocomplemented lattice**, then under Definition C.1, functions (1)–(5) below are valid *implication* functions with *strong entailment* and *weak modus ponens*. The *relevance implication* (6) in this lattice is *not* a valid implication: It does have *weak modus ponens*, but it does not have weak or strong entailment. However, if $\textbf{L}$ is an **orthomodular lattice** (Definition A.23 page 26, a special case of an orthocomplemented lattice), then (6) is also a valid implication function with *strong entailment*.

1.   $x \overset{c}{\to} y \triangleq \neg x \lor y \qquad \forall x, y \in X \qquad$ (*classical implication / material implication / horseshoe*)
2.   $x \overset{s}{\to} y \triangleq \neg x \lor (x \land y) \qquad\qquad\qquad \forall x, y \in X \qquad$ (*Sasaki hook / quantum implication*)
3.   $x \overset{d}{\to} y \triangleq y \lor (\neg x \land \neg y) \qquad\qquad \forall x, y \in X \qquad$ (*Dishkant implication*)
4.   $x \overset{k}{\to} y \triangleq (\neg x \land y) \lor (\neg x \land \neg y) \lor (x \land (\neg x \lor y)) \quad \forall x, y \in X \qquad$ (*Kalmbach implication*)
5.   $x \overset{n}{\to} y \triangleq (\neg x \land y) \lor (x \land y) \lor ((\neg x \lor y) \land \neg y) \quad \forall x, y \in X \qquad$ (*non-tollens implication*)
6.   $x \overset{r}{\to} y \triangleq (\neg x \land y) \lor (x \land y) \lor (\neg x \land \neg y) \qquad \forall x, y \in X \qquad$ (*relevance implication*)

Moreover, if $\textbf{L}$ is a **Boolean lattice**, then all of these implications are equivalent to $\overset{c}{\to}$, and all of them have *strong entailment* and *strong modus ponens*.

Note that $\forall x, y \in X$, $\quad x \overset{d}{\to} y = \neg y \overset{s}{\to} \neg x \quad$ and $\quad x \overset{n}{\to} y = \neg y \overset{k}{\to} \neg x$. The values for the six implications on an orthocomplemented $O_6$ *lattice* (Definition A.45 page 31) are listed in Example C.14 (page 54).

✎PROOF:

(1)   Proofs for the *classical implication* $\overset{c}{\to}$:

    (a)   Proof that on an *orthocomplemented lattice*, $\overset{c}{\to}$ is an *implication*:

$$x \le y \implies x \overset{c}{\to} y \triangleq \neg x \lor y \qquad \text{by definition of } \overset{c}{\to}$$
$$\ge \neg y \lor y \qquad \text{by } x \le y \text{ and } \textit{antitone} \text{ property of } \neg \text{ (Definition B.3 page 35)}$$
$$= 1 \qquad \text{by } \textit{excluded middle} \text{ property of } \neg \text{ (Theorem B.15 page 37)}$$
$$\implies \textit{strong entailment} \qquad \text{by definition of } \textit{strong entailment}$$
$$x \land (\neg x \lor y) \le \neg x \lor y \qquad \text{by definition of } \land \text{ (Definition A.9 page 24)}$$
$$\implies \textit{weak modus ponens} \qquad \text{by definition of } \textit{weak modus ponens}$$

Note that in general for an *orthocomplemented lattice*, the bound cannot be tightened to *strong modus ponens* because, for example in the $O_6$ *lattice* (Definition A.45 page 31) illustrated to the right

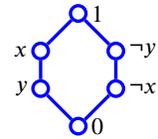

$x \land (\neg x \lor y) = x \land 1 = x \nleq y \implies \textit{not strong modus ponens}$

    (b)   Proof that on a *Boolean lattice*, $\overset{c}{\to}$ is an *implication*:

$$x \land (\neg x \lor y) = (x \land \neg x) \lor (x \land y) \qquad \text{by } \textit{distributive} \text{ property (Definition A.41 page 30)}$$
$$= 1 \lor (x \land y) \qquad \text{by } \textit{excluded middle} \text{ property of } \textit{Boolean lattice} \text{s}$$
$$= x \land y \qquad \text{by definition of 1}$$
$$\le y \qquad \text{by definition of } \land \text{ (Definition A.9 page 24)}$$
$$\implies \textit{strong modus ponens} \qquad \text{by definition of } \textit{strong modus ponens}$$

(2)   Proofs for *Sasaki implication* $\overset{s}{\to}$:

   (a)   Proof that on an *orthocomplemented lattice*, $\overset{s}{\to}$ is an *implication*:

$$x \le y \implies x \overset{s}{\to} y$$
$$\triangleq \neg x \vee (x \wedge y) \qquad\qquad \text{by definition of } \overset{k}{\to}$$
$$= \neg x \vee x \qquad\qquad\qquad \text{by } x \le y \text{ hypothesis}$$
$$= 1 \qquad\qquad\qquad\qquad \text{by } excluded\ middle \text{ prop. } \text{(Theorem B.15 page 37)}$$
$$\implies strong\ entailment \qquad \text{by definition of } strong\ entailment$$
$$x \wedge (x \overset{s}{\to} y) \triangleq x \wedge [\neg x \vee (x \wedge y)] \quad \text{by definition of } \overset{s}{\to}$$
$$\le [\neg x \vee (x \wedge y)] \qquad\qquad \text{by definition of } \wedge \text{ (Definition A.9 page 24)}$$
$$\le \neg x \vee y \qquad\qquad\qquad \text{by definition of } \wedge \text{ (Definition A.9 page 24)}$$
$$\implies weak\ modus\ ponens$$

   (b)   Proof that on a *Boolean lattice*, $\overset{s}{\to} = \overset{c}{\to}$:

$$x \overset{s}{\to} y \triangleq \neg x \vee (x \wedge y) \qquad\qquad \text{by definition of } \overset{s}{\to}$$
$$= \neg x \vee y \qquad\qquad\qquad\quad \text{by Lemma A.43 (page 30)}$$
$$= x \overset{c}{\to} y \qquad\qquad\qquad\quad \text{by definition of } \overset{c}{\to}$$

(3)   Proofs for *Dishkant implication* $\overset{d}{\to}$:

   (a)   Proof that $x \overset{d}{\to} y \equiv \neg y \overset{s}{\to} \neg x$:

$$x \overset{d}{\to} y \triangleq y \vee (\neg x \wedge \neg y) \qquad \text{by definition of } \overset{d}{\to}$$
$$= y \vee (\neg y \wedge \neg x) \qquad\qquad \text{by } commutative \text{ property of } lattice\text{s (Theorem A.14 page 25)}$$
$$= \neg\neg y \vee (\neg y \wedge \neg x) \qquad \text{by } involutory \text{ property of } ortho\ negation\text{s (Definition B.3 page 35)}$$
$$\triangleq \neg y \overset{s}{\to} \neg x \qquad\qquad\qquad \text{by definition of } \overset{s}{\to}$$

   (b)   Proof that on an *orthocomplemented lattice*, $\overset{d}{\to}$ is an *implication*:

$$x \le y \implies x \overset{d}{\to} y$$
$$\triangleq y \vee (\neg x \wedge \neg y) \qquad\qquad \text{by definition of } \overset{d}{\to}$$
$$= y \vee \neg y \qquad\qquad\qquad\quad \text{by } x \le y \text{ hypothesis and } antitone \text{ property}$$
$$= 1 \qquad\qquad\qquad\qquad\quad \text{by } excluded\ middle \text{ property of } ortho\ negation$$
$$\implies strong\ entailment \qquad \text{by definition of } strong\ entailment$$
$$x \wedge (x \overset{d}{\to} y) \triangleq y \vee (\neg x \wedge \neg y) \quad \text{by definition of } \overset{d}{\to}$$
$$= y \vee \neg x \qquad\qquad\qquad \text{by definition of } \wedge \text{ (Definition A.9 page 24)}$$
$$\implies weak\ modus\ ponens$$

   (c)   Proof that on a *Boolean lattice*, $\overset{d}{\to} = \overset{c}{\to}$:

$$x \overset{d}{\to} y \triangleq y \vee (\neg x \wedge \neg y) \qquad\qquad \text{by definition of } \overset{d}{\to}$$
$$= \neg x \vee y \qquad\qquad\qquad\qquad \text{by Lemma A.43 (page 30)}$$
$$= x \overset{c}{\to} y \qquad\qquad\qquad\qquad \text{by definition of } \overset{c}{\to}$$

(4)   Proofs for the *Kalmbach implication* $\overset{k}{\to}$:

(a) Proof that on an *orthocomplemented lattice*, $\overset{k}{\rightarrow}$ is an *implication*:

$$x \le y \implies x \overset{k}{\rightarrow} y$$

$$
\begin{aligned}
&\triangleq (\neg x \wedge y) \vee (\neg x \wedge \neg y) \vee [x \wedge (\neg x \vee y)] && \text{by definition of } \overset{k}{\rightarrow} \\
&= (\neg x \wedge y) \vee (\neg y) \vee [x \wedge (\neg x \vee y)] && \text{by } \textit{antitone} \text{ property } \text{(Definition B.3 page 35)} \\
&= (\neg x \wedge y) \vee \neg y \vee [x \wedge (1)] && \text{by definition of 1 } \text{(Definition A.8 page 24)} \\
&= (\neg x \wedge y) \vee (x \vee \neg y) && \\
&= \neg\neg(\neg x \wedge y) \vee (x \vee \neg y) && \text{by } \textit{involutory} \text{ property } \text{(Definition B.3 page 35)} \\
&= \neg(\neg\neg x \vee \neg y) \vee (x \vee \neg y) && \text{by } \textit{de Morgan} \text{ property } \text{(Theorem B.15 page 37)} \\
&= \neg(x \vee \neg y) \vee (x \vee \neg y) && \text{by } \textit{involutory} \text{ property } \text{(Definition B.3 page 35)} \\
&= 1 && \text{by } \textit{excluded middle} \text{ prop. } \text{(Theorem B.15 page 37)} \\
&\implies \textit{strong entailment} &&
\end{aligned}
$$

$$x \wedge (x \overset{k}{\rightarrow} y)$$

$$
\begin{aligned}
&\triangleq x \wedge [(\neg x \wedge y) \vee (\neg x \wedge \neg y) \vee [x \wedge (\neg x \vee y)]] && \text{by definition of } \overset{k}{\rightarrow} \\
&\le (\neg x \wedge y) \vee (\neg x \wedge \neg y) \vee [x \wedge (\neg x \vee y)] && \text{by definition of } \wedge \text{ (Definition A.9 page 24)} \\
&\le (\neg x \wedge y) \vee (\neg x \wedge \neg y) \vee (\neg x \vee y) && \text{by definition of } \wedge \text{ (Definition A.9 page 24)} \\
&\le y \vee (\neg x \wedge \neg y) \vee \neg x \vee y && \text{by definition of } \wedge \text{ (Definition A.9 page 24)} \\
&= y \vee \neg x \vee (\neg x \wedge \neg y) && \text{by } \textit{idempotent} \text{ p. } \text{(Theorem A.14 page 25)} \\
&\le y \vee \neg x \vee \neg x && \text{by definition of } \wedge \text{ (Definition A.9 page 24)} \\
&= \neg x \vee y && \text{by } \textit{idempotent} \text{ p. } \text{(Theorem A.14 page 25)} \\
&\implies \textit{weak modus ponens} &&
\end{aligned}
$$

(b) Proof that on a *Boolean lattice*, $\overset{k}{\rightarrow} = \overset{c}{\rightarrow}$:

$$
\begin{aligned}
&x \overset{k}{\rightarrow} y && \\
&\triangleq (\neg x \wedge y) \vee (\neg x \wedge \neg y) \vee [x \wedge (\neg x \vee y)] && \text{by definition of } \overset{k}{\rightarrow} \\
&= (\neg x \wedge y) \vee (\neg x \wedge \neg y) \vee [(x \wedge \neg x) \vee (x \wedge y)] && \text{by } \textit{distributive} \text{ property } \text{(Definition A.41 page 30)} \\
&= (\neg x \wedge y) \vee (\neg x \wedge \neg y) \vee [(0) \vee (x \wedge y)] && \text{by } \textit{non-contradiction} \text{ property} \\
&= (\neg x \wedge y) \vee (\neg x \wedge \neg y) \vee (x \wedge y) && \text{by } \textit{bounded} \text{ property } \text{(Definition A.19 page 25)} \\
&= \neg x \wedge (y \vee \neg y) \vee (x \wedge y) && \text{by } \textit{distributive} \text{ property } \text{(Definition A.41 page 30)} \\
&= \neg x \wedge 1 \vee (x \wedge y) && \text{by } \textit{excluded middle} \text{ property} \\
&= \neg x \vee (x \wedge y) && \text{by definition of 1 } \text{(Definition A.8 page 24)} \\
&= \neg x \vee y && \text{by Lemma A.43 (page 30)} \\
&\triangleq x \overset{c}{\rightarrow} y && \text{by definition of } \overset{c}{\rightarrow}
\end{aligned}
$$

(5) Proofs for the *non-tollens implication* $\overset{n}{\rightarrow}$:

(a) Proof that $x \overset{n}{\rightarrow} y \equiv \neg y \overset{k}{\rightarrow} \neg x$:

$$
\begin{aligned}
&x \overset{n}{\rightarrow} y \triangleq (\neg x \wedge y) \vee (x \wedge y) \vee [(\neg x \vee y) \wedge \neg y] && \text{by definition of } \overset{n}{\rightarrow} \\
&= (y \wedge \neg x) \vee (y \wedge x) \vee [\neg y \wedge (y \vee \neg x)] && \\
&= (\neg\neg y \wedge \neg x) \vee (\neg\neg y \wedge \neg\neg x) \vee [\neg y \wedge (\neg\neg y \vee \neg x)] && \\
&\triangleq \neg y \overset{k}{\rightarrow} \neg x && \text{by definition of } \overset{k}{\rightarrow}
\end{aligned}
$$

(b) Proof that on an *orthocomplemented lattice*, $\overset{n}{\to}$ is an *implication*:

$$x \le y \implies x \overset{n}{\to} y$$

$$\equiv \neg y \overset{k}{\to} \neg x \qquad\qquad \text{by item (5a) page 48}$$

$$= 1 \qquad\qquad\qquad \text{by item (4a) page 48}$$

$$\implies \textit{strong entailment}$$

$$x \wedge (x \overset{n}{\to} y) = x \wedge (\neg y \overset{k}{\to} \neg x) \qquad \text{by item (5a) page 48}$$

$$\le \neg\neg y \vee \neg x \qquad\qquad \text{by item (4a) page 48}$$

$$= y \vee \neg x \qquad\qquad\qquad \text{by \textit{involutory} property of } \neg \text{ (Definition B.3 page 35)}$$

$$= \neg x \vee y \qquad\qquad\qquad \text{by \textit{commutative} property of \textit{lattices}}$$

$$\implies \textit{weak modus ponens}$$

(c) Proof that on a *Boolean lattice*, $\overset{n}{\to} = \overset{c}{\to}$:

$$x \overset{n}{\to} y = \neg y \overset{k}{\to} \neg x \qquad \text{by item (5a) page 48}$$

$$= \neg\neg y \vee \neg x \qquad\qquad \text{by item (4b) page 48}$$

$$= y \vee \neg x \qquad\qquad\qquad \text{by \textit{involutory} property of } \neg \text{ (Definition B.3 page 35)}$$

$$= \neg x \vee y \qquad\qquad\qquad \text{by \textit{commutative} property of \textit{lattices} (Definition A.11 page 24)}$$

$$\triangleq x \overset{c}{\to} y \qquad\qquad\qquad \text{by definition of } \overset{c}{\to}$$

(6) Proofs for the *relevance implication* $\overset{r}{\to}$:

(a) Proof that on an *orthocomplemented lattice*, $\overset{r}{\to}$ does *not* have *weak entailment*:
In the *orthocomplemented lattice* to the right…

$$x \le y \implies x \overset{r}{\to} y$$

$$\triangleq (\neg x \wedge y) \vee (x \wedge y) \vee (\neg x \wedge \neg y) \qquad \text{by definition of } \overset{r}{\to}$$

$$= 0 \vee x \vee \neg y$$

$$= x \vee \neg y$$

$$\ne x \vee y$$

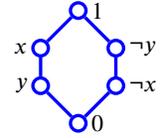

(b) Proof that on an *orthomodular lattice*, $\overset{r}{\to}$ *does* have *strong entailment*:

$$x \le y \implies x \overset{r}{\to} y$$

$$\triangleq (\neg x \wedge y) \vee (x \wedge y) \vee (\neg x \wedge \neg y) \qquad \text{by definition of } \overset{r}{\to}$$

$$= (\neg x \wedge y) \vee x \vee (\neg x \wedge \neg y) \qquad \text{by } x \le y \text{ hypothesis}$$

$$= (\neg x \wedge y) \vee x \vee \neg y \qquad\qquad \text{by } x \le y \text{ and \textit{antitone} property (Definition B.3 page 35)}$$

$$= y \vee \neg y \qquad\qquad\qquad\qquad \text{by \textit{orthomodular identity} (Definition B.3 page 35)}$$

$$= 1 \qquad\qquad\qquad\qquad\qquad \text{by \textit{excluded middle} property of } \neg \text{ (Theorem B.15 page 37)}$$

(c) Proof that on an *orthocomplemented lattice*, $\overset{r}{\to}$ *does* have *weak modus ponens*:

$$x \wedge (x \overset{r}{\to} y) \triangleq x \wedge [(\neg x \wedge y) \vee (x \wedge y) \vee (\neg x \wedge \neg y)] \quad \text{by definition of } \overset{r}{\to}$$

$$\le [(\neg x \wedge y) \vee (x \wedge y) \vee (\neg x \wedge \neg y)] \qquad \text{by definition of } \wedge \text{ (Definition A.9 page 24)}$$

$$\le \neg x \vee (x \wedge y) \vee (\neg x \wedge \neg y) \qquad\qquad \text{by definition of } \wedge \text{ (Definition A.9 page 24)}$$

$$\le \neg x \vee y \vee (\neg x \wedge \neg y) \qquad\qquad\qquad \text{by definition of } \wedge \text{ (Definition A.9 page 24)}$$

$$\le \neg x \vee y \qquad\qquad\qquad\qquad\qquad \text{by \textit{absorption} property (Theorem A.14 page 25)}$$

$$\implies \textit{weak modus ponens}$$

(d)　Proof that on a *Boolean lattice*, $\overset{r}{\hookrightarrow}=\overset{s}{\hookrightarrow}$:

$$x \overset{r}{\hookrightarrow} y \triangleq (\neg x \wedge y) \vee (x \wedge y) \vee (\neg x \wedge \neg y) \qquad \text{by definition of } \overset{r}{\hookrightarrow}$$

$$= [\neg x \wedge (y \vee \neg y)] \vee (x \wedge y) \qquad \text{by } \textit{distributive} \text{ property (Definition A.41 page 30)}$$

$$= [\neg x \wedge 1] \vee (x \wedge y) \qquad \text{by } \textit{excluded middle} \text{ property of } \neg \text{ (Theorem B.15 page 37)}$$

$$= \neg x \vee (x \wedge y) \qquad \text{by definition of 1 and } \wedge \text{ (Definition A.9 page 24)}$$

$$= \neg x \vee y \qquad \text{by property of } \textit{Boolean lattice}\text{s (Lemma A.43 page 30)}$$

$$\triangleq x \overset{s}{\hookrightarrow} y \qquad \text{by definition of } \overset{s}{\hookrightarrow}$$

## C.2　Logics

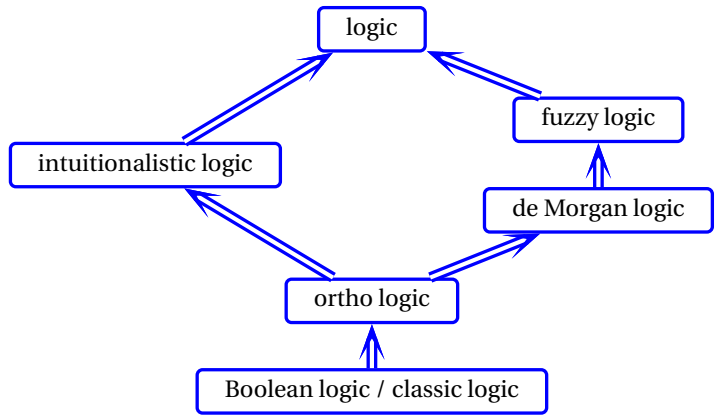

Figure 13: lattice of logics

**Definition C.5**　[87] Let $\rightarrow$ be an *implication* (Definition C.1 page 45) defined on a *lattice with negation* $\boldsymbol{L} \triangleq$ ( $X$, $\vee$, $\wedge$, $\neg$, 0, 1 ; $\leq$) (Definition B.5 page 35).

| | |
|---|---|
| ( $X$, $\vee$, $\wedge$, $\neg$, 0, 1 ; $\leq$, $\rightarrow$) is a **logic** | if $\neg$ is a *minimal negation*. |
| ( $X$, $\vee$, $\wedge$, $\neg$, 0, 1 ; $\leq$, $\rightarrow$) is a **fuzzy logic** | if $\neg$ is a *fuzzy negation*. |
| ( $X$, $\vee$, $\wedge$, $\neg$, 0, 1 ; $\leq$, $\rightarrow$) is an **intuitionalistic logic** | if $\neg$ is an *intuitionalistic negation*. |
| ( $X$, $\vee$, $\wedge$, $\neg$, 0, 1 ; $\leq$, $\rightarrow$) is a **de Morgan logic** | if $\neg$ is a *de Morgan negation*. |
| ( $X$, $\vee$, $\wedge$, $\neg$, 0, 1 ; $\leq$, $\rightarrow$) is a **Kleene logic** | if $\neg$ is a *Kleene negation*. |
| ( $X$, $\vee$, $\wedge$, $\neg$, 0, 1 ; $\leq$, $\rightarrow$) is an **ortho logic** | if $\neg$ is an *ortho negation*. |
| ( $X$, $\vee$, $\wedge$, $\neg$, 0, 1 ; $\leq$, $\rightarrow$) is a **Boolean logic** | if $\neg$ is an *ortho negation* and $\boldsymbol{L}$ is *Boolean*. |

**Example C.6**　(Aristotelian logic / classical logic)　[88] The *classical bi-variate logic* is defined below. It is a 2 element *Boolean logic* (Definition C.5 page 50). with $\boldsymbol{L} \triangleq$ ( $\{1, 0\}$, $\wedge$, $\neg$, 0, 1, $\leq$ ; $\vee$) and a *classical implication* $\rightarrow$ with *strong entailment* and *strong modus ponens*. The value 1 represents "*true*" and 0 represents

---

[87]　📖　[Straßburger(2005)] page 136 ⟨Definition 2.1⟩, 📖 [de Vries(2007)] page 11 ⟨Definition 16⟩
[88]　📄　[Novák et al.(1999)Novák, Perfilieva, and Močkoř] pages 17–18 ⟨Example 2.1⟩

"*false*".

$$x \to y \triangleq \left\{ \begin{array}{ll} 1 & \forall x \le y \\ y & \text{otherwise} \end{array} \right\} = \left\{ \begin{array}{c|cc} \to & 1 & 0 \\ \hline 1 & 1 & 0 \\ 0 & 1 & 1 \end{array} \quad \forall x,y \in X \right\} = \neg x \lor y$$

✎ PROOF:

(1) Proof that $\neg$ is an *ortho negation*: by Definition B.3 (page 35)

(2) Proof that $\to$ is an *implication* with *strong entailment* and *strong modus ponens*:

    (a) **L** is *Boolean* and therefore is *orthocomplemented*.

    (b) $\to$ is equivalent to the *classical implication* $\overset{c}{\to}$ (Example C.4 page 46).

    (c) By Example C.4 (page 46), $\to$ has *strong entailment* and *strong modus ponens*.

The *classical logic* (previous example) can be generalized in several ways. Arguably one of the simplest of these is the 3-valued logic due to Kleene (next example).

**Example C.7** [89] The **Kleene 3-valued logic** ( $X$, $\lor$, $\land$, $\neg$, $0$, $1$ ; $\le$, $\to$ ) is defined below. The function $\neg$ is a *Kleene negation* (Definition B.3 page 35) and is presented in Example B.25 (page 41). The function $\to$ is the *classic implication* $x \to y \triangleq \neg x \lor y$. The values 1 represents "*true*", 0 represents "*false*", and $n$ represents "*neutral*" or "*undecided*".

$$x \to y \triangleq \left\{ \neg x \lor y \quad \forall x \in X \right\} = \left\{ \begin{array}{c|ccc} \to & 1 & n & 0 \\ \hline 1 & 1 & n & 0 \\ n & 1 & n & n \\ 0 & 1 & 1 & 1 \end{array} \quad \forall x,y \in X \right\}$$

✎ PROOF:

(1) Proof that $\neg$ is a *Kleene negation*: see Example B.25 (page 41)

(2) Proof that $\to$ is an *implication*: This follows directly from the definition of $\to$ and the definition of an *implication* (Definition C.1 page 45).

(3) Proof that $\to$ does not have *strong entailment*: $n \to n = n = n \lor n \ne 1$.

(4) Proof that $\to$ does not have *strong modus ponens*: $n \to 0 = n = \neg n \lor 0 \not\le 0$.

A lattice and negation alone do not uniquely define a logic. Łukasiewicz also introduced a 3-valued logic with identical lattice structure to Kleene, but with a different implication relation (next example). Historically, Łukasiewicz's logic was introduced before Kleene's.

---

[89] 📖 [Kleene(1938)] page 153, ✎ [Kleene(1952)], pages 332–339 ⟨§64. The 3-valued logic⟩, 📖 [Avron(1991)] page 277

**Example C.8** [90]

The **Łukasiewicz 3-valued logic** ( $X$, $\vee$, $\wedge$, $\neg$, 0, 1 ; $\leq$, $\rightarrow$ ) is defined to the right and below. The function $\neg$ is a *Kleene negation* (Definition B.3 page 35) and is presented in Example B.25 (page 41). The implication has *strong entailment* but *weak modus ponens*. In the implication table below, values that differ from the classical $x \rightarrow y \triangleq \neg x \vee y$ are shaded.

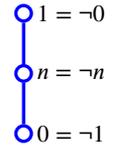

$$x \rightarrow y \triangleq \left\{ \begin{array}{ll} 1 & \forall x \leq y \\ \neg x \vee y & \text{otherwise} \end{array} \right\} = \left\{ \begin{array}{c|ccc} \rightarrow & 1 & n & 0 \\ \hline 1 & 1 & n & 0 \\ n & 1 & \boxed{1} & n \\ 0 & 1 & 1 & 1 \end{array} \right\} \; \forall x,y \in X = \left\{ \begin{array}{ll} 1 & \text{for } x = y = n \\ \neg x \vee y & \text{otherwise} \end{array} \right\}$$

✎Proof:

(1) Proof that $\neg$ is a *Kleene negation*: see Example B.25 (page 41)

(2) Proof that $\rightarrow$ is an *implication*: This follows directly from the definition of $\rightarrow$ and the definition of an *implication* (Definition C.1 page 45).

(3) Proof that $\rightarrow$ does not have *strong modus ponens*: $n \rightarrow 0 = n = \neg n \vee 0 \not\leq 0$.

✏

**Example C.9** [91] The **RM₃ logic** ( $X$, $\vee$, $\wedge$, $\neg$, 0, 1 ; $\leq$, $\rightarrow$ ) is defined below. The function $\neg$ is a *Kleene negation* (Definition B.3 page 35) and is presented in Example B.25 (page 41). The implication function has *weak entailment* but *strong modus ponens*. In the implication table below, values that differ from the classical $x \rightarrow y \triangleq \neg x \vee y$ are shaded.

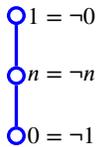

$$x \rightarrow y \triangleq \left\{ \begin{array}{ll} 1 & \forall x < y \\ n & \forall x = y \\ 0 & \forall x > y \end{array} \right\} = \left\{ \begin{array}{c|ccc} \rightarrow & 1 & n & 0 \\ \hline 1 & 1 & n & 0 \\ n & 1 & n & 0 \\ 0 & 1 & 1 & 1 \end{array} \right\} \; \forall x,y \in X$$

✎Proof:

(1) Proof that $\neg$ is a *Kleene negation*: see Example B.25 (page 41)

(2) Proof that $\rightarrow$ is an *implication*: This follows directly from the definition of $\rightarrow$ and the definition of an *implication* (Definition C.1 page 45).

(3) Proof that $\rightarrow$ does not have *strong entailment*: $n \rightarrow n = n = n \vee n \neq 1$.

✏

In a 3-valued logic, the negation does not necessarily have to be as in the previous three examples. The next example offers a different negation.

**Example C.10** (Heyting 3-valued logic/Jaśkowski's first matrix) [92]

The **Heyting 3-valued logic** ( $X$, $\vee$, $\wedge$, $\neg$, 0, 1 ; $\leq$, $\rightarrow$ ) is defined below. The negation $\neg$ is both *intuitionistic* and *fuzzy* (Definition B.2 page 35), and is defined on a 3 element *linearly ordered lattice* (Definition A.3 page 23).

The implication function has both *strong entailment* and *strong modus ponens*. In the implication table below, values that differ from the classical $x \to y \triangleq \neg x \vee y$ are shaded.

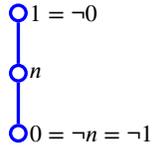

$$x \to y \triangleq \left\{ \begin{array}{ll} 1 & \forall x \leq y \\ y & \text{otherwise} \end{array} \right\} = \left\{ \begin{array}{c|ccc} \to & 1 & n & 0 \\ \hline 1 & 1 & n & 0 \\ n & 1 & 1 & 0 \\ 0 & 1 & 1 & 1 \end{array} \quad \forall x,y \in X \right\}$$

✎ Proof:

(1) Proof that ¬ is a *Kleene negation*: see Example B.26 (page 41)

(2) Proof that → is an *implication*: by definition of *implication* (Definition C.1 page 45)

Of course it is possible to generalize to more than 3 values (next example).

**Example C.11** [93] The **Łukasiewicz 5-valued logic** ( $X$, $\vee$, $\wedge$, $\neg$, $0$, $1$ ; $\leq$, $\to$ ) is defined below. The implication function has *strong entailment* but *weak modus ponens*. In the implication table below, values that differ from the classical $x \to y \triangleq \neg x \vee y$ are shaded.

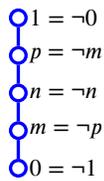

$$x \to y \triangleq \begin{array}{c|ccccc} \to & 1 & p & n & m & 0 \\ \hline 1 & 1 & p & n & m & 0 \\ p & 1 & 1 & n & m & m \\ n & 1 & 1 & 1 & m & n \\ m & 1 & 1 & 1 & 1 & p \\ 0 & 1 & 1 & 1 & 1 & 1 \end{array} \quad \forall x,y \in X$$

✎ Proof:

All the previous examples in this section are *linearly ordered*. The following examples employ logics that are not.

**Example C.12** [94] The **Boolean 4-valued logic** is defined below. The negation function ¬ is an *ortho negation* (Example B.28 page 42) defined on an $M_2$ lattice. The value 1 represents "*true*", 0 represents "*false*", and $m$ and $n$ represent some intermediate values.

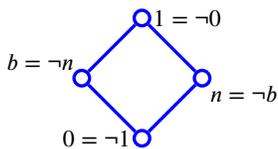

$$x \to y \triangleq \neg x \vee y = \begin{array}{c|cccc} \to & 1 & b & n & 0 \\ \hline 1 & 1 & b & n & 0 \\ b & 1 & 1 & n & n \\ n & 1 & b & 1 & b \\ 0 & 1 & 1 & 1 & 1 \end{array} \quad \forall x,y \in X$$

All the previous examples in this section are *distributive*; the previous example was *Boolean*. The next example is *non-distributive*, and *de Morgan* (but *non-Boolean*). Note for a given order structure, the method of negation may not be unique; in the previous and following examples both have identical lattices, but are negated differently.

**Example C.13**   [95] The **BN₄ logic** is defined below. The function ¬ is a *de Morgan negation* (Example B.29 page 43) defined on a 4 element $M_2$ *lattice*. The value 1 represents "*true*", 0 represents "*false*", $b$ represents "*both*" (both true and false), and $n$ represents "*neither*". In the implication table below, the values that differ from those of the *classical implication* $\overset{c}{\to}$ are shaded.

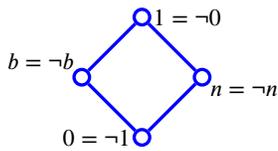

$$x \to y \triangleq \begin{array}{c|cccc} \to & 1 & n & b & 0 \\ \hline 1 & 1 & n & 0 & 0 \\ n & 1 & 1 & n & n \\ b & 1 & n & b & 0 \\ 0 & 1 & 1 & 1 & 1 \end{array} \qquad \forall x,y \in X$$

**Example C.14**

    The tables that follow are the 6 implications defined in Example C.4 (page 46) on the $O_6$ *lattice with ortho negation* (Definition B.3 page 35), or the $O_6$ *orthocomplemented lattice* (Definition A.45 page 31), illustrated to the right. In the tables, the values that differ from those of the *classical implication* $\overset{c}{\to}$ are shaded.

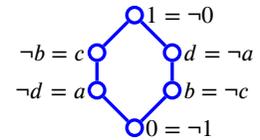

| $\overset{c}{\to}$ | 1 | d | c | b | a | 0 |
|---|---|---|---|---|---|---|
| 1 | 1 | d | c | b | a | 0 |
| d | 1 | 1 | c | 1 | a | a |
| c | 1 | d | 1 | b | 1 | b |
| b | 1 | 1 | c | 1 | c | c |
| a | 1 | d | 1 | d | 1 | d |
| 0 | 1 | 1 | 1 | 1 | 1 | 1 |

| $\overset{s}{\to}$ | 1 | d | c | b | a | 0 |
|---|---|---|---|---|---|---|
| 1 | 1 | d | c | b | a | 0 |
| d | 1 | 1 | a | 1 | a | a |
| c | 1 | b | 1 | b | 1 | b |
| b | 1 | 1 | c | 1 | c | c |
| a | 1 | 1 | 1 | d | 1 | d |
| 0 | 1 | 1 | 1 | 1 | 1 | 1 |

| $\overset{d}{\to}$ | 1 | d | c | b | a | 0 |
|---|---|---|---|---|---|---|
| 1 | 1 | d | c | b | a | 0 |
| d | 1 | 1 | c | 1 | a | a |
| c | 1 | d | 1 | b | 1 | b |
| b | 1 | 1 | c | 1 | a | c |
| a | 1 | d | 1 | b | 1 | d |
| 0 | 1 | 1 | 1 | 1 | 1 | 1 |

| $\overset{k}{\to}$ | 1 | d | c | b | a | 0 |
|---|---|---|---|---|---|---|
| 1 | 1 | d | c | b | a | 0 |
| d | 1 | 1 | a | 1 | a | a |
| c | 1 | b | 1 | b | 1 | b |
| b | 1 | 1 | c | 1 | a | c |
| a | 1 | d | 1 | b | 1 | d |
| 0 | 1 | 1 | 1 | 1 | 1 | 1 |

| $\overset{n}{\to}$ | 1 | d | c | b | a | 0 |
|---|---|---|---|---|---|---|
| 1 | 1 | d | c | b | a | 0 |
| d | 1 | 1 | a | 1 | a | a |
| c | 1 | b | 1 | b | 1 | b |
| b | 1 | 1 | c | 1 | a | c |
| a | 1 | d | 1 | b | 1 | d |
| 0 | 1 | 1 | 1 | 1 | 1 | 1 |

| $\overset{r}{\to}$ | 1 | d | c | b | a | 0 |
|---|---|---|---|---|---|---|
| 1 | 1 | d | c | b | a | 0 |
| d | 1 | 1 | a | 1 | a | a |
| c | 1 | b | 1 | b | 1 | b |
| b | 1 | 1 | c | 1 | a | c |
| a | 1 | d | 1 | b | 1 | d |
| 0 | 1 | 1 | 1 | 1 | 1 | 1 |

**Example C.15**   [96] A 6 element logic is defined below. The function ¬ is a *Kleene negation* (Example B.32 page 44). The implication has *strong entailment* but *weak modus ponens*. In the implication table below, the values that differ from those of the *classical implication* $\overset{c}{\to}$ are shaded.

---

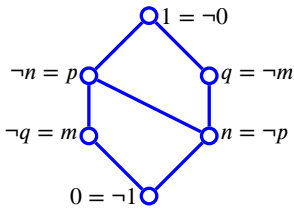

| $\rightarrow$ | 1 | $p$ | $q$ | $m$ | $n$ | 0 |
|---|---|---|---|---|---|---|
| 1 | 1 | $p$ | $q$ | $m$ | $n$ | 0 |
| $p$ | 1 | **1** | $q$ | $p$ | **$q$** | $n$ |
| $q$ | 1 | $p$ | 1 | $m$ | $p$ | $m$ |
| $m$ | 1 | 1 | $q$ | 1 | $q$ | $q$ |
| $n$ | 1 | **1** | 1 | $p$ | **1** | $p$ |
| 0 | 1 | 1 | 1 | 1 | 1 | 1 |

$$x \rightarrow y \triangleq \qquad \forall x,y \in X$$

✎Proof:

(1) Proof that ¬ is a *Kleene negation*: see Example B.32 (page 44)

(2) Proof that → is an *implication*: This follows directly from the definition of → and the definition of an *implication* (Definition C.1 page 45).

(3) Proof that → does not have *strong modus ponens*:

$$\neg p \wedge (p \rightarrow m) \; = \; n \wedge p \; = \; n \; \leq \; p \; = \; \neg p \vee m \; \nleq \; m$$
$$\neg n \wedge (n \rightarrow m) \; = \; n \wedge p \; = \; n \; \leq \; p \; = \; \neg p \vee m \; \nleq \; m$$
$$\neg p \wedge (p \rightarrow 0) \; = \; n \wedge n \; = \; n \; \leq \; n \; = \; \neg p \vee 0 \; \nleq \; 0$$
$$\neg n \wedge (n \rightarrow 0) \; = \; p \wedge n \; = \; n \; \leq \; p \; = \; \neg n \vee 0 \; \nleq \; 0$$

For an example of an 8-valued logic, see  [Kamide(2013)]. For examples of 16-valued logics, see 📖 [Shramko and Wansing(2005)].

# Reference Index

# Subject Index

*Telecommunications Engineering Department, National Chiao-Tung University, Hsinchu, Taiwan;* 國立交通大學 *(Gúo Lì Jiāo Tōng Dà Xúe)* 電信工程學系 *(Diàn Xìn Gōng Chéng Xúe Xì)* 新竹，台灣 *(Xīn Zhú, Tái Wān)*

[dgreenhoe@gmail.com](mailto:dgreenhoe@gmail.com)